\documentclass[a4paper,reqno,12pt]{amsart}
\usepackage{mathrsfs}
\usepackage{stmaryrd}
\usepackage{cases}
\usepackage{epic}
\usepackage{amsfonts}
\usepackage{graphicx}
\usepackage{fancyhdr}
\usepackage{amsmath}
\usepackage{amssymb, upgreek}
\usepackage{mathtools}
\usepackage{bm}
\usepackage{array}
\usepackage{chngcntr} 
\counterwithin{figure}{section}
\counterwithin{table}{section}
\usepackage{float}
\usepackage{latexsym,todonotes}
\usepackage{pdflscape}
\usepackage[all]{xypic}
\usepackage[all]{xy}
\usepackage{color}
\usepackage{colordvi}
\usepackage{multicol}
\usepackage{tikz}
\usepackage{tikz-cd}
\usepackage[linktocpage=true]{hyperref}
\hypersetup{colorlinks,linkcolor=blue,urlcolor=cyan,citecolor=blue}
\newcommand{\bvs}{\mathbf{\varsigma}}
\newcommand{\vs}{\varsigma}
\allowdisplaybreaks

\DeclareMathAlphabet{\mathbbb}{U}{bbold}{m}{n}

\begin{document}
	\input xy
	\xyoption{all}

	\newtheorem{innercustomthm}{{\bf Main Theorem}}
	\newenvironment{customthm}[1]
	{\renewcommand\theinnercustomthm{#1}\innercustomthm}
	{\endinnercustomthm}
	
	\newtheorem{innercustomcor}{{\bf Corollary}}
	\newenvironment{customcor}[1]
	{\renewcommand\theinnercustomcor{#1}\innercustomcor}
	{\endinnercustomthm}
	
	\newtheorem{innercustomprop}{{\bf Proposition}}
	\newenvironment{customprop}[1]
	{\renewcommand\theinnercustomprop{#1}\innercustomprop}
	{\endinnercustomthm}
	\def \wt{\mathrm{wt}}
	\def \wItau{\I_{\circ,\tau}}
	\def \brW{\mathrm{Br}(W_\btau)}
	\def \br{\mathbf{r}}
	\newcommand{\iadd}{\operatorname{iadd}\nolimits}
	\newcommand{\Gr}{\operatorname{Gr}\nolimits}
	\newcommand{\FGS}{\operatorname{FGS}\nolimits}
	\def \tbU{\tU_{\bullet}}
	\def \bU{\U_{\bullet}}
	\renewcommand{\mod}{\operatorname{mod^{\rm nil}}\nolimits}
	\newcommand{\proj}{\operatorname{proj}\nolimits}
	\newcommand{\inj}{\operatorname{inj.}\nolimits}
	\newcommand{\rad}{\operatorname{rad}\nolimits}
	\newcommand{\Span}{\operatorname{Span}\nolimits}
	\newcommand{\soc}{\operatorname{soc}\nolimits}
	\newcommand{\ind}{\operatorname{inj.dim}\nolimits}
	\newcommand{\Ginj}{\operatorname{Ginj}\nolimits}
	\newcommand{\res}{\operatorname{res}\nolimits}
	\newcommand{\np}{\operatorname{np}\nolimits}
	\newcommand{\Fac}{\operatorname{Fac}\nolimits}
	\newcommand{\Aut}{\operatorname{Aut}\nolimits}
	\newcommand{\DTr}{\operatorname{DTr}\nolimits}
	\newcommand{\TrD}{\operatorname{TrD}\nolimits}
	\newcommand{\Mod}{\operatorname{Mod}\nolimits}
	\newcommand{\R}{\operatorname{R}\nolimits}
	\newcommand{\End}{\operatorname{End}\nolimits}
	\newcommand{\lf}{\operatorname{l.f.}\nolimits}
	\newcommand{\Iso}{\operatorname{Iso}\nolimits}
	\newcommand{\aut}{\operatorname{Aut}\nolimits}
	\newcommand{\Ui}{{\mathbf U}^\imath}	\newcommand{\UU}{{\mathbf U}\otimes {\mathbf U}}
	\newcommand{\UUi}{(\UU)^\imath}
	\newcommand{\tUU}{{\tU}\otimes {\tU}}
	\newcommand{\tUUi}{(\tUU)^\imath}
	\newcommand{\tUi}{\widetilde{{\mathbf U}}^\imath}
	\newcommand{\sqq}{{\bf v}}
	\newcommand{\sqvs}{\sqrt{\vs}}
	\newcommand{\dbl}{\operatorname{dbl}\nolimits}
	\newcommand{\swa}{\operatorname{swap}\nolimits}
	\newcommand{\Gp}{\operatorname{Gp}\nolimits}
	\def \tTD{\widetilde{T}}
	\def \fX{\Upsilon}
	\newcommand{\U}{{\mathbf U}}
	\newcommand{\tU}{\widetilde{\mathbf U}}
	\newcommand{\fgm}{{\rm mod}^{{\rm fg}}}
	\newcommand{\fgmz}{\mathrm{mod}^{{\rm fg},\Z}}
	\newcommand{\fdmz}{\mathrm{mod}^{{\rm nil},\Z}}
	
	\newcommand{\ov}{\overline}
	\newcommand{\und}{\underline}
	\newcommand{\tk}{\widetilde{k}}
	\newcommand{\tK}{K}
	\newcommand{\tTT}{\operatorname{\widetilde{\mathbf{T}}}\nolimits}
	
	\def \hU{\widehat{\U}}
	\def \hUi{\widehat{\U}^\imath}
	\def \bfk{\mathbf{k}}
	\def \wI{\I_\circ}
	\def \cR{\mathcal{R}}
	\newcommand{\tH}{\operatorname{{\ch}_{\rm{tw}}}\nolimits}
	
	\newcommand{\utM}{\operatorname{\cm\ch}\nolimits}
	\newcommand{\tM}{\operatorname{\cs\cd\widetilde{\ch}}\nolimits}
	\newcommand{\rM}{\operatorname{\cm\ch_{\rm{red}}}\nolimits}
	\newcommand{\utMH}{\cs\cd\ch(\Lambda^\imath)}
	\newcommand{\tMH}{\cs\cd\widetilde{\ch}(\Lambda^\imath)}
	\newcommand{\tCMH}{{\cc\widetilde{\ch}(\K Q,\btau)}}
	
	\newcommand{\rMH}{\operatorname{\cs\cd\ch_{\rm{red}}(\Lambda^\imath)}\nolimits}
	\newcommand{\utMHg}{\operatorname{\ch(Q,\btau)}\nolimits}
	\newcommand{\tMHg}{\operatorname{\widetilde{\ch}(Q,\btau)}\nolimits}
	\newcommand{\tMHk}{{\widetilde{\ch}(\K Q,\btau)}}
	\newcommand{\rMHg}{\operatorname{\ch_{\rm{red}}(Q,\btau)}\nolimits}
	
	\newcommand{\rMHd}{\operatorname{\cm\ch_{\rm{red}}(\Lambda^\imath)_{\bvsd}}\nolimits}
	\newcommand{\tMHd}{\operatorname{\cs\cd\widetilde{\ch}(\Lambda^\imath)_{\bvsd}}\nolimits}
	
	\def \h{h} 
	\newcommand{\tMHl}{\cs\cd\widetilde{\ch}({\bs}_\ell\Lambda^\imath)}
	\newcommand{\rMHl}{\cm\ch_{\rm{red}}({\bs}_\ell\Lambda^\imath)_{\bvsd}}
	\newcommand{\tMHi}{\cs\cd\widetilde{\ch}({\bs}_i\Lambda^\imath)}
	\newcommand{\rMHi}{\cm\ch_{\rm{red}}({\bs}_i\Lambda^\imath)_{\bvsd}}
	\newcommand{\tMHgi}{\widetilde{\ch}({\bs}_i Q,\btau)}
	\def \bI{\I_\bullet}
	\newcommand{\utGpg}{\operatorname{\ch^{\rm Gp}(Q,\btau)}\nolimits}
	\newcommand{\tGpg}{\operatorname{\widetilde{\ch}^{\rm Gp}(Q,\btau)}\nolimits}
	\newcommand{\rGpg}{\operatorname{\ch_{red}^{\rm Gp}(Q,\btau)}\nolimits}

	\newcommand{\colim}{\operatorname{colim}\nolimits}
	\newcommand{\gldim}{\operatorname{gl.dim}\nolimits}
	\newcommand{\cone}{\operatorname{cone}\nolimits}
	\newcommand{\rep}{\operatorname{rep}\nolimits}
	\newcommand{\Ext}{\operatorname{Ext}\nolimits}
	\newcommand{\Tor}{\operatorname{Tor}\nolimits}
	\newcommand{\Hom}{\operatorname{Hom}\nolimits}
	\newcommand{\Top}{\operatorname{top}\nolimits}
	\newcommand{\Coker}{\operatorname{Coker}\nolimits}
	\newcommand{\thick}{\operatorname{thick}\nolimits}
	\newcommand{\rank}{\operatorname{rank}\nolimits}
	\newcommand{\Gproj}{\operatorname{Gproj}\nolimits}
	\newcommand{\Len}{\operatorname{Length}\nolimits}
	\newcommand{\RHom}{\operatorname{RHom}\nolimits}
	\renewcommand{\deg}{\operatorname{deg}\nolimits}
	\renewcommand{\Im}{\operatorname{Im}\nolimits}
	\newcommand{\Ker}{\operatorname{Ker}\nolimits}
	\newcommand{\Coh}{\operatorname{Coh}\nolimits}
	\newcommand{\Id}{\operatorname{Id}\nolimits}
	\newcommand{\Qcoh}{\operatorname{Qch}\nolimits}
	\newcommand{\CM}{\operatorname{CM}\nolimits}
	\newcommand{\sgn}{\operatorname{sgn}\nolimits}
	\newcommand{\Gdim}{\operatorname{G.dim}\nolimits}
	\newcommand{\fpr}{\operatorname{\mathcal{P}^{\leq1}}\nolimits}
	
	\newcommand{\For}{\operatorname{{\bf F}or}\nolimits}
	\newcommand{\coker}{\operatorname{Coker}\nolimits}
	\renewcommand{\dim}{\operatorname{dim}\nolimits}
	\newcommand{\rankv}{\operatorname{\underline{rank}}\nolimits}
	\newcommand{\dimv}{{\operatorname{\underline{dim}}\nolimits}}
	\newcommand{\diag}{{\operatorname{diag}\nolimits}}
	\newcommand{\qbinom}[2]{\begin{bmatrix} #1\\#2 \end{bmatrix} }
	
	\renewcommand{\Vec}{{\operatorname{Vec}\nolimits}}
	\newcommand{\pd}{\operatorname{proj.dim}\nolimits}
	\newcommand{\gr}{\operatorname{gr}\nolimits}
	\newcommand{\id}{\operatorname{Id}\nolimits}
	\newcommand{\Res}{\operatorname{Res}\nolimits}
	\def \tT{\widetilde{\mathscr{T} }}
	
	\def \bwi{w_{\bullet,i}}
	\def \tTL{\tT(\Lambda^\imath)}
	
	\newcommand{\mbf}{\mathbf}
	\newcommand{\mbb}{\mathbb}
	\newcommand{\mrm}{\mathrm}
	\newcommand{\cbinom}[2]{\left\{ \begin{matrix} #1\\#2 \end{matrix} \right\}}
	\newcommand{\dvev}[1]{{B_1|}_{\ev}^{{(#1)}}}
	\newcommand{\dv}[1]{{B_1|}_{\odd}^{{(#1)}}}
	\newcommand{\dvd}[1]{t_{\odd}^{{(#1)}}}
	\newcommand{\dvp}[1]{t_{\ev}^{{(#1)}}}
	\newcommand{\ev}{\bar{0}}
	\newcommand{\odd}{\bar{1}}
	\newcommand{\Iblack}{\I_{\bullet}}
	\newcommand{\wb}{w_\bullet}
	\newcommand{\Uidot}{\dot{\bold{U}}^{\imath}}
	
	\newcommand{\kk}{h}
	\newcommand{\la}{\lambda}
	\newcommand{\LR}[2]{\left\llbracket \begin{matrix} #1\\#2 \end{matrix} \right\rrbracket}
	\newcommand{\ff}{\mathbf{f}}
	\newcommand{\pdim}{\operatorname{proj.dim}\nolimits}
	\newcommand{\idim}{\operatorname{inj.dim}\nolimits}
	\newcommand{\Gd}{\operatorname{G.dim}\nolimits}
	\newcommand{\Ind}{\operatorname{Ind}\nolimits}
	\newcommand{\add}{\operatorname{add}\nolimits}
	\newcommand{\ad}{\operatorname{ad}\nolimits}
	\newcommand{\pr}{\operatorname{pr}\nolimits}
	\newcommand{\oR}{\operatorname{R}\nolimits}
	\newcommand{\oL}{\operatorname{L}\nolimits}
	\newcommand{\ext}{{ \mathfrak{Ext}}}
	\newcommand{\Perf}{{\mathfrak Perf}}
	\def\scrP{\mathscr{P}}
	\newcommand{\bk}{{\mathbb K}}
	\newcommand{\cc}{{\mathcal C}}
	\newcommand{\gc}{{\mathcal GC}}
	\newcommand{\dg}{{\rm dg}}
	\newcommand{\ce}{{\mathcal E}}
	\newcommand{\cs}{{\mathcal S}}
	\newcommand{\cl}{{\mathcal L}}
	\newcommand{\cf}{{\mathcal F}}
	\newcommand{\cx}{{\mathcal X}}
	\newcommand{\cy}{{\mathcal Y}}
	\newcommand{\ct}{{\mathcal T}}
	\newcommand{\cu}{{\mathcal U}}
	\newcommand{\cv}{{\mathcal V}}
	\newcommand{\cn}{{\mathcal N}}
	\newcommand{\mcr}{{\mathcal R}}
	\newcommand{\ch}{{\mathcal H}}
	\newcommand{\ca}{{\mathcal A}}
	\newcommand{\cb}{{\mathcal B}}
	\newcommand{\ci}{{\I}_{\btau}}
	\newcommand{\cj}{{\mathcal J}}
	\newcommand{\cm}{{\mathcal M}}
	\newcommand{\cp}{{\mathcal P}}
	\newcommand{\cg}{{\mathcal G}}
	\newcommand{\cw}{{\mathcal W}}
	\newcommand{\co}{{\mathcal O}}
	\newcommand{\cq}{{Q^{\rm dbl}}}
	\newcommand{\cd}{{\mathcal D}}
	\newcommand{\ck}{\widetilde{\mathcal K}}
	\newcommand{\calr}{{\mathcal R}}
	\newcommand{\iLa}{\Lambda^{\imath}}
	\newcommand{\La}{\Lambda}
	\newcommand{\ol}{\overline}
	\newcommand{\ul}{\underline}
	\newcommand{\st}{[1]}
	\newcommand{\ow}{\widetilde}
	\renewcommand{\P}{\mathbf{P}}
	\newcommand{\pic}{\operatorname{Pic}\nolimits}
	\newcommand{\Spec}{\operatorname{Spec}\nolimits}
	
	\newtheorem{theorem}{Theorem}[section]
	\newtheorem{acknowledgement}[theorem]{Acknowledgement}
	\newtheorem{conjecture}[theorem]{Conjecture}
	\newtheorem{corollary}[theorem]{Corollary}
	\newtheorem{definition}[theorem]{Definition}
	\newtheorem{example}[theorem]{Example}
	\newtheorem{lemma}[theorem]{Lemma}
	\newtheorem{notation}[theorem]{Notation}
	\newtheorem{problem}[theorem]{Problem}
	\newtheorem{proposition}[theorem]{Proposition}
	\newtheorem{summary}[theorem]{Summary}
	\numberwithin{equation}{section}

	\newtheorem{alphatheorem}{Theorem}
	\newtheorem{alphacorollary}[alphatheorem]{Corollary}
	\newtheorem{alphaproposition}[alphatheorem]{Proposition}
	\renewcommand*{\thealphatheorem}{\Alph{alphatheorem}}
	
	\theoremstyle{remark}
	\newtheorem{remark}[theorem]{Remark}

	\newcommand{\Pd}{\pi_*}
	\def \bvs{{\boldsymbol{\varsigma}}}
	\def \bvsd{{\boldsymbol{\varsigma}_{\diamond}}}
	\def \btau{{{\tau}}}
	
	\def \Br{\mathrm{Br}}
	\def \bp{{\mathbf p}}
	\def \bq{{\bm q}}
	\def \bv{{v}}
	\def \bs{{r}}
	
	\def \bfK{{\mathbf K}}
	
	\def \bw{w_\bullet}
	
	\newcommand{\tCMHg}{\cc\widetilde{\ch}(Q,\btau)}
	\newcommand{\bfv}{\mathbf{v}}
	\def \bA{{\mathbf A}}
	\def \ba{{\mathbf a}}
	\def \bi{{\mathbf i}}
	\def \bL{{\mathbf L}}
	\def \bF{{\mathbf F}}
	\def \bS{{\mathbf S}}
	\def \bC{{\mathbf C}}
	\def \bU{{\mathbf U}}
	\def \bc{{\mathbf c}}
	\def \fpi{\mathfrak{P}^\imath}
	\def \Ni{N^\imath}
	\def \fp{\mathfrak{P}}
	\def \fg{\mathfrak{g}}
	\def \fk{\fg^\theta}  
	
	\def\Inv{\mathrm{Inv}}
	\def \bbW{W^{\circ}}
	\def \bbw{{\boldsymbol{w}}}
	\def \BB{\mathbf{B}}
	\def \tB{\widetilde{\mathbf{B}}}
	\def \hB{\widehat{\mathbf{B}}}
	\def \reW{W^\circ} 
	\def \fn{\mathfrak{n}}
	\def \fh{\mathfrak{h}}
	\def \fu{\mathfrak{u}}
	\def \fv{\mathfrak{v}}
	\def \fa{\mathfrak{a}}
	\def \fq{\mathfrak{q}}
	\def \Z{{\Bbb Z}}
	\def \F{{\Bbb F}}
	\def \D{{\Bbb D}}
	\def \C{{\Bbb C}}
	\def \N{{\Bbb N}}
	\def \Q{{\Bbb Q}}
	\def \G{{\Bbb G}}
	\def \P{{\Bbb P}}
	\def \K{{\mathbb K}}
	\def \bK{{\Bbb K}}
	\def \J{{\Bbb J}}
	\def \E{{\Bbb E}}
	\def \A{{\Bbb A}}
	\def \L{{\Bbb L}}
	\def \I{{\Bbb I}}
	\def \BH{{\Bbb H}}
	\def \T{{\Bbb T}}
	\def \upLa {\raisebox{8pt}{\rotatebox{180}{$\Lambda$}}}
	\newcommand{\TT}
	{\operatorname{\mathbf{T}}\nolimits}
	\newcommand {\lu}[1]{\textcolor{red}{$\clubsuit$: #1}}
	
	\newcommand{\nc}{\newcommand}
	\newcommand{\browntext}[1]{\textcolor{brown}{#1}}
	\newcommand{\greentext}[1]{\textcolor{green}{#1}}
	\newcommand{\redtext}[1]{\textcolor{red}{#1}}
	\newcommand{\bluetext}[1]{\textcolor{blue}{#1}}
	\newcommand{\brown}[1]{\browntext{ #1}}
	\newcommand{\green}[1]{\greentext{ #1}}
	\newcommand{\red}[1]{\redtext{ #1}}
	\newcommand{\blue}[1]{\bluetext{ #1}}
	
	\newcommand{\wtodo}{\todo[inline,color=orange!20, caption={}]}
	\newcommand{\lutodo}{\todo[inline,color=green!20, caption={}]}

	\title[$\mathrm i$Quantum groups and $\mathrm i$Hopf algebras I]{$\mathrm i$Quantum groups and $\mathrm i$Hopf algebras I: foundation}
	
	\author[Jiayi Chen]{Jiayi Chen}
	\address{Department of Mathematics, Shantou University, Shantou 515063, P.R.China}
	\email{chenjiayi@stu.edu.cn}
	
	\author[Ming Lu]{Ming Lu}
	\address{Department of Mathematics, Sichuan University, Chengdu 610064, P.R.China}
	\email{luming@scu.edu.cn}

	\author[Xiaolong Pan]{Xiaolong Pan}
	\address{Department of Mathematics, Sichuan University, Chengdu 610064, P.R.China}
	\email{xiaolong\_pan@stu.scu.edu.cn}

	\author[Shiquan Ruan]{Shiquan Ruan}
	\address{ School of Mathematical Sciences,
		Xiamen University, Xiamen 361005, P.R.China}
	\email{sqruan@xmu.edu.cn}

	\author[Weiqiang Wang]{Weiqiang Wang}
	\address{Department of Mathematics, University of Virginia, Charlottesville, VA 22904, USA}
	\email{ww9c@virginia.edu}
	
	\begin{abstract}
		We introduce the notion of iHopf algebra, a new associative algebra structure defined on a Hopf algebra equipped with a Hopf pairing. The iHopf algebra on a Borel quantum group endowed with a $\tau$-twisted Hopf pairing is shown to be a quasi-split universal iquantum group. In particular, the Drinfeld double quantum group is realized as the iHopf algebra on the double Borel. This iHopf approach allows us to develop new connections between Lusztig's braid group action and ibraid group action. It will further lead to the construction of dual canonical basis in a sequel. 
	\end{abstract}
	
	\maketitle
	\setcounter{tocdepth}{1}
	\tableofcontents

	\section{Introduction}
	
	
	Let $\U$ be the Drinfeld-Jimbo quantum group, $\tB$ be a Borel subalgebra of $\U$, and $\tU$ be the Drinfeld double on $\tB$. Associated to any quasi-split Satake diagram $(\I,\tau)$ (with no $\tau$-fixed edge), a universal quasi-split iquantum group $\tUi$ arises from an iHall algebra realization \cite{LW22a}. The Drinfeld double $\tU$ can be identified as a universal iquantum group associated to the diagonal Satake diagram $(\I \sqcup \I, \swa)$. Here are some basic facts about $\tUi$: it has the same size as the Borel $\tB$ and admits a $\Z\I$-grading. In contrast, the iquantum groups with parameters introduced originally by G.~Letzter \cite{Let99} (see Kolb \cite{Ko14}) have a smaller size, do not have a $\Z\I$-grading, and can be recovered by central reductions from $\tUi$. 
	
	The presentations for $\tU$ and $\tUi$ in terms of  ``dual" Chevalley generators used in this paper look a little different from the ones commonly used in the literature; they are most suitable  toward dual canonical bases in our sequel \cite{CLPRW} as well as in Berenstein-Greenstein \cite{BG17a} and in \cite{LP25}. 
	The algebras $\tU$ (and resp. $\tB$, or $\tUi$) admits variants $\hU$ (and resp. $\hB$, or $\hUi$) where the generators $K_i, K_i'$ (and resp. $K_i$, or $\K_i$) of the Cartan subalgebras are not required to be invertible. When discussing about braid group symmetries, these Cartan generators need to be invertible and we work with the tilde versions. The hat variant is natural from Hall algebra viewpoint and will be mostly used in our sequel on dual canonical bases, and so we can formulate both versions hand-in-hand in this paper; the reader can choose to skip the hat variant throughout.
	
	\vspace{2mm}
	
	In this paper, we shall formulate a notion of iHopf algebras, a new associative algebra structure defined on Hopf algebras with Hopf pairings. Our main motivation comes from Theorem \ref{thm:iHopfB=iQG} below: the iHopf algebra defined on the Borel $\tB$
	provides a realization of quasi-split universal iquantum groups. 
	In particular, the iHopf algebra defined on the double Borel $\tB\otimes \tB$ provides a realization of the Drinfeld double $\tU$. As a first application of the iHopf approach, we connect normalized Lusztig's braid group action on $\tU$ and relative braid (= ibraid) group action on $\tUi$.
	
	\vspace{2mm}
	Let us explain in detail. 
	Given a Hopf algebra $H$ with a Hopf pairing and a Hopf endomorphism $\tau$ on $H$, we define a new associative algebra structure on $H$, denoted by $H^\imath_\tau$. We also consider $(H\otimes H)^\imath$ associated to the tensor product Hopf algebra $H\otimes H$ with a twisted Hopf pairing \eqref{eq:twisted}. It is shown in Proposition \ref{prop:double-diagonal} that $(H\otimes H)^\imath$ can be identified with the Drinfeld double $D(H)$ of $H$. The notion of $\tau$-twisted compatible linear map $\chi:H\rightarrow \F$ can be found in Definition~ \ref{def:chi}. 
	
	\begin{alphatheorem} [Propositions \ref{prop:iHopf}, \ref{prop:embedding-iHopf-diagonal} and \ref{prop:coideal}]
		\label{thm:iHopfcoideal}
		There exist natural algebra homomorphisms 
		$\xi_{\tau,\chi}: H^\imath_\tau\rightarrow (H\otimes H)^\imath$ 
		and 
		$\Psi: H_\tau^\imath\rightarrow H^\imath_\tau\otimes (H\otimes H)^\imath$. 
	\end{alphatheorem}
	
	Under some mild assumption on the Hopf algebras (including quantum groups or more generally Nichols algebras of diagonal type), $\xi_{\tau,\chi}$ and $\Psi$ are injective. Hence according to Theorem~ \ref{thm:iHopfcoideal}, $H_\tau^\imath$ can be viewed as a coideal subalgebra of $D(H)$. 
	
	Given a quasi-split Satake diagram $(\I,\tau)$, we specialize the Hopf algebra $H$ to be the Borel subalgebra $\tB$ of a quantum group $\U$. In this case, we work with a specific $\chi$ which is determined by $\tau$ and we write $\widetilde{\xi}_{\tau}$ for $\xi_{\tau,\chi}$ below. 
	
	\begin{alphatheorem} [Lemmas \ref{lem:U=iHopfBB} and \ref{lem:iHembed}, Theorem \ref{thm:iH=iQG}]
		\label{thm:iHopfB=iQG}
		We have an algebra embedding $\widetilde{\xi}_\tau: \tB^\imath_\tau \rightarrow(\tB\otimes\tB)^\imath$ and
		algebra isomorphisms $\widetilde{\Phi}^\imath: \tB^\imath_\tau\stackrel{\cong}{\longrightarrow} \tUi$ and $\widetilde{\Phi}_{\sharp}:\tU \stackrel{\cong}{\longrightarrow} (\tB\otimes\tB )^\imath $ which fit into the following commutative diagram:
		\begin{equation*}  
			\begin{tikzcd}
				\tB^\imath_\tau\ar[r,"\widetilde{\xi}_\tau"]\ar[d,"\widetilde{\Phi}^\imath"]&(\tB\otimes\tB)^\imath\ar[d,"\widetilde{\Phi}_\sharp^{-1}"]\\
				\tUi\ar[r,hook]&\tU
			\end{tikzcd} 
		\end{equation*}
		Moreover, this commutative diagram admits an  integral form version over $\Z[v^{\frac12},v^{-\frac12}]$. 
	\end{alphatheorem}
	
	Let us compare our construction with another deformation construction of iquantum groups; for more details see Remark \ref{rem:KY}. Kolb and Yakimov \cite{KY20} used the quasi R-matrix to define a star product on a partial Bosonization  (which is $\U^-$ in split types) by a closed formula, and show that the resulting algebra is isomorphic to iquantum groups. Both constructions lead to almost the same recursive relations, but otherwise seem to be quite different; our iHopf construction requires only some basic input. 
	
	The iHopf realization of iquantum groups has several applications, and we present one to braid group action in this paper. There exist braid group symmetries, denoted by $\widetilde{T}_{i}$ and $\widetilde{T}_{w}$ for $i\in \I$ and $w\in W$, on quantum groups due to Lusztig \cite{Lus90b, Lus93}. In this paper, we use normalized variants $\widetilde{T}_{i}$ on the Drinfeld double $\tU$ which are compatible with dual Chevalley generators and a bar (anti)-involution on $\tU$; see Berenstein-Greenstein \cite{BG17a} or Lemma~ \ref{lem:QGbraid-bar}. Similarly, universal iquantum groups $\tUi$ admit ibraid group symmetries $\tTT_i$ and $\tTT_w$, for $i\in \I_\tau$ and $w\in W_\tau$ \cite{LW21a, LW22b, WZ23, Z23}; also cf. \cite{KP11}. We shall use variants $\tTT_i$ which are normalized to be compatible with dual Chevalley generators of $\tUi$; it turns out that $\tTT_i$ is compatible with a bar (anti-)involution on $\tUi$ as well; see Lemma \ref{lem:bar-invar-braid-Ui}.  
	
	There is a natural embedding of Lusztig's algebra $\ff$ into the Borel subalgebra, $\iota:\ff \rightarrow \tB$, and also an embedding $\hB \rightarrow \tB^\imath_\tau\equiv \tUi$. Recall $\vartheta_i$ in \eqref{rescale:theta} is a rescaling of a generator $\theta_i\in \ff$. For notation $r_i\in W$, see \eqref{def:ri}.
	
	\begin{alphatheorem} [Theorems \ref{thm:braidiHopf}]
		\label{thm:ibraid}
		For $i,j\in \I$ such that $i\neq j,\tau j$, we have $\tTT_i(\vartheta_j)=\iota\big(\widetilde{T}_{r_i}(\vartheta_j)\big).$
	\end{alphatheorem}
	
	Reflection functors in Hall algebras \cite{Rin96} and in iHall algebras \cite{LW21a,LW22b} have provided realizations of braid group action on quantum groups and ibraid group action on iquantum groups, respectively. 
	The reflection functors in iHall algebras coincide with their counterparts in Hall algebras when acting on modules of quivers, even though the underlying (Hall vs iHall) algebra structures are different. Theorem \ref{thm:ibraid} is a first indication that iHopf algebra construction captures the essence of iHall algebras to a large extent. 
	
	The proof of Theorem \ref{thm:ibraid} is carried out case-by-case,  since there are 3 rank one cases for quasi-split iquantum groups depending on whether the Cartan integers $c_{i,\tau i}=2, 0$ or $-1$. The constructions of $v$-root vectors in $\tUi$ via ibraid group action are quite involved \cite{CLW21b, Z23} as there are many local rank one and rank two types. It sheds new light on how $v$-root vectors in $\tUi$ and $\tU$ are identified via iHopf. 
	\vspace{2mm}
	
	The iHopf construction will be applied to provide a construction of dual canonical bases in the sequel \cite{CLPRW}. Via the iHopf mechanism, we shall be able to generalize the main results of the geometric approach in \cite{LW21b, LP25, LP26a, LP26b} of dual canonical basis of $\tUi$ from quasi-split ADE type (excluding type AIII$_{2r}$) to arbitrary finite type; see also \cite{Qin16} for dual canonical basis on the Drinfeld double $\tU$ of ADE type. In particular, we shall settle several major conjectures made by Berenstein-Greenstein \cite{BG17a} on dual canonical bases for $\tU$ of any finite type.

	\vspace{2mm}
	
	The paper is organized as follows. 
	In Section \ref{sec:i-Hopf}, we present a general construction of iHopf algebras and establish some of the basic properties.
	In Section \ref{sec:QG and iQG}, we review some constructions in quantum groups and iquantum groups, including (relative) braid group actions. 
	Then we show in Section \ref{sec:iQG=iHopfB} that the $\tau$-twisted iHopf algebra defined on the Borel subalgebra $\hB$ can be identified with the iquantum group $\hUi$. 
	In Section \ref{sec:braid}, we identify the braid and ibraid group actions on Chevalley generators through the iHopf construction.

	\vspace{2mm}
	\noindent{\bf Acknowledgments} 
	This paper owes its existence to Fan Qin who first suggested the idea of iHopf algebras. We thank Haicheng Zhang for communicating to us that he has also anticipated the iHopf realization of iquantum groups. We thank an anonymous expert for very helpful suggestions regarding the formula for the antipode $S^\imath$ and the inverse formula for the isomorphism $\Phi^{\texttt D}$.
	JC is partially supported by the National Natural Science Foundation of China Tianyuan Fund for Mathematics (No. 12526562). ML is partially supported by the National Natural Science Foundation of China (No. 12171333). SR is partially supported by
	Fundamental Research Funds for Central Universities of China (No. 20720250059), Fujian Provincial Natural Science Foundation of China (No. 2024J010006) and
	the National Natural Science Foundation of China (Nos. 12271448 and 12471035). WW is partially supported by the NSF grant DMS-2401351, and he thanks National University of Singapore (Department of Mathematics and IMS) for providing an excellent research environment and support during his sabbatical leave.  
	
	%

	\section{iHopf algebras}
	\label{sec:i-Hopf}
	
	In this section, we construct a new associative algebra structure on a Hopf algebra $H$ with a Hopf pairing, which shall be called an iHopf algebra and denoted by $H^\imath$. Under some suitable conditions, $H^\imath$ is further shown to be a (right coideal) subalgebra of $(H\otimes H)^\imath$.

	\subsection{Definition}
	All algebras in this section are over an arbitrary field $\F$. 
	
	Assume that $(H,m,\Delta,1,\varepsilon)$, or simply $H$ below, is a Hopf algebra with multiplication $m$, unit $1$, comultiplication $\Delta$ and counit $\varepsilon$. The antipode of $H$ is denoted by $S$, which will be assumed to be invertible. We write $ab=m(a,b)$ and use the Sweedler notation $\Delta(a)=\sum a_{(1)}\otimes a_{(2)}$ for $a,b\in H$. Define  $\Delta^{(2)}:=(\Delta\otimes\Id)\circ\Delta $, and then $\Delta^{(n)}$ by iteration. Denote 
	$\Delta^{(n)}(a)=\sum a_{(1)}\otimes a_{(2)}\otimes\cdots\otimes a_{(n+1)}$. 
	
	Let $\varphi$ be a Hopf pairing of $H$ (see e.g. \cite[\S3.2.1]{Jo95}), that is, $\varphi$ is a bilinear form satisfying the following conditions, for any $a,b,a',b'\in H$:
	\begin{itemize}
		\item[(1)] $\varphi(a,1)=\varepsilon(a),\;\varphi(1,b)=\varepsilon(b)$;
		\item[(2)] $\varphi(a,bb')=\varphi(\Delta(a),b\otimes b')$; 
		\item[(3)] $\varphi(aa',b)=\varphi(a\otimes a',\Delta(b))$;
		\item[(4)] $\varphi(a,S(b))=\varphi(S(a),b)$.
	\end{itemize}
	Here the pairing on the tensor product space is given by $\varphi(a \otimes a', b \otimes b') = \varphi(a, b)\varphi(a', b')$.
	
	\begin{definition}
		The {\em iHopf algebra}, denoted by $H^\imath =\text{iHopf}\, (H,\varphi)$, is the same vector space as $H$ equipped with a new multiplication: 
		\begin{align}\label{star product}
			a\ast b:=\sum \varphi(b_{(2)},a_{(1)})\cdot a_{(2)}b_{(1)},\quad \forall a,b\in H,
		\end{align}
		where $\Delta(a)=\sum a_{(1)}\otimes a_{(2)}$, $\Delta(b)=\sum b_{(1)}\otimes b_{(2)}$.
	\end{definition}

	\begin{proposition} \label{prop:iHopf}
		Let $(H,m,\Delta,1,\varepsilon)$ be a Hopf algebra with a Hopf pairing $\varphi$. Then $H^\imath$ is an associative algebra with the unit $1$.
	\end{proposition}
	
	\begin{proof}
		For any $a,b,c\in H^\imath$, we have
		\begin{align*}
			(a\ast b)\ast c&=\sum \varphi(b_{(2)},a_{(1)})\cdot (a_{(2)}b_{(1)})\ast c
			\\&=\sum \varphi(b_{(3)},a_{(1)})  \cdot\varphi(c_{(2)},a_{(2)}b_{(1)})\cdot a_{(3)}b_{(2)}c_{(1)}
			\\&=\sum \varphi(b_{(3)},a_{(1)})  \cdot\varphi(c_{(2)},a_{(2)})\cdot\varphi(c_{(3)},b_{(1)})\cdot a_{(3)}b_{(2)}c_{(1)}.
		\end{align*}
		On the other hand, we have
		\begin{align*}
			a\ast (b\ast c)&=\sum \varphi(c_{(2)},b_{(1)})\cdot a\ast (b_{(2)}c_{(1)})\\
			&=\sum \varphi(c_{(3)},b_{(1)})  \cdot \varphi(b_{(3)}c_{(2)},a_{(1)})\cdot a_{(2)}b_{(2)}c_{(1)}
			\\&=\sum \varphi(c_{(3)},b_{(1)}) \cdot\varphi (b_{(3)},a_{(1)}) \cdot\varphi(c_{(2)},a_{(2)})\cdot a_{(3)}b_{(2)}c_{(1)}.
		\end{align*}
		Hence $(a\ast b)\ast c=a\ast (b\ast c)$. It remains to prove that $1$ is the unit for $H^\imath$. In fact, by the definition of Hopf pairing, we have $a\ast 1=\sum\varphi(1,a_{(1)})\,  a_{(2)} =\sum\varepsilon(a_{(1)})\, a_{(2)}=a,$ and  
		$1\ast a=\sum\varphi(a_{(2)},1)\,  a_{(1)} =\sum\varepsilon(a_{(2)})\, a_{(1)}=a.$
	\end{proof}
	
	\begin{remark}
		The formulation of iHopf algebra $H^\imath$ is applicable to bialgebras $H$ equipped with bialgebra pairings 
		$\varphi$, meaning 
		$\varphi$
		satisfies only axioms (1)–(3) from the definition of Hopf pairings.
	\end{remark}
	
	\begin{remark}
		\label{rem:opposite}
		The space $H$ can be endowed with a different algebra structure with multiplication  
		$$a\ast' b=\sum\varphi(a_{(2)},b_{(1)})\cdot b_{(2)} a_{(1)}.$$
		One can show that 
		$$(H^\imath,\ast', 1)\cong (H^\imath,\ast, 1)^{\mathrm{op}}.$$
	\end{remark}
	
	Given a Hopf pairing $\varphi$, we shall construct several variants of Hopf pairings by twistings. Denote by $P:H\otimes H\to H\otimes H$, $a\otimes b\mapsto b\otimes a$ the permutation map.

	\begin{lemma}
		Assume that $H$ is a Hopf algebra with a Hopf pairing $\varphi$. Let $\tau:H\to H$ be any Hopf algebra endomorphism. Then
		
		\begin{itemize}
			\item[(1)] $\varphi\circ P$ is a Hopf pairing; 
			\item[(2)] 
			$\varphi \circ (\tau\otimes 1)$ and $\varphi \circ (1\otimes \tau)$ are Hopf pairings.
		\end{itemize}
	\end{lemma}
	
	\begin{proof}
		The statement (1) is clear. For (2), we only verify that $\varphi \circ (\tau\otimes 1)$ is a Hopf pairing. Indeed, we have $\varphi(\tau(1),b)=\varphi(1,b)=\varepsilon(b)$ and $\varphi(\tau(a),1)=\varepsilon(\tau(a))=\varepsilon(a)$. Moreover, we have
		\begin{align*}
			\varphi(\tau(a),bb')&=\varphi(\Delta(\tau(a)),b\otimes b')=\sum\varphi(\tau(a)_{(1)},b)\cdot\varphi(\tau(a)_{(2)},b')\\
			&=\sum\varphi(\tau(a_{(1)}),b)\cdot\varphi(\tau(a_{(2)}),b'),
		\end{align*}
		and
		\begin{align*}
			\varphi(\tau(aa'),b) &=\varphi(\tau(a) \tau (a'),b)=\varphi(\tau(a)\otimes \tau (a'),\Delta(b))
			\\
			&=\sum\varphi(\tau(a),b_{(1)})\cdot\varphi(\tau(a'),b_{(2)}),
		\end{align*}
		for any $a,a',b,b'\in H$. The proof is completed.
	\end{proof}
	
	Let $H$ and $'H$ be two Hopf algebras with Hopf pairings $\varphi$ and $\varphi'$, respectively, and denote by $H^\imath$, $'H^\imath$ the corresponding iHopf algebras. A Hopf algebra homomorphism $f: H\rightarrow{} 'H$ is said to preserve the Hopf pairings if  $\varphi(a,b)=\varphi'(f(a),f(b))$ for any $a,b\in H$. We denote by $f:H^\imath \rightarrow {}'H^\imath$ the same linear map as 
	$f: H\rightarrow{} 'H$.
	
	\begin{lemma}
		\label{lem:hom-iHopf}
		If $f: H\rightarrow{} 'H$ is a Hopf algebra homomorphism preserving the Hopf pairings, then $f: H^\imath \rightarrow {}'H^\imath$ is a homomorphism of algebras. 
	\end{lemma}
	
	\begin{proof}
		By definition, we have
		\begin{align*}
			f(a\ast b)=\sum \varphi(b_{(2)},a_{(1)})\cdot f(a_{(2)}b_{(1)})
			=\sum \varphi(b_{(2)},a_{(1)})\cdot f(a_{(2)})f(b_{(1)}).
		\end{align*}
		On the other hand, we have
		\begin{align*}
			f(a)\ast f(b)&=\sum \varphi'(f(b)_{(2)},f(a)_{(1)})\cdot f(a)_{(2)}f(b)_{(1)}
			\\&=\sum \varphi'(f(b_{(2)}),f(a_{(1)}))\cdot f(a_{(2)})f(b_{(1)})
			\\&=\sum \varphi(b_{(2)},a_{(1)})\cdot f(a_{(2)})f(b_{(1)}).
		\end{align*}
		Hence $f(a)\ast f(b)=f(a\ast b)$ and it is clear $f(1)=1$. Therefore, $f:H^\imath\rightarrow {}'H^\imath$ is an algebra homomorphism.
	\end{proof}



	%
	%
	\subsection{iHopf algebras of diagonal type}
	\label{subsec:i-Hopf-diagonal}
	
	Let $H$ be a Hopf algebra with a Hopf pairing $\varphi$. Then $H\otimes H$ is a Hopf algebra with its product, coproduct and counit given by
	\begin{align*}
		(a\otimes b)(c\otimes d) &=ac\otimes bd,
		\\
		\Delta(a\otimes b) &=\sum a_{(1)}\otimes b_{(1)}\otimes a_{(2)}\otimes b_{(2)},
		\\
		\varepsilon(a\otimes b) &=\varepsilon(a)\, \varepsilon(b),\qquad \forall a,b\in H.
	\end{align*}
	Moreover, $H\otimes H$ is endowed with a (twisted!) Hopf pairing $\varphi_{\sharp}$ given by  
	\begin{align} \label{eq:twisted}
		\varphi_{\sharp}(c\otimes d,a\otimes b)=\varphi(a,d)\, \varphi(c,b).
	\end{align}
	The iHopf algebra $(H\otimes H)^\imath$ associated with $(H\otimes H,\varphi_{\sharp})$ is called the iHopf algebra of diagonal type. Its multiplication is explicitly given by
	\[(a\otimes b)\ast (c\otimes d)=\sum\varphi (a_{(1)},d_{(2)})\cdot \varphi(c_{(2)},b_{(1)})\cdot a_{(2)}c_{(1)}\otimes b_{(2)}d_{(1)},\qquad\forall a,b,c,d\in H.
	\]
	
	\begin{proposition}
		\label{prop:ihopf-Drinfelddouble}
		Let $H$ be a Hopf algebra with a Hopf pairing $\varphi$. 
		Then $((H\otimes H)^\imath,\ast,1)$ can be made into a Hopf algebra with coproduct, counit and antipode given by 
		\begin{align*}
			\Delta^\imath(a\otimes b)&=\sum\varphi(a_{(2)},b_{(2)})\cdot (a_{(1)}\otimes b_{(3)})\otimes (a_{(3)}\otimes b_{(1)}),
			\\
			\varepsilon^\imath(a\otimes b)&=\varphi(a,S^{-1}(b)),\\
			S^\imath(a\otimes b)&=S(a)\otimes S^{-1}(b).
		\end{align*}
	\end{proposition}
	
	\begin{proof}
		We first verify that $\Delta^\imath$ is an algebra homomorphism. By definition, we have
		\begin{align*}
			&\Delta^\imath\big((a\otimes b)\ast (c\otimes d)\big)\\ &=\sum\varphi(a_{(1)},d_{(2)})\cdot \varphi(c_{(2)},b_{(1)})\cdot \Delta^\imath(a_{(2)}c_{(1)}\otimes b_{(2)}d_{(1)})       \\
			&=\sum\varphi(a_{(1)},d_{(4)})\cdot\varphi(c_{(4)},b_{(1)}) \cdot\varphi\big(a_{(3)}c_{(2)}, b_{(3)}d_{(2)}\big)        \\&\qquad\quad\cdot(a_{(2)}c_{(1)}\otimes b_{(4)}d_{(3)})\otimes
			(a_{(4)}c_{(3)}\otimes b_{(2)}d_{(1)})      \\&=\sum\varphi(a_{(1)},d_{(5)})\cdot\varphi(c_{(5)},b_{(1)}) \cdot\varphi(a_{(3)}, b_{(3)})\cdot
			\varphi(a_{(4)}, d_{(2)})\cdot
			\varphi(c_{(2)}, b_{(4)})\cdot
			\varphi(c_{(3)}, d_{(3)})\\
			&\qquad\quad\cdot(a_{(2)}c_{(1)}\otimes b_{(5)}d_{(4)})\otimes
			(a_{(5)}c_{(4)}\otimes b_{(2)}d_{(1)}).
		\end{align*}
		On the other hand, we have
		\begin{align*}
			&\Delta^\imath(a\otimes b)\ast \Delta^\imath(c\otimes d) \\
			&=\sum\big(\varphi(a_{(2)},b_{(2)})\cdot (a_{(1)}\otimes b_{(3)})\otimes (a_{(3)}\otimes b_{(1)})\big)
			\\&\qquad\quad\ast \big(\varphi(c_{(2)},d_{(2)})\cdot (c_{(1)}\otimes d_{(3)})\otimes (c_{(3)}\otimes d_{(1)})\big)        \\ &=\sum\varphi(a_{(2)},b_{(2)})\cdot\varphi(c_{(2)},d_{(2)})
			\\ &\quad\qquad\cdot\big((a_{(1)}\otimes b_{(3)})\ast (c_{(1)}\otimes d_{(3)})\big)
			\otimes\big((a_{(3)}\otimes b_{(1)})\ast (c_{(3)}\otimes d_{(1)})\big)
			\\ 
			&=\sum\varphi(a_{(3)},b_{(3)})\cdot\varphi(c_{(3)},d_{(3)})
			\cdot\big(\varphi(a_{(1)},d_{(5)})
			\varphi({c_{(2)}},b_{(4)})
			\cdot a_{(2)}c_{(1)}
			\otimes b_{(5)}
			d_{(4)}\big)
			\\ &\qquad\quad\otimes\big(\varphi(a_{(4)},d_{(2)})
			\cdot\varphi({c_{(5)}},b_{(1)})
			\cdot a_{(5)}c_{(4)}
			\otimes b_{(2)}
			d_{(1)}\big).
		\end{align*}
		Therefore $\Delta^\imath\big((a\otimes b)\ast (c\otimes d)\big) =\Delta^\imath(a\otimes b)\ast \Delta^\imath(c\otimes d)$.
		
		Next we show that $\Delta^\imath$ is coassociative. Indeed, we have 
		\begin{align*}
			&\big((\Delta^\imath\otimes 1\otimes1)\circ\Delta^\imath\big)(a\otimes b)
			\\
			&=\sum\varphi(a_{(2)},b_{(2)})\cdot (\Delta^\imath\otimes 1\otimes1)\big((a_{(1)}\otimes b_{(3)})\otimes (a_{(3)}\otimes b_{(1)})\big)
			\\
			&=\sum\varphi\big(a_{(4)},b_{(2)}\big) \cdot\varphi\big(a_{(2)},b_{(4)}\big)\cdot \big(a_{(1)}\otimes b_{(5)}\big)\otimes \big(a_{(3)}\otimes b_{(3)}\big)
			\otimes \big(a_{(5)}\otimes b_{(1)}\big),
		\end{align*}
		while
		\begin{align*}
			&\big((1\otimes1\otimes\Delta^\imath)\circ\Delta^\imath\big)(a\otimes b)
			\\
			&=\sum\varphi(a_{(2)},b_{(2)})\cdot (1\otimes1\otimes\Delta^\imath) \big((a_{(1)}\otimes b_{(3)})\otimes (a_{(3)}\otimes b_{(1)})\big)       \\ 
			&=\sum\varphi(a_{(2)},b_{(4)})\cdot(a_{(1)}\otimes b_{(5)})\otimes \Big(\varphi(a_{(4)},b_{(2)})\cdot (a_{(3)}\otimes b_{(3)})\otimes (a_{(5)}\otimes b_{(1)})\Big),
		\end{align*}
		which are equal to each other.
		
		To see that $\varepsilon^\imath$ is a counit for $\Delta^\imath$, we compute that
		\begin{align*}
			(\varepsilon^\imath\otimes 1)\Delta^\imath(a\otimes b)&=\sum\varphi(a_{(2)},b_{(2)})\cdot\varepsilon^\imath(a_{(1)}\otimes b_{(3)})\otimes(a_{(3)}\otimes b_{(1)})\\
			&=\sum\varphi(a_{(2)},b_{(2)})\cdot \varphi(a_{(1)},S^{-1}(b_{(3)}))\cdot a_{(3)}\otimes b_{(1)}\\
			&=\sum\varphi(a_{(1)},S^{-1}(b_{(3)})b_{(2)})\cdot a_{(2)}\otimes b_{(1)}\\
			&=\sum\varphi(a_{(1)},\varepsilon(b_{(2)}))\cdot a_{(2)}\otimes b_{(1)}
			\\
			&=\sum\varepsilon(a_{(1)})\cdot\varepsilon(b_{(2)})\cdot a_{(2)}\otimes b_{(1)}\\
			&=a\otimes b,
		\end{align*}
		and similarly $(1\otimes\varepsilon^\imath)\Delta^\imath(a\otimes b)=a\otimes b$. Moreover, $\varepsilon^\imath$ is an algebra homomorphism since
		\begin{align*}
			&\varepsilon^\imath((a\otimes b)\ast(c\otimes d))\\
			&=\sum\varphi(a_{(1)},d_{(2)})\cdot\varphi(c_{(2)},b_{(1)})\cdot\varphi(a_{(2)}c_{(1)},S^{-1}(b_{(2)}d_{(1)}))\\
			&=\sum\varphi(a_{(1)},d_{(2)})\cdot\varphi(c_{(2)},b_{(1)})\cdot\varphi(a_{(2)}c_{(1)},S^{-1}(d_{(1)})S^{-1}(b_{(2)}))\\
			&=\sum\varphi(a_{(1)},d_{(3)})\cdot\varphi(c_{(3)},b_{(1)})\cdot\varphi(S^{-1}(c_{(1)}),d_{(1)})\cdot\varphi(S^{-1}(a_{(2)}),d_{(2)})
			\\
			&\qquad\quad\cdot \varphi(S^{-1}(c_{(2)}),b_{(2)})\cdot\varphi(S^{-1}(a_{(3)}),b_{(3)})\\
			&=\sum\varphi(S^{-1}(a_{(2)})a_{(1)},d_{(2)})\cdot\varphi(c_{(3)}S^{-1}(c_{(2)}),b_{(1)})\cdot\varphi(S^{-1}(c_{(1)}),d_{(1)})\cdot\varphi(S^{-1}(a_{(3)}),b_{(2)})\\
			&=\sum\varepsilon(a_{(1)})\cdot\varepsilon(d_{(2)})\cdot\varepsilon(c_{(2)})\cdot\varepsilon(b_{(1)})\cdot\varphi(S^{-1}(c_{(1)}),d_{(1)})\cdot\varphi(S^{-1}(a_{(2)}),b_{(2)})\\
			&=\varphi(a,S^{-1}(b))\cdot\varphi(c,S^{-1}(d)).
		\end{align*}
		Therefore, $((H\otimes H)^\imath,\ast,1,\Delta^\imath,\varepsilon^\imath)$ forms a bialgebra.
		
		Finally, it remains to verify that $S^\imath$ is an antipode. To this end let us first expand
		\begin{align*}
			&(S^\imath\otimes 1)\Delta^\imath(a\otimes b)\\
			&=\sum\varphi(a_{(2)},b_{(2)})\cdot S^\imath(a_{(1)}\otimes b_{(3)})\otimes(a_{(3)}\otimes b_{(1)})\\
			&=\sum\varphi(a_{(2)},b_{(2)})\cdot S(a_{(1)})\otimes S^{-1}(b_{(3)})\otimes(a_{(3)}\otimes b_{(1)}).
		\end{align*}
		Let $m:(H\otimes H)^\imath\otimes(H\otimes H)^\imath\to (H\otimes H)^\imath$ denote the multiplication map. Then an immediate calculation, using the properties of the antipode $S$ of $H$, gives
		\begin{align*}
			&m(S^\imath\otimes 1)\Delta^\imath(a\otimes b)\\
			&=\sum\varphi(a_{(2)},b_{(2)})\cdot (S(a_{(1)})\otimes S^{-1}(b_{(3)}))\ast(a_{(3)}\otimes b_{(1)})\\
			&=\sum\varphi(a_{(3)},b_{(3)})\varphi(S(a_{(2)}),b_{(2)})\varphi(a_{(5)},S^{-1}(b_{(5)}))\cdot S(a_{(1)})a_{(4)}\otimes S^{-1}(b_{(4)})b_{(1)}\\
			&=\sum\varphi(S(a_{(2)})a_{(3)},b_{(2)})\varphi(a_{(5)},S^{-1}(b_{(4)}))\cdot S(a_{(1)})a_{(4)}\otimes S^{-1}(b_{(3)})b_{(1)}\\
			&=\sum\varepsilon(a_{(2)})\varepsilon(b_{(2)})\varphi(a_{(4)},S^{-1}(b_{(4)}))\cdot S(a_{(1)})a_{(3)}\otimes S^{-1}(b_{(3)})b_{(1)}\\
			&=\sum\varphi(a_{(3}),S^{-1}(b_{(3)}))\cdot S(a_{(1)})a_{(2)}\otimes S^{-1}(b_{(2)})b_{(1)}\\
			&=\sum\varphi(a_{(2)},S^{-1}(b_{(2)}))\cdot\varepsilon(a_{(1)})\varepsilon(b_{(1)})\cdot 1\otimes 1\\
			&=\varphi(a,S^{-1}(b))\cdot 1\otimes 1\\
			&=\varepsilon^\imath(a\otimes b)\cdot 1\otimes 1.
		\end{align*}
		A similar computation gives $m(1\otimes S^\imath)\Delta^\imath(a\otimes b)=\varepsilon^\imath(a\otimes b)\cdot 1\otimes 1$, hence $S^\imath$ is an antipode. This completes the proof of the proposition.
	\end{proof}
	
	\subsection{Connection to Drinfeld doubles}
	
	There is a construction of Drinfeld doubles; see \cite{Dr87}, \cite[\S3.2]{Jo95}, where one starts with a skew-Hopf pairing. We will follow a variant used in \cite{Cr10,LP21}, which starts with a Hopf pairing. Given a Hopf algebra $H$ with a Hopf pairing $\varphi$, its Drinfeld double $D(H)$ is the associative algebra defined on $H\otimes H$ subject to the following relations for any $a,b,a',b'\in H$:
	\begin{itemize}
		\item[(D1)]  $(a\otimes 1)(a'\otimes1)=aa'\otimes1$;
		\item[(D2)] $(1\otimes b)(1\otimes b')=1\otimes bb'$;
		\item[(D3)] $(a\otimes 1)(1\otimes b)=a\otimes b$;
		\item[(D4)] $\sum \varphi(a_{(1)},b_{(2)})\cdot(1\otimes b_{(1)})(a_{(2)}\otimes 1)=\sum \varphi(a_{(2)},b_{(1)})\cdot (a_{(1)}\otimes b_{(2)})$.
	\end{itemize}
	See \cite[Remark 3.2.2]{Jo95}.
	
	Moreover, $D(H)$ is a Hopf algebra with comultiplication, counit and antipode given by
	\begin{align*}
		\Delta(a\otimes b)&=\sum (a_{(1)}\otimes b_{(2)})\otimes (a_{(2)}\otimes b_{(1)}),\\
		\varepsilon(a\otimes b)&=\varepsilon(a)\varepsilon(b),&\\
		S(a\otimes b)&=(1\otimes S^{-1}(b))(S(a)\otimes 1).
	\end{align*}

	\begin{proposition}
		\label{prop:double-diagonal} 
		Let $H$ be a Hopf algebra with a Hopf pairing $\varphi$. There exists an isomorphism of Hopf algebras
		\begin{align}
			\Phi^{\textup{\texttt{D}}}:D(H)\longrightarrow (H\otimes H)^\imath,\quad (a\otimes b)\mapsto (a\otimes 1)\ast(1\otimes b),
		\end{align}
		with inverse given by
		\begin{align}  \label{inversePhi}     (\Phi^{\textup{\texttt{D}}})^{-1}:(H\otimes H)^\imath\longrightarrow D(H),\quad a\otimes b\mapsto \sum\varphi(a_{(1)},S^{-1}(b_{(2)}))\cdot a_{(2)}\otimes b_{(1)}.
		\end{align}
	\end{proposition}
	
	\begin{proof}
		We first verify that the relations (D1)--(D4) for $D(H)$ are preserved by $\Phi^\texttt{D}$. 
		
		For any $a,a'\in H$ we note that 
		\begin{align*}
			\Phi^\texttt{D}(a\otimes 1)\ast \Phi^\texttt{D}(a'\otimes 1)&=(a\otimes 1)\ast(a'\otimes 1)=aa'\otimes 1=\Phi^\texttt{D}((a\otimes 1)(a'\otimes 1)).
		\end{align*}
		Similarly, the relation $\Phi^\texttt{D}(1\otimes b)\ast \Phi^\texttt{D}(1\otimes b')=\Phi^\texttt{D}((1\otimes b)*(1\otimes b'))$.
		
		As for (D4), note that by the definition of iHopf algebra, we have
		\[
		\sum \varphi(a_{(1)},b_{(2)})\cdot(1\otimes b_{(1)})\ast(a_{(2)}\otimes 1)=\sum \varphi(a_{(1)},b_{(3)})\varphi(a_{(3)},b_{(1)})\cdot a_{(2)}\otimes b_{(2)}.
		\]
		On the other hand,
		\[
		\sum \varphi(a_{(2)},b_{(1)})\cdot(a_{(1)}\otimes 1)\ast(1\otimes b_{(2)})=\sum \varphi(a_{(3)},b_{(1)})\varphi(a_{(1)},b_{(3)})\cdot a_{(2)}\otimes b_{(2)}.
		\]
		The above right-hand sides are equal. Therefore, $\Phi^\texttt{D}$ is an algebra homomorphism.
		
		Note that $\Phi^\texttt{D}$ is also a coalgebra homomorphism because
		\begin{align*}
			\Delta^\imath\Phi^\texttt{D}(a\otimes b)&=\sum\varphi(a_{(1)},b_{(2)})\Delta^\imath(a_{(2)}\otimes b_{(1)})\\
			&=\sum\varphi(a_{(1)},b_{(4)})\cdot \varphi(a_{(3)},b_{(2)})\cdot (a_{(2)}\otimes b_{(3)})\otimes(a_{(4)}\otimes b_{(1)}),
		\end{align*}
		while
		\begin{align*}
			(\Phi^\texttt{D}\otimes \Phi^\texttt{D})\Delta(a\otimes b)&=\sum\Phi^\texttt{D}(a_{(1)}\otimes b_{(2)})\otimes\Phi^\texttt{D}(a_{(2)}\otimes b_{(1)})\\
			&=\sum\varphi(a_{(1)},b_{(4)})\cdot\varphi(a_{(3)},b_{(2)})\cdot(a_{(2)}\otimes b_{(3)})\otimes(a_{(4)},b_{(1)}). 
		\end{align*}
		Also, it is compatible with the counits of both sides:
		\begin{align*}
			\varepsilon^\imath\Phi^\texttt{D}(a\otimes b)&=\sum\varphi(a_{(1)},b_{(2)})\cdot\varphi(a_{(2)},S^{-1}(b_{(1)}))\\
			&=\sum\varphi(a,b_{(2)}S^{-1}(b_{(1)}))=\varphi(a,\varepsilon(b))\\
			&=\varepsilon(a)\cdot\varepsilon(b)\\
			&=\varepsilon(a\otimes b).
		\end{align*}
		Therefore, $\Phi^\texttt{D}:D(H)\to (H\otimes H)^\imath$ is a homomorphism of Hopf algebras.
		
		It then remains to show that $\Phi^\texttt{D}$ is an isomorphism with inverse given by \eqref{inversePhi}. Let us denote the map in \eqref{inversePhi} for now by $\Psi^{\textup{\texttt{D}}}$. A direct calculation gives
		\begin{align*}
			\Psi^{\textup{\texttt{D}}}\Phi^{\textup{\texttt{D}}}(a\otimes b)&=\sum\varphi(a_{(1)},b_{(2)})\Psi^{\textup{\texttt{D}}}(a_{(2)}\otimes b_{(1)})\\
			&=\sum\varphi(a_{(1)},b_{(3)})\varphi(a_{(2)},S^{-1}(b_{(2)}))\cdot a_{(3)}\otimes b_{(1)}\\
			&=\sum\varphi(a_{(1)},b_{(3)}S^{-1}(b_{(2)}))\cdot a_{(3)}\otimes b_{(1)}\\
			&=\sum\varphi(a_{(1)},\varepsilon(b_{(2)}))\cdot a_{(2)}\otimes b_{(1)}\\
			&=\sum\varepsilon(a_{(1)})\varepsilon(b_{(2)})\cdot a_{(2)}\otimes b_{(1)} =a\otimes b,
		\end{align*}
		while
		\begin{align*}
			\Phi^{\textup{\texttt{D}}}\Psi^{\textup{\texttt{D}}}(a\otimes b)&=\sum\varphi(a_{(1)},S^{-1}(b_{(2)}))\Phi^{\textup{\texttt{D}}}(a_{(2)}\otimes b_{(1)})\\
			&=\sum\varphi(a_{(1)},S^{-1}(b_{(3)}))\varphi(a_{(2)},b_{(2)})\cdot a_{(3)}\otimes b_{(1)}\\
			&=\sum\varphi(a_{(1)},S^{-1}(b_{(3)})b_{(2)})\cdot a_{(3)}\otimes b_{(1)}\\
			&=\sum\varphi(a_{(1)},\varepsilon(b_{(2)}))\cdot a_{(2)}\otimes b_{(1)}\\
			&=\sum\varepsilon(a_{(1)})\varepsilon(b_{(2)})\cdot a_{(2)}\otimes b_{(1)} =a\otimes b.
		\end{align*}
		This completes the proof of the proposition.
	\end{proof}

	\subsection{A coideal subalgebra structure via twists}
	
	\begin{definition} \label{def:chi}
		Let $H$ be a Hopf algebra with a Hopf pairing $\varphi$, and $\tau$ a Hopf algebra endomorphism of $H$. A linear map $\chi:H\rightarrow \F$ is {\em $\tau$-twisted compatible} if the identity $\chi(ab)=\sum\chi(a_{(1)})\chi(b_{(2)})\varphi(\tau (a_{(2)}),b_{(1)})$ holds for all $a,b\in H$.
	\end{definition}
	We just call  $\chi$  a compatible map if it is  $\tau$-twisted compatible with $\tau=\Id$.
	We refer to Lemma \ref{lem:fvartheta} for a natural example of $\chi$ when $H$ is a Borel quantum group.

	
	
	For a Hopf algebra endomorphism $\tau$ of $H$, we denote by $H^\imath_\tau$  the iHopf algebra associated with $H$ and the Hopf pairing $\varphi \circ(\tau \otimes 1)$. 
	
	\begin{proposition}
		\label{prop:embedding-iHopf-diagonal}
		Let $H$ be a Hopf algebra with a Hopf pairing $\varphi$, and $\chi:H\rightarrow \F$ be a $\tau$-twisted compatible map for a Hopf algebra endomorphism $\tau$ of $H$. Then there is an algebra homomorphism 
		\begin{align} \label{xi}
			\xi_{\tau,\chi}: H^\imath_\tau\longrightarrow (H\otimes H)^\imath,
			\qquad  a \mapsto \sum \chi(a_{(2)})\cdot \tau(a_{(3)})\otimes a_{(1)}.
		\end{align}
	\end{proposition}
	
	\begin{proof}
		Since $a\ast b=\sum\varphi(\tau(b_{(2)}),a_{(1)})\cdot  a_{(2)}b_{(1)}$ in $H^\imath_\tau$, we have 
		\begin{align*}
			\xi_{\tau,\chi}(a\ast b)&=\sum\varphi(\tau(b_{(2)}),a_{(1)})\cdot  \xi_{\tau,\chi}(a_{(2)}b_{(1)})
			\\&=\sum\varphi(\tau(b_{(4)}),a_{(1)})\chi(a_{(3)}b_{(2)})\cdot \tau(a_{(4)}b_{(3)})\otimes a_{(2)}b_{(1)}
			\\&=\sum\varphi(\tau(b_{(5)}),a_{(1)})\cdot \chi(a_{(3)})\chi(b_{(3)})\varphi(\tau(a_{(4)}),b_{(2)})\cdot \tau(a_{(5)}b_{(4)})\otimes a_{(2)}b_{(1)}.
		\end{align*}
		On the other hand, we have
		\begin{align*}
			\xi_{\tau,\chi}(a)\ast\xi_{\tau,\chi}( b) 
			&=\sum\big( \chi(a_{(2)})\cdot \tau(a_{(3)})\otimes a_{(1)}\big)\ast \big( \chi(b_{(2)})\cdot \tau(b_{(3)})\otimes b_{(1)}\big)\\
			&=\sum\chi(a_{(2)})\chi(b_{(2)})\cdot\big( \tau(a_{(3)})\otimes a_{(1)}\big)\ast \big(\tau( b_{(3)})\otimes b_{(1)}\big)
			\\&=\sum\chi(a_{(3)})\chi(b_{(3)}) \varphi(\tau(a_{(4)}),b_{(2)})\varphi(\tau(b_{(5)}),a_{(1)})\cdot \tau (a_{(5)}) \tau(b_{(4)})\otimes a_{(2)}b_{(1)}.
		\end{align*}
		It follows that $\xi_{\tau,\chi}(a)\ast\xi_{\tau,\chi}( b) =\xi_{\tau,\chi}(a\ast b)$.
	\end{proof}
	
	\begin{proposition}
		\label{prop:coideal}
		Let $(H,m,\Delta,1,\varepsilon)$ be a Hopf algebra with a Hopf pairing $\varphi$, and $\chi$ be a $\tau$-twisted compatible map for a Hopf algebra endomorphism $\tau$ of $H$. Then there is a right $(H\otimes H)^\imath$-coaction
		\begin{align*}
			\Psi: H_\tau^\imath&\longrightarrow H^\imath_\tau\otimes (H\otimes H)^\imath,
			\\
			a&\mapsto \sum \varphi(\tau(a_{(4)}),a_{(2)})\cdot a_{(3)}\otimes (\tau(a_{(5)})\otimes a_{(1)})
		\end{align*}
		such that $H^\imath_\tau$ is a comodule algebra. Moreover, we have
		$$ ( \xi_{\tau,\chi}\otimes 1\otimes1)\circ\Psi=\Delta^\imath\circ \xi_{\tau,\chi}.$$
	\end{proposition}
	
	\begin{proof}
		We first show that $\Psi$ is a right coaction and an algebra homomorphism. For this we compute that
		\begin{align*}
			(1\otimes\Delta^\imath)\Psi(a)&=\sum\varphi(\tau a_{(4)},a_{(2)})\cdot a_{(3)}\otimes\Delta^\imath(\tau a_{(5)}\otimes a_{(1)})\\
			&=\sum\varphi(\tau a_{(6)},a_{(4)})\varphi(\tau a_{(8)},a_{(2)})\cdot a_{(5)}\otimes(\tau a_{(7)}\otimes a_{(3)})\otimes(\tau a_{(9)}\otimes a_{(1)}),
		\end{align*}
		\begin{align*}
			(\Psi\otimes 1)\Psi(a)&=\sum\varphi(\tau a_{(4)},a_{(2)})\cdot\Psi(a_{(3)})\otimes(\tau a_{(5)}\otimes a_{(1)})\\
			&=\sum\varphi(\tau a_{(8)},a_{(2)})\varphi(\tau a_{(6)},a_{(4)})\cdot a_{(5)}\otimes(\tau a_{(7)}\otimes a_{(3)})\otimes(\tau a_{(9)}\otimes a_{(1)}),
		\end{align*}
		and
		\begin{align*}
			(1\otimes\varepsilon)\Psi(a)&=\sum \varphi(\tau(a_{(4)}),a_{(2)})\cdot a_{(3)}\otimes \varepsilon(\tau(a_{(5)})\otimes a_{(1)})\\
			&=\sum \varphi(\tau(a_{(4)}),a_{(2)})\varphi(\tau a_{(5)},S^{-1}(a_{(1)})\cdot a_{(3)}\\
			&=\sum\varphi(\tau a_{(4)},a_{(2)}S^{-1}(a_{(1)}))\cdot a_{(3)}\\
			&=\sum\varepsilon(\tau a_{(3)})\varepsilon(a_{(1)})\cdot a_{(2)}\\
			&=a.
		\end{align*}
		So $\Psi$ is a right coaction. On the other hand, note that
		\begin{align*}
			\Psi(a\ast b)&=\Psi\Big(\sum\varphi(\tau(b_{(2)}),a_{(1)})\cdot  a_{(2)}b_{(1)}\Big)\\
			&=\sum\varphi(\tau(b_{(6)}),a_{(1)})\varphi(\tau(a_{(5)}b_{(4)}),a_{(3)}b_{(2)})\cdot a_{(4)}b_{(3)}\otimes (\tau( a_{(6)}b_{(5)})\otimes a_{(2)}b_{(1)} )\\
			&=\sum\varphi(\tau(b_{(7)}),a_{(1)})\varphi(\tau(a_{(5)}b_{(4)}),a_{(3)})\varphi(\tau(a_{(6)}b_{(5)}),b_{(2)})\\&\qquad\cdot a_{(4)}b_{(3)}\otimes ( \tau(a_{(7)}b_{(6)})\otimes a_{(2)}b_{(1)}\\
			&=\sum\varphi(\tau(b_{(8)}),a_{(1)})\varphi(\tau(a_{(6)}),a_{(3)})\varphi(\tau(b_{(5)}),a_{(4)})\varphi(\tau(a_{(7)}),b_{(2)})\varphi(\tau(b_{(6)}),b_{(3)})\\&\qquad\cdot a_{(5)}b_{(4)}\otimes ( \tau(a_{(8)}b_{(7)})\otimes a_{(2)}b_{(1)} ),
		\end{align*}
		while 
		\begin{align*}
			&\Psi(a)\ast\Psi(b)
			\\
			&=\sum\varphi(\tau(a_{(4)}),a_{(2)})\varphi(\tau(b_{(4)}),b_{(2)})\cdot (a_{(3)}\ast b_{(3)})\otimes \big(( \tau(a_{(5)})\otimes a_{(1)})\ast ( \tau(b_{(5)})\otimes b_{(1)})\big)\\
			&=\sum\varphi(\tau(a_{(6)}),a_{(3)})\varphi(\tau(b_{(6)}),b_{(3)})\varphi(\tau(b_{(5)}),a_{(4)}\varphi(\tau(b_{(8)}),a_{(1)})\varphi(\tau(a_{(7)}),b_{(2)})\\&\qquad\cdot a_{(5)}b_{(4)}\otimes ( \tau(a_{(8)}b_{(7)})\otimes a_{(2)}b_{(1)}).
		\end{align*}
		Therefore, $\Psi$ is an algebra homomorphism and $H^\imath_\tau$ is a comodule algebra. 
		
		Now the second statement follows from
		\begin{align*}
			&      \Big( ( \xi_{\tau,\chi}\otimes 1\otimes1)\circ\Psi\Big)(a)
			\\     &=(\xi_{\tau,\chi}\otimes1\otimes 1)\Big(\sum \varphi(\tau(a_{(4)}),a_{(2)})\cdot a_{(3)}\otimes (\tau(a_{(5)})\otimes a_{(1)})\Big)
			\\&=\sum \varphi(\tau(a_{(6)}),a_{(2)}) 
			\chi(a_{(4)})\cdot ( \tau(a_{(5)})\otimes a_{(3)})
			\otimes (\tau(a_{(7)})\otimes a_{(1)})
		\end{align*}
		and
		\begin{align*}
			&       ( \Delta^\imath\circ \xi_{\tau,\chi})(a)
			\\      &=\sum \chi(a_{(2)})\cdot \Delta^\imath(\tau(a_{(3)})\otimes a_{(1)})
			\\&=\sum \chi(a_{(4)})\varphi(\tau(a_{(6)}),a_{(2)})\cdot (\tau(a_{(5)})\otimes a_{(3)} )\otimes (\tau(a_{(7)})\otimes a_{(1)}).
		\end{align*}
		This completes the proof.
	\end{proof}
	
	\begin{remark}
		If the algebra homomorphism $\xi_{\tau,\chi}:H^\imath_\tau\rightarrow (H\otimes H)^\imath$ is injective, then we can identify $H^\imath_\tau$ as a subalgebra of $(H\otimes H)^\imath$.
		In this case, Proposition \ref{prop:coideal}
		says that $H_\tau^\imath$ can be identified as a right coideal subalgebra of $(H\otimes H)^\imath$; that is, the homomorphism  $\Psi:H_\tau^\imath\rightarrow H^\imath_\tau\otimes (H\otimes H)^\imath$ is the restriction of the coproduct $\Delta^\imath$ of $(H\otimes H)^\imath$ to $H^\imath_\tau$. 
	\end{remark}

	\section{Quantum groups and $\mathrm{i}$quantum groups} 
	\label{sec:QG and iQG}
	
	In this preliminary section, we formulate quantum groups and iquantum groups in terms of (less standard) dual generators. We also review and formulate (relative) braid group action in this rescaled setting. 
	
	\subsection{Quantum groups}
	\label{subsec:QG}
	
	Let $\I=\{1,\dots,n\}$ be the index set. 
	Let $C=(c_{ij})_{i,j \in \I}$ be the symmetrizable Cartan matrix of a Kac-Moody Lie algebra $\fg$. Its Dynkin diagram is  denoted by $\Gamma$. Let $D=\diag(d_i\mid i\in \I)$ with $d_i\in\Z_{>0}$ be the symmetrizer of $C$, i.e., $DC$ is symmetric. 
	Let $\Delta^+=\{\alpha_i\mid i\in\I\}$ be the set of simple roots of $\fg$, and denote the root lattice by $\Z^{\I}:=\Z\alpha_1\oplus\cdots\oplus\Z\alpha_n$. We define a symmetric bilinear form on $\Z^\I$ by setting
	\begin{align} \label{BilForm}
		(\alpha_i,\alpha_j)=d_ic_{ij},\quad \forall i,j\in\I.
	\end{align}

	Let $\Phi^+$ be the set of positive roots. The simple reflection $s_i:\Z^{\I}\rightarrow\Z^{\I}$ is defined to be $s_i(\alpha_j)=\alpha_j-c_{ij}\alpha_i$, for $i,j\in \I$.
	Denote the Weyl group by $W =\langle s_i\mid i\in \I\rangle$.
	
	Let $v$ be an indeterminate. Let 
	$$v_i=v^{d_i},\qquad \forall i\in\I.$$
	For $A,B$ in a $\Q(v^{\frac12})$-algebra, we write $[A, B]=AB-BA$, and $[A,B]_q=AB-qBA$ for any $q\in\Q(v^{\frac12})^\times$. Denote, for $r\in\N,m \in \Z$,
	\[
	[r]_{v_i}=\frac{v_i^r-v_i^{-r}}{v_i-v_i^{-1}},
	\quad
	[r]_{v_i}^!=\prod_{i=1}^r [i]_{v_i}, \quad \qbinom{m}{r}_{v_i} =\frac{[m]_{v_i}[m-1]_{v_i}\ldots [m-r+1]_{v_i}}{[r]_{v_i}^!}.
	\]
	Following \cite{Dr87, BG17a}, the (Drinfeld double) quantum group $\hU := \hU_v(\fg)$ is defined to be the $\Q(v^{\frac12})$-algebra generated by $E_i,F_i, \tK_i,\tK_i'$, $i\in \I$, subject to the following relations:  for $i, j \in \I$,
	\begin{align}
		[E_i,F_j]= \delta_{ij}(v_i^{-1}-v_i) (\tK_i-\tK_i'),  &\qquad [\tK_i,\tK_j]=[\tK_i,\tK_j']  =[\tK_i',\tK_j']=0,
		\label{eq:KK}
		\\
		\tK_i E_j=v_i^{c_{ij}} E_j \tK_i, & \qquad \tK_i F_j=v_i^{-c_{ij}} F_j \tK_i,
		\label{eq:EK}
		\\
		\tK_i' E_j=v_i^{-c_{ij}} E_j \tK_i', & \qquad \tK_i' F_j=v_i^{c_{ij}} F_j \tK_i',
		\label{eq:K2}
	\end{align}
	and for $i\neq j \in \I$,
	\begin{align}
		& \sum_{r=0}^{1-c_{ij}} (-1)^r \left[ \begin{array}{c} 1-c_{ij} \\r \end{array} \right]_{v_i}  E_i^r E_j  E_i^{1-c_{ij}-r}=0,
		\label{eq:serre1} 
		\\ 
		&\sum_{r=0}^{1-c_{ij}} (-1)^r \left[ \begin{array}{c} 1-c_{ij} \\r \end{array} \right]_{v_i}  F_i^r F_j  F_i^{1-c_{ij}-r}=0.
		\label{eq:serre2}
	\end{align}
	
	We define $\tU=\tU_v(\fg)$ as the $\Q(v^{\frac12})$-algebra with   generators and relations of $\widehat{\U}$ above, but in addition requiring $\tK_i,\tK_i'$ ($i\in\I$) to be invertible. Then $\tU$ and $\hU$ are $\Z^\I$-graded algebras by setting 
	$$\deg E_i=\alpha_i, \qquad\deg F_i=-\alpha_i,\qquad \deg K_i=0=\deg K_i'.$$ 
	Let $\tU_\mu$ be the homogeneous subspace of degree $\mu$.  Then $\tU=\oplus_{\mu\in\Z^\I} \tU_\mu$ and $\hU=\oplus_{\mu\in\Z^\I} \hU_\mu$. 
	
	The Drinfeld-Jimbo quantum group $\bU$ is defined to the $\Q(v^{\frac12})$-algebra generated by $E_i,F_i, K_i, K_i^{-1}$, $i\in \I$, subject to the relations modified from \eqref{eq:KK}--\eqref{eq:serre2} with $\tK_i'$ replaced by $K_i^{-1}$. We can also view $\bU$ as the quotient algebra of $\hU$ (or $\tU$) modulo the ideal generated by $K_iK_i'-1$ ($i\in\I$); see \cite{Dr87}. 
	
	\begin{remark}
		Note that the presentation of $\tU$ here looks a bit different as it uses ``dual" generators; it is isomorphic to the usual Drinfeld double quantum group (see e.g. \cite{LW22a}). In fact, denote 
		\begin{align}
			\label{eq:Udj-gen}
			\ce_i=\frac{E_i}{v_i-v_i^{-1}},\qquad \cf_i=\frac{F_i}{v_i^{-1}-v_i},\qquad \forall i\in\I.
		\end{align} 
		(Our sign convention in \eqref{eq:Udj-gen} is opposite to the one adopted in \cite{LP25}.) 
		The presentation of $\tU$ given in \cite{LW22a} is in terms of generators $\ce_i,\cf_i,K_i,K_i'$ ($i\in\I$) and by a central deduction it descends to $\U$ in the standard presentation. Similarly, the presentation $\U$ given here is different from the usual presentation of the Drinfeld-Jimbo quantum group. 
		This non-standard presentation of the quantum groups $\tU$, $\hU$ and $\U$ can be viewed as quantized coordinate algebras of the Poisson dual groups of Lie groups; see \cite{Dr87} and also \cite[Remark 1.2]{BG17a}.
	\end{remark}
	
	By a slight abuse of notation, let $\U^+$ be the subalgebra of $\hU$ (and also $\tU$, $\U$) generated by $E_i$ $(i\in \I)$, and let $\U^-$ be the subalgebra generated by $F_i$ $(i\in\I)$. Let $\hU^0$ and $\tU^0$ be the subalgebras of $\hU$ and $\tU$ generated by $\tK_i, \tK_i'$ $(i\in \I)$, and $\U^0$ be the subalgebra of $\bU$ generated by $\tK_i^{\pm 1}$ $(i\in \I)$. Then the algebras $\hU$, $\widetilde{\bU}$ and $\bU$ have triangular decompositions:
	\begin{align*}
		\hU=\U^+\otimes\hU^0\otimes \U^-,\qquad 
		\widetilde{\bU} =\U^+\otimes \widetilde{\bU}^0\otimes\U^-,
		\qquad
		\bU &=\bU^+\otimes \bU^0\otimes\bU^-.
	\end{align*}
	Clearly, 
	$\hU^0\cong\Q(v^{\frac12})[\tK_i,\tK_i'\mid i\in\I]$, $\tU^0\cong \Q(v^{\frac12})[\tK_i^{\pm1},(\tK_i')^{\pm1}\mid i\in\I]$ and $\bU^0\cong\Q(v^{\frac12})[K_i^{\pm1}\mid i\in\I]$. Note that 
	${\bU}^0 \cong \hU^0/\langle\tK_i \tK_i' -1 \mid   i\in \I\rangle$. 
	For any $\mu=\sum_{i\in\I}m_i\alpha_i\in\Z^\I$, we denote $K_\mu=\prod_{i\in\I} K_i^{m_i}$, $K_\mu'=\prod_{i\in\I} (K_i')^{m_i}$.

	The algebras $\tU$ is a Hopf algebra, with the coproduct $\Delta$, the counit $\varepsilon$ and the antipode $S$ defined by
	\begin{align}
		\begin{split}
			\Delta(E_i)=E_i\otimes 1+K_i\otimes E_i,\quad &\Delta(F_i)=1\otimes F_i+F_i\otimes K_i',
			\\
			\Delta(K_i)=K_i\otimes K_i,\quad &\Delta(K_i')=K_i'\otimes K_i';
			\\
			\varepsilon(E_i)=0=\varepsilon(F_i),\qquad &\varepsilon(K_i)=1=\varepsilon(K_i');
			\\
			S(E_i)=-K_i^{-1}E_i,\quad S(F_i)=-F_iK_i'^{-1},\quad &S(K_i)=K_i^{-1},\quad S(K_i')=K_i'^{-1}.
		\end{split}
	\end{align}
	The algebra $\U$ is also a Hopf algebra such that the natural projection $\widetilde{\U}\rightarrow \U$ is a morphism of Hopf algebras. 
	And the algebra $\widehat{\U}$ is a bialgebra, whose coproduct $\Delta$ and the counit $\varepsilon$ are defined via identical formulas.
	
	\begin{lemma} \label{QG bar-involution def}
		There exists an anti-involution (called the bar-involution) $u\mapsto \ov{u}$ on $\hU$ (and also $\tU$, $\U$) given by $\ov{v^{1/2}}=v^{-1/2}$, $\ov{E_i}=E_i$, $\ov{F_i}=F_i$, and $\ov{K_i}=K_i$, $\ov{K_i'}=K_i'$, for $i\in\I$.
	\end{lemma}
	
	\begin{proof}
		Clear. 
	\end{proof}
	
	\begin{lemma}[cf. \cite{Lus93}]
		\label{lem:anti-involut-QG}
		There exists an anti-involution $\sigma$ on $\hU$ (also $\tU$, $\U$) given by $\sigma(E_i)=E_i$, $\sigma(F_i)=F_i$, and $\sigma(K_i)=K_i'$, for $i\in\I$.
	\end{lemma}
	
	\begin{lemma} \label{lem:twisting Psi}
		Let $\F$ be the algebraic closure of $\Q(v^{\frac{1}{2}})$ and $\F^\times=\F\setminus\{0\}$. For scalars $\ba=(a_i)_{i\in\I}\in (\F^{\times})^\I$, we have an automorphism $\widetilde{\Psi}_{\ba}$ on the $\F$-algebra $\tU$ such that
		\[
		\widetilde{\Psi}_{\ba}:K_i\mapsto a_i^{\frac{1}{2}}K_i,\quad K'_i\mapsto a_i^{\frac{1}{2}}K'_i,\quad E_i\mapsto a_i^{\frac{1}{2}}E_i,\quad F_i\mapsto F_i.
		\]
	\end{lemma}
	
	Let $\Br(W)$ be the braid group associated to the Weyl group $W$, generated by simple reflections $t_i$ ($i\in\I$). Lusztig introduced 4 variants of braid group symmetries on the quantum group $\U$ in \cite{Lus90b} and \cite[\S37.1.3]{Lus93}. These braid group symmetries can be lifted to the Drinfeld double $\tU$; see, e.g., \cite[Propositions 6.20–6.21]{LW22b}, which are denoted by $\widetilde{T}_{i,e}',\widetilde{T}_{i,-e}''$, $e=\pm1$. When converting these symmetries on the ``dual" generators, we have some flexibility on choosing renormalization scalars such that the braid group actions of $\tU$ are bar invariant; cf. \cite[\S5]{BG17a}.
	
	\begin{proposition} {\rm(cf. \cite[\S5]{BG17a})}
		\label{prop:BG1U}
		For $i\in \I$, there exist automorphisms $\widetilde{T}_{i,-e}''$ on $\tU$ such that
		\begin{align*}
			&\widetilde{T}_{i,-e}''(K_\mu)= K_{s_i(\mu)},
			\qquad \widetilde{T}_{i,-e}''(K'_\mu)= K'_{s_i(\mu)},\;\;\forall \mu\in \Z^\I,\\
			&\widetilde{T}_{i,-1}''(E_i)=v_iF_iK_i^{-1},\qquad \widetilde{T}_{i,-1}''(F_i)=v_i^{-1} (K_i')^{-1}E_i,\\
			&\widetilde{T}_{i,1}''(E_i)=v_i^{-1}F_i(K_i')^{-1},\qquad \widetilde{T}_{i,1}''(F_i)=v_i K_i^{-1}E_i,\\
			&\widetilde{T}_{i,-e}''(E_j)
			=\sum_{r+s=-c_{ij}} (-1)^{r}v_i^{e(r+\frac{1}{2}c_{ij})}(v_i-v_i^{-1})^{c_{ij}} E_i^{(s)}  E_j E_i^{(r)}, \quad\forall j\neq i,
			\\
			&\widetilde{T}_{i,-e}''(F_j)
			= \sum_{r+s=-c_{ij}} (-1)^{r}v_i^{e(r+\frac{1}{2}c_{ij})}(v_i-v_i^{-1})^{c_{ij}} F_{i}^{(s)} F_j F_i^{(r)}, \quad\forall j\neq i.
		\end{align*}
	\end{proposition}
	
	\begin{proof}
		For fixed $e=\pm 1$, consider the elements (compare the notations in \eqref{eq:Udj-gen}): 
		\begin{align*}
			\mathcal{E}_i=\frac{v_i^{-\frac{1}{2}e}}{v_i-v_i^{-1}}E_i,\quad \mathcal{F}_i=\frac{v_i^{\frac{1}{2}e}}{v_i^{-1}-v_i}F_i,
		\end{align*}
		Then $\mathcal{E}_i,\mathcal{F}_i,K_i,K_i'$ ($i\in\I$) satisfy the usual defining relations of $\tU$ given in \cite{LW22a}. We then have braid group symmetries $\widetilde{T}_{i,-e}''$ ($i\in\I$) acting by
		\begin{align*}
			&\widetilde{T}_{i,-e}''(K_\mu)= K_{s_i(\mu)},
			\qquad \widetilde{T}_{i,-e}''(K'_\mu)= K'_{s_i(\mu)},\;\;\forall \mu\in \Z^\I,\\
			&\widetilde{T}_{i,-1}''(\mathcal{E}_i)=-\mathcal{F}_iK_i^{-1},\qquad \widetilde{T}_{i,-1}''(\mathcal{F}_i)=-(K_i')^{-1}\mathcal{E}_i,\\
			&\widetilde{T}_{i,1}''(\mathcal{E}_i)=-\mathcal{F}_i(K_i')^{-1},\qquad \widetilde{T}_{i,1}''(\mathcal{F}_i)=-K_i^{-1}\mathcal{E}_i,\\
			&\widetilde{T}_{i,-e}''(\mathcal{E}_j)
			=\sum_{r+s=-c_{ij}} (-1)^{r}v_i^{er} \mathcal{E}_i^{(s)}  \mathcal{E}_j \mathcal{E}_i^{(r)}, \quad\forall j\neq i,
			\\
			&\widetilde{T}_{i,-e}''(\mathcal{F}_j)
			= \sum_{r+s=-c_{ij}} (-1)^{r}v_i^{-er} \mathcal{F}_{i}^{(r)} \mathcal{F}_j \mathcal{F}_i^{(s)}, \quad\forall j\neq i.
		\end{align*}
		Writing $\widetilde{T}_{i,-e}''$ in terms of the generators $E_i,F_i,K_i,K_i'$ then gives the desired formulas: for example,
		\begin{align*}
			\widetilde{T}_{i,-e}''(E_j)&=v_j^{\frac{1}{2}e}(v_j-v_j^{-1})\widetilde{T}_{i,-e}''(\mathcal{E}_j)\\
			&=v_j^{\frac{1}{2}e}(v_j-v_j^{-1})\sum_{r+s=-c_{ij}}(-1)^rv_i^{er}\mathcal{E}_i^{(s)} \mathcal{E}_j \mathcal{E}_i^{(r)}\\
			&=v_i^{\frac{1}{2}ec_{ij}}(v_i-v_i^{-1})^{c_{ij}}\sum_{r+s=-c_{ij}}(-1)^{r}v_i^{er}E_i^{(s)} E_j E_i^{(r)}\\
			&=(v_i-v_i^{-1})^{c_{ij}}\sum_{r+s=-c_{ij}}(-1)^{r}v_i^{e(r+\frac{1}{2}c_{ij})}E_i^{(s)}E_jE_i^{(r)};
		\end{align*}
		and
		\begin{align*}
			\widetilde{T}_{i,-e}''(F_j)&=v_j^{-\frac{1}{2}e}(v_j^{-1}-v_j)\widetilde{T}_{i,-e}''(\mathcal{F}_i)\\
			&=v_j^{-\frac{1}{2}e}(v_j^{-1}-v_j)\sum_{r+s=-c_{ij}} (-1)^{r}v_i^{-er} \mathcal{F}_{i}^{(r)} \mathcal{F}_j \mathcal{F}_i^{(s)}\\
			&=v_i^{-\frac{1}{2}ec_{ij}}(v_i^{-1}-v_i)^{c_{ij}}\sum_{r+s=-c_{ij}} (-1)^{r}v_i^{-er}F_i^{(r)}F_jF_i^{(s)}\\
			&=(v_i-v_i^{-1})^{c_{ij}}\sum_{r+s=-c_{ij}}(-1)^{r+c_{ij}}v_i^{-e(r+\frac{1}{2}c_{ij})}F_i^{(r)}F_jF_i^{(s)}\\
			&=(v_i-v_i^{-1})^{c_{ij}}\sum_{r+s=-c_{ij}}(-1)^{s}v_i^{e(s+\frac{1}{2}c_{ij})}F_i^{(r)}F_jF_i^{(s)}\\
			&=(v_i-v_i^{-1})^{c_{ij}}\sum_{r+s=-c_{ij}}(-1)^{r}v_i^{e(r+\frac{1}{2}c_{ij})}F_i^{(s)}F_jF_i^{(r)}.
			\qedhere
		\end{align*}
	\end{proof}
	
	\begin{proposition}
		\label{prop:BG2U}
		For $i\in \I$, the automorphisms $\widetilde{T}_{i,e}'$ on $\tU$ given by $\widetilde{T}_{i,e}'=\sigma\circ\widetilde{T}_{i,-e}''\circ\sigma$ satisfy that
		\begin{align*}
			&\widetilde{T}_{i,e}'(K_\mu)= K_{s_i(\mu)},
			\qquad \widetilde{T}_{i,e}'(K'_\mu)= K'_{s_i(\mu)},\;\;\forall \mu\in \Z^\I,\\
			&\widetilde{T}_{i,1}'(E_i)=v_i (K_i')^{-1}F_i,\qquad \widetilde{T}_{i,1}'(F_i)=v_i^{-1} E_iK_i^{-1},\\
			&\widetilde{T}_{i,-1}'(E_i)=v_i^{-1} K_i^{-1}F_i,\qquad \widetilde{T}_{i,-1}'(F_i)=v_i E_i(K_i')^{-1},\\
			&\widetilde{T}_{i,e}'(E_j)
			=\sum_{r+s=-c_{ij}} (-1)^{r}v_i^{e(r+\frac{1}{2}c_{ij})}(v_i-v_i^{-1})^{c_{ij}} E_i^{(r)}  E_j E_i^{(s)}, \quad\forall j\neq i,
			\\
			&\widetilde{T}_{i,e}'(F_j)
			= \sum_{r+s=-c_{ij}} (-1)^{r}v_i^{e(r+\frac{1}{2}c_{ij})}(v_i-v_i^{-1})^{c_{ij}} F_i^{(r)}  F_j F_i^{(s)}, \quad\forall j\neq i.
		\end{align*}
	\end{proposition}
	
	\begin{proof}
		This follows from Proposition \ref{prop:BG1U} by a direct computation.
	\end{proof}
	
	Recall the bar involution on $\tU$ from Lemma \ref{QG bar-involution def}. 
	\begin{lemma}\label{lem:QGbraid-bar}
		The braid group actions $\widetilde{T}_{i,e}'$ and $\widetilde{T}_{i,-e}''$ commute with the bar-involution, i.e., $\ov{\widetilde{T}_{i,e}'(u)}=\widetilde{T}_{i,e}'(\ov{u})$ and $\ov{\widetilde{T}_{i,-e}''(u)}=\widetilde{T}_{i,-e}''(\ov{u})$ for any $u\in\tU$.
	\end{lemma}
	
	\begin{proof}
		Define $\ov{\widetilde{T}_{i,e}'}$ and $\ov{\widetilde{T}_{i,-e}''}$ by $\ov{\widetilde{T}_{i,e}'}(u):=\ov{\widetilde{T}_{i,e}'(\ov{u})}$ and $\ov{\widetilde{T}_{i,-e}''}(u):=\ov{\widetilde{T}_{i,-e}''(\ov{u})}$ for $u\in\tU$. Then $\ov{\widetilde{T}_{i,e}'}$ and $\ov{\widetilde{T}_{i,-e}''}$ are homomorphisms of algebras. In order to prove that $\ov{\widetilde{T}_{i,e}'}=\widetilde{T}_{i,e}'$ and $\ov{\widetilde{T}_{i,-e}''}=\widetilde{T}_{i,-e}''$, we only need to check the identities on generators $E_j,F_j,K_j,K_j'$, which follow by a direct computation.
	\end{proof}
	
	We shall often use the shorthand notation 
	\begin{align}
		\label{braid-shorthand}
		\widetilde{T}_i:=\widetilde{T}_{i,1}'',\qquad \widetilde{T}_i^{-1}:=\widetilde{T}_{i,-1}'.
	\end{align}
	The $\widetilde{T}_i$'s satisfy the braid group relations and so 
	$\widetilde{T}_w 
	:= \widetilde{T}_{i_1}\cdots
	\widetilde{T}_{i_r} \in \Aut(\tU)$ is well defined,
	where $w = s_{i_1}\cdots s_{i_r}$ is any reduced expression of $w \in W$.

	\subsection{iQuantum groups}
	
	For a Cartan matrix $C=(c_{ij})_{i,j\in \I}$, let $\text{Inv}(C)$ be the group of permutations $\btau$ of the set $\I$ such that $c_{ij}=c_{\btau i,\btau j}$, for all $i,j$, and $\btau^2=\Id$. Then  $\btau \in \text{Inv}(C)$ can be viewed as an involution (which is allowed to be the identity) of the corresponding Dynkin diagram (which is identified with $\I$ by abuse of notation). We shall refer to the pair $(\I,\btau)$ as a (quasi-split) Satake diagram.

	We denote by $\bs_{i}$ the following element of order 2 in the Weyl group $W$, i.e.,
	\begin{align}
		\label{def:ri}
		\bs_i= \left\{
		\begin{array}{ll}
			s_{i}, & \text{ if } c_{i,\tau i}=2,
			\\
			s_is_{\btau i}, & \text{ if } c_{i,\tau i}=0,
			\\
			s_is_{\btau i}s_i, & \text{ if } c_{i,\tau i}=-1.
		\end{array}
		\right.
	\end{align}
	It is well known that the {\rm restricted Weyl group} associated to $(\I,\tau)$ can be identified with the following subgroup $W_\btau$ of $W$:
	\begin{align}
		\label{eq:Wtau}
		W_{\btau} =\{w\in W\mid \btau w =w \btau\},
	\end{align}
	where $\btau$ is regarded as an automorphism of the root lattice $\Z^\I$. Moreover, the restricted Weyl group $W_{\btau}$ can be identified with a Weyl group with $\bs_i$ ($i\in \I_\btau$) as its simple reflections. 
	
	Associated with the Satake diagram $(\I,\tau)$, following \cite{LW22a} we define the universal iquantum groups ${\hU}^\imath$ (resp. $\tUi$) to be the $\Q(v^{\frac12})$-subalgebra of $\hU$ (resp. $\tU$) generated by
	\begin{equation}
		\label{eq:Bi}
		B_i= F_i +  E_{\btau i} \tK_i',
		\qquad \tk_i = \tK_i \tK_{\btau i}', 
		\quad \forall i \in \I, 
	\end{equation}
	(with $\tk_i$ invertible in $\tUi$). Let $\hU^{\imath 0}$ be the $\Q(v^{\frac12})$-subalgebra of $\hUi$ generated by $\tk_i$, for $i\in \I$. Similarly, let $\tU^{\imath 0}$ be the $\Q(v^{\frac12})$-subalgebra of $\tUi$ generated by $\tk_i^{\pm1}$, for $i\in \I$. 
	The algebra $\widetilde{\bU}^\imath$ (resp. $\hUi$) is a right coideal subalgebra of $\widetilde{\bU}$ (resp. $\hU$); the pairs $(\widetilde{\bU}, \widetilde{\bU}^\imath)$ and $(\hU,\hUi)$ are called quantum symmetric pairs, and $\hUi$ and $\tUi$ are called the universal {\em (quasi-split) iquantum groups}; they are {\em split} if $\btau =\Id$.
	
	Let $\bvs=(\vs_i)\in(\Q(v^{\frac12})^\times)^\I$ be such that $\vs_i=\vs_{\tau i}$ for each $i\in\I$ which satisfies $c_{i,\tau i}=0$. The iquantum groups \`a la Letzter-Kolb \cite{Let99,Ko14} $\Ui=\Ui_\bvs$ is the $\Q(v^{\frac12})$-subalgebra of $\U$ generated by
	$$B_i=F_i+\vs_i E_{\tau i}K_i^{-1},\quad k_i=K_iK_{\tau i}^{-1},\quad \forall i\in\I.$$ 
	By \cite[Proposition 6.2]{LW22a}, the $\Q(v^{\frac12})$-algebra $\Ui$ is isomorphic to the quotient of $\tUi$ by the ideal generated by $\tk_i-\vs_i$ (for $i=\tau i$) and $\tk_i\tk_{\tau i}-\vs_i\vs_{\tau i}$ (for $i\neq \tau i$).
	
	For $i\in \I$ with $\btau i= i$, adapting the definitions in \cite{BW18b, BeW18, CLW21a} to our $\tU$ with ``dual" generators, we define the {\em idivided powers} of $B_i$ to be  
	\begin{align}
		B_{i,\odd}^{(m)}&=\frac{1}{[m]_{v_i}^!} \begin{cases} B_i\prod_{s=1}^k (B_i^2+v_i(v_i-v_i^{-1})^2\tk_i[2s-1]_{v_i}^2 ) & \text{if }m=2k+1,\\
			\quad\prod_{s=1}^k (B_i^2+v_i(v_i-v_i^{-1})^2\tk_i[2s-1]_{v_i}^2) &\text{if }m=2k; \end{cases}
		\\
		B_{i,\ev}^{(m)}= &\frac{1}{[m]_{v_i}^!} \begin{cases} B_i\prod_{s=1}^k (B_i^2+v_i(v_i-v_i^{-1})^2\tk_i[2s]_{v_i}^2 ) & \text{if }m=2k+1,\\
			\quad\prod_{s=1}^{k} (B_i^2+v_i(v_i-v_i^{-1})^2\tk_i[2s-2]_{v_i}^2) &\text{if }m=2k. \end{cases}
	\end{align}
	On the other hand, for $i\in \I$ with $i \neq \tau i$, we define the divided powers as usual: for $m\in \N$,
	\begin{align}
		\label{eq:DP}
		B_{i}^{(m)}= \frac{B_i^m}{[m]_{v_i}^!}.
	\end{align}
	
	Denote 
	$$
	(a;x)_0=1,\quad (a;x)_n=(1-a)(1-ax)\cdots(1-ax^{n-1}),\quad \forall n\geq1.
	$$
	A Serre presentation of $\Ui$ was obtained in all finite type by Letzter \cite{Let02} and in certain Kac-Moody type \cite{Ko14}; for $\Ui$ of arbitrary quasi-split type a Serre relation was given in \cite[Theorem 3.1]{CLW21a} which contains the new relation \eqref{relation6} and incorporates the relation \eqref{relation5} from \cite[Theorem 3.6]{BK15}, and this presentation has been adapted to $\tUi$ for general quasi-split type in \cite[Proposition~ 6.4]{LW22a}. We reformulate it below for $\tUi$ in terms ``dual" generators (the presentation also holds for $\hUi$ by omitting that $\tk_i$ ($i\in\I$) are invertible).
	
	\begin{proposition} [cf. \text{\cite[Theorem 3.1]{CLW21a}}]
		\label{prop:Serre}
		Fix $\ov{p}_i\in \Z_2$ for each $i\in \I$.
		The $\Q(v^{\frac12})$-algebra $\tUi$ has a presentation with generators $B_i$, $\tk_i$ $(i\in \I)$, where $\tk_i$ are invertible, subject to the relations \eqref{relation1}--\eqref{relation6} below: for $\ell \in \I$, and $i\neq j \in \I$,
		\begin{align}
			\tk_i \tk_\ell =\tk_\ell \tk_i,
			&\qquad
			\tk_i B_\ell  = v_i^{c_{\btau i,\ell} -c_{i \ell}} B_\ell \tk_i,
			\label{relation1}
			\\
			B_iB_{j}-B_jB_i &=0, \quad \text{ if }c_{ij} =0 \text{ and }\btau i\neq j,\label{relation2}
			\\
			\sum_{r=0}^{1-c_{ij}} (-1)^rB_i^{(r)}&B_jB_i^{(1-c_{ij}-r)} =0, \quad \text{ if } j \neq i\neq \btau i, \label{relation3}
			\\
			\sum_{r=0}^{1-c_{i,\btau i}} (-1)^{r+c_{i,\btau i}}&B_i^{(r)}B_{\btau i}B_i^{(1-c_{i,\btau i}-r)} =(v_i^{-1}-v_i)\times
			\notag
			\\
			\Big(v^{c_{i,\btau i}} (&v_i^{-2};v_i^{-2})_{-c_{i,\btau i}}    
			B_i^{(-c_{i,\btau i})} \tk_i
			-(v_i^{2};v_i^{2})_{-c_{i,\btau i}}B_i^{(-c_{i,\tau i})} \tk_{\btau i}  \Big),
			\text{ if } \btau i \neq i, 
			\label{relation5} \\
			\sum_{r=0}^{1-c_{ij}} (-1)^r  B_{i, \overline{p_i}}^{(r)}&B_j B_{i,\overline{c_{ij}}+\overline{p}_i}^{(1-c_{ij}-r)} =0,\quad   \text{ if }i=\btau i.
			\label{relation6}
		\end{align}
	\end{proposition}
	
	
	For $i\in\I$, for any $\alpha=\sum_{i\in\I} a_i\alpha_i\in\Z^\I$, we set
	\begin{align}
		\label{eq:bbKi}
		\K_i=v_i^{\frac{1}{2}c_{i,\tau i}}\tk_i,\qquad \K_\alpha:=\prod_{i\in\I}\K_i^{a_i}.
	\end{align}
	
	We formulate some basic symmetries of $\hUi$.
	
	\begin{lemma}  
		\label{lem:involution-iQG}
		\quad
		\begin{enumerate}
			\item 
			There exists an anti-involution $\sigma^\imath$ on $\hUi$ (and also $\tUi$) given by $\sigma^\imath(B_i)=B_i$, $\sigma^\imath(\tk_i)=\tk_{\tau i}$, for $i\in\I$.
			\item There exists an anti-involution (called bar-involution) $:u\mapsto \ov{u}$ on $\hUi$ (and also $\tUi$) given by $\ov{v^{1/2}}=v^{-1/2}$, $\ov{B_i}=B_i$, and $\ov{\K_i}=\K_i$,  for $i\in\I$. In particular, $\ov{\tk_i}=v_i^{c_{i,\tau i}}\tk_{i}$. 
			\item There exists an involution $\psi^\imath$ of $\tUi$ such that $\psi^\imath(v^{1/2})=v^{-1/2}$, $\psi^\imath(B_i)=B_i$, $\psi^\imath(\tk_i)=v_i^{c_{i,\tau i}}\tk_{\tau i}$, for $i\in\I$.
		\end{enumerate}
	\end{lemma}
	
	\begin{proof}
		Parts (1) and (3) can be found in \cite[Lemma 6.9]{LW22b}, while (2) follows from Proposition \ref{prop:Serre} by a direct computation.
	\end{proof}

	\begin{example} 
		Quasi-split iquantum groups of finite type are associated to the Satake diagrams of finite type, i.e., Dynkin graphs with involution. The split Satake diagrams formally look the same as Dynkin diagrams. A complete list of such connected quasi-split Satake diagrams of finite type with $\tau \neq \id$ is given in Table \ref{tab:Satakediag} below.  
		\begin{table}[h]  
			\begin{center}
				\centering
				\begin{tabular}{|m{5cm}<{\centering}|m{8cm}<{\centering}|}
					\hline
					Types &  Satake diagrams  \\
					\hline

					${\rm AIII}_{2r-1}$ $(r\geq 1)$ & \setlength{\unitlength}{0.6mm}			\begin{picture}(70,35)(0,5)
						
						\put(0,10){$\circ$}
						\put(0,30){$\circ$}
						\put(50,10){$\circ$}
						\put(50,30){$\circ$}
						\put(72,10){$\circ$}
						\put(72,30){$\circ$}
						\put(92,20){$\circ$}
						\put(-5,5.5){$2r-1$}
						\put(-2,34){${1}$}
						\put(47,6){\small $r+2$}
						\put(47,34){\small $r-2$}
						\put(69,6){\small $r+1$}
						\put(69,34){\small $r-1$}
						\put(92,16){\small $r$}
						
						\put(3,11.5){\line(1,0){16}}
						\put(3,31.5){\line(1,0){16}}
						\put(23,10){$\cdots$}
						\put(23,30){$\cdots$}
						\put(33.5,11.5){\line(1,0){16}}
						\put(33.5,31.5){\line(1,0){16}}
						\put(53,11.5){\line(1,0){18.5}}
						\put(53,31.5){\line(1,0){18.5}}
						
						\put(75,12){\line(2,1){17}}
						\put(75,31){\line(2,-1){17}}

						\color{red}
						\put(-7,20){$\tau$}
						\qbezier(0,13.5)(-4,21.5)(0,29.5)
						\put(-0.25,14){\vector(1,-2){0.5}}
						\put(-0.25,29){\vector(1,2){0.5}}
						
						\qbezier(50,13.5)(46,21.5)(50,29.5)
						\put(49.75,14){\vector(1,-2){0.5}}
						\put(49.75,29){\vector(1,2){0.5}}

						\qbezier(72,13.5)(68,21.5)(72,29.5)
						\put(71.75,14){\vector(1,-2){0.5}}
						\put(71.75,29){\vector(1,2){0.5}}
					\end{picture}\\
					\hline
					${\rm AIII}_{2r}$ $(r\geq 1)$ & \setlength{\unitlength}{0.6mm}			\begin{picture}(70,35)(0,5)
						
						\put(0,10){$\circ$}
						\put(0,30){$\circ$}
						\put(50,10){$\circ$}
						\put(50,30){$\circ$}
						\put(72,10){$\circ$}
						\put(72,30){$\circ$}
						\put(-4,5){$2r$}
						\put(-2,34){${1}$}
						\put(47,6){\small $r+2$}
						\put(47,34){\small $r-1$}
						\put(69,6){\small $r+1$}
						\put(69,34){\small $r$}
						
						\put(3,11.5){\line(1,0){16}}
						\put(3,31.5){\line(1,0){16}}
						\put(23,10){$\cdots$}
						\put(23,30){$\cdots$}
						\put(33.5,11.5){\line(1,0){16}}
						\put(33.5,31.5){\line(1,0){16}}
						\put(53,11.5){\line(1,0){18.5}}
						\put(53,31.5){\line(1,0){18.5}}
						
						\put(73.5,13.6){\line(0,1){16}}
						
						\color{red}
						\put(-7,20){$\tau$}
						\qbezier(0,13.5)(-4,21.5)(0,29.5)
						\put(-0.25,14){\vector(1,-2){0.5}}
						\put(-0.25,29){\vector(1,2){0.5}}
						
						\qbezier(50,13.5)(46,21.5)(50,29.5)
						\put(49.75,14){\vector(1,-2){0.5}}
						\put(49.75,29){\vector(1,2){0.5}}

						\qbezier(76,13.5)(80,21.5)(76,29.5)
						\put(75.85,14){\vector(-1,-2){0.5}}
						\put(75.85,29){\vector(-1,2){0.5}}
					\end{picture}\\
					\hline
					${\rm DI}_r$ $(r\geq4)$ & \setlength{\unitlength}{0.6mm}	\begin{picture}(40,35)(20,-15)
						\put(0,-1){$\circ$}
						\put(0,-6){\small$1$}
						
						\put(3,0){\line(1,0){16.5}}
						\put(20,-1){$\circ$}
						\put(20,-6){\small$2$}
						\put(64,-1){$\circ$}
						\put(56,-6){\small$r-2$}
						\put(84,-10){$\circ$}
						\put(80,-13){\small${r-1}$}
						\put(84,9.5){$\circ$}
						\put(84,13.5){\small${r}$}

						\put(38,0){\line(-1,0){15.5}}
						\put(64,0){\line(-1,0){15}}
						
						\put(40,-1){$\cdots$}
						
						\put(83.5,9.5){\line(-2,-1){16.5}}
						\put(83.5,-8.5){\line(-2,1){16.5}}
						
						\put(12,-20.5){\begin{picture}(100,40)	\color{red}
								\put(79,20){$\tau$}
								\qbezier(75,13.5)(79,21.5)(75,29.5)
								\put(75.25,14){\vector(-1,-2){0.5}}
								\put(75.25,29){\vector(-1,2){0.5}}
						\end{picture}}
					\end{picture}
					\\
					\hline
					${\rm EII}_6$  &	\setlength{\unitlength}{0.6mm}		\begin{picture}(70,35)(20,7)
						\put(97,36){\small${6}$}
						\put(75,36){\small${5}$}
						\put(97,6.5){\small${1}$}
						\put(75,6.5){\small${2}$}
						\put(33,16){\small $4$}	
						\put(55,16){\small $3$}	\put(10,35){\rotatebox[origin=c]{180}{\begin{picture}(100,10)
									
									\put(10,10){$\circ$}
									
									\put(32,10){$\circ$}
									
									\put(10,30){$\circ$}
									
									\put(32,30){$\circ$}

									
									\put(31.5,11){\line(-1,0){19}}
									\put(31.5,31){\line(-1,0){19}}
									
									\put(52,22){\line(-2,1){17.5}}
									\put(52,20){\line(-2,-1){17.5}}
									
									\put(54.7,21.2){\line(1,0){19}}

									\put(52,20){$\circ$}
									
									\put(74,20){$\circ$}
							\end{picture}}

							\put(-37,-32.5){\begin{picture}(100,40)	\color{red}
									\qbezier(5,13.5)(9,21.5)(5,29.5)
									\put(5.25,14){\vector(-1,-2){0.5}}
									\put(5.25,29){\vector(-1,2){0.5}}
							\end{picture}}
							\put(-15,-32.5){\begin{picture}(100,40)	\color{red}
									\put(9,20){$\tau$}
									\qbezier(5,13.5)(9,21.5)(5,29.5)
									\put(5.25,14){\vector(-1,-2){0.5}}
									\put(5.25,29){\vector(-1,2){0.5}}
							\end{picture}}
						}
						
					\end{picture}
					\\\hline
				\end{tabular} 
			\end{center}
			\vspace{0.5cm}
			\caption{Quasi-split Satake diagrams of finite type with $\tau \neq {\rm Id}$}
			\label{tab:Satakediag}
		\end{table}
	\end{example}
	
	\begin{example}
		\label{ex:QGvsiQG}
		{\rm (Quantum groups as iquantum groups of diagonal type)} 
		Consider the $\Q(v^{\frac12})$-subalgebra $\tUUi$ of $\tUU$
		generated by
		\[
		\ck_i:=\tK_{i} \tK_{i^{\diamond}}', \quad
		\ck_i':=\tK_{i^{\diamond}} \tK_{i}',  \quad
		\cb_{i}:= F_{i}+ E_{i^{\diamond}} \tK_{i}', \quad
		\cb_{i^{\diamond}}:=F_{i^{\diamond}}+ E_{i} \tK_{i^{\diamond}}',
		\qquad \forall i\in \I.
		\]
		Here we drop the tensor product notation and use instead $i^\diamond$ to index the generators of the second copy of $\tU$ in $\tUU$. There exists a $\Q(v^{\frac12})$-algebra isomorphism $\widetilde{\phi}: \tU \rightarrow \tUUi$ such that
		\[
		\widetilde{\phi}(E_i)= \cb_{i},\quad \widetilde{\phi}(F_i)= \cb_{i^{\diamond}}, \quad \widetilde{\phi}(\tK_i)= \ck_i', \quad \widetilde{\phi}(\tK_i')= \ck_i, \qquad \forall  i\in \I.
		\]
		In this case, the Satake diagram is $(\Gamma\sqcup\Gamma^\diamond,\swa)$, where $\Gamma^\diamond$ is a copy of the Dynkin diagram $\Gamma$ of $\tU$. 
	\end{example}

	\subsection{Relative braid group symmetries}
	
	Choose one representative for each $\btau$-orbit on $\I$, and let
	\begin{align}\label{eq:ci}
		\ci = \{ \text{the chosen representatives of $\btau$-orbits in $\I$} \}.
	\end{align} 
	The braid group associated to the relative Weyl group $W_\tau$ is denoted
	\begin{equation}
		\label{eq:braidCox}
		\brW =\langle \br_i \mid i\in \I_\btau \rangle,
	\end{equation}
	where $\br_i$ satisfy the same braid relations as for $\bs_i$ in $W_{\tau }$. The relative braid (or ibraid) group symmetries $\tTT_{i,e}'$ and $\tTT''_{i,e}$ ($i\in\I$, $e\in\{+1,-1\}$) on $\tUi$ are established in \cite{LW22a, WZ23} (and \cite{Z23}); see \cite{KP11} for earlier conjectures on iquantum groups with specific parameters. In this paper, we shall also use the bar-equivariant versions of these ibraid group symmetries of $\tUi$. 
	
	Quasi K-matrix appeared earlier in different formulations; see \cite{BW18a, BK19,AV22}. We shall need the following.
	
	\begin{proposition}
		[\text{\cite[Theorem~3.6]{WZ23}}]
		\label{prop:Kmatrix}
		There exists a unique element $\widetilde{\Upsilon}=\sum_{\mu\in\N^\I}\widetilde{\Upsilon}^\mu$ (called quasi K-matrix) such that $\widetilde{\Upsilon}^0=1$, $\widetilde{\Upsilon}^\mu\in\U^+_{\mu}$ and the following identities hold:
		\begin{align}
			B_i\widetilde{\Upsilon}&=\widetilde{\Upsilon}B_i^\sigma,\qquad (i\in\I),
			\\
			x\widetilde{\Upsilon}&=\widetilde{\Upsilon}x,\qquad (x\in\tU^{\imath0}),
		\end{align}
		where $B_i^\sigma:=\sigma(B_i)=F_i+K_iE_{\tau i}$.  Moreover, $\widetilde{\Upsilon}^\mu=0$ unless $\tau(\mu)=\mu$.
	\end{proposition}
	
	Denote by $\tU_{i,\tau i}$ the quantum group associated to $\I_i=\{i,\tau i\}$. Let $\widetilde{\Upsilon}_i$ be the rank one quasi K-matrix associated to $\I_i=\{i,\tau i\}$, i.e., $\widetilde{\Upsilon}_i=\sum_{\mu\in\N^\I}\widetilde{\Upsilon}_{i}^{\mu}$ with $\widetilde{\Upsilon}_{i}^{\mu}\in\tU_{i,\tau i}^+$, and $\widetilde{\Upsilon}_{i}^{0}=1$. Define a distinguished parameter $\bm{\varsigma}_\diamond=(\varsigma_{i,\diamond})_{i\in\I}$ by
	\begin{align}
		\label{eq:disting-para}
		\varsigma_{i,\diamond}=v^{-\frac{1}{2}(\alpha_i,\alpha_{\tau i})}.
	\end{align}
	Recall the automorphism $\widetilde{\Psi}_{\varsigma_{\diamond}}$ of $\tU$ from Lemma~\ref{lem:twisting Psi}. We set
	\begin{equation}\label{eq: T_i twisted def}
		\widetilde{\mathscr{T}}_i:=\widetilde{\Psi}_{\bm{\varsigma}_\diamond}^{-1}\circ\widetilde{T}_i\circ \widetilde{\Psi}_{\bm{\varsigma}_\diamond},
		\qquad \widetilde{\mathscr{T}}_i^{-1}:=\widetilde{\Psi}_{\bm{\varsigma}_\diamond}^{-1}\circ\widetilde{T}_i^{-1}\circ \widetilde{\Psi}_{\bm{\varsigma}_\diamond}.
	\end{equation}
	Clearly $\widetilde{\mathscr{T}}_i$ and $\widetilde{\mathscr{T}}_i^{-1}$, for $i\in\I$, are automorphisms of $\tU$ and satisfy the braid group relations. Hence, we can define
	$\widetilde{\mathscr{T}}_w:=\widetilde{\mathscr{T}}_{i_1}\cdots \widetilde{\mathscr{T}}_{i_r},$
	where $w=s_{i_1}\cdots s_{i_r}$ is any reduced expression. 
	
	\begin{theorem}[{cf. \cite{WZ23,Z23}}]
		\label{thm: relative T_i conjugate}
		For $i\in\I$, there are mutually inverse automorphisms $\widetilde{\TT}_i$ and $\widetilde{\TT}_i^{-1}$ on $\tUi$ such that
		\begin{align}
			\label{eq:relbraid1}
			\widetilde{\TT}_i^{-1}(x)\widetilde{\Upsilon}_i &=\widetilde{\Upsilon}_i\widetilde{\mathscr{T}}_{r_i}^{-1}(x),
			\\
			\label{eq:relbraid2}
			\widetilde{\TT}_i(x)\widetilde{\mathscr{T}}_i^{-1}(\widetilde{\Upsilon}_i)^{-1} &=\widetilde{\mathscr{T}}_i^{-1}(\widetilde{\Upsilon}_i)^{-1}\widetilde{\mathscr{T}}_{r_i}(x).
		\end{align}
		Moreover, we have $\widetilde{\TT}_i^{-1}=\sigma^\imath\circ \widetilde{\TT}_i\circ \sigma^\imath$, and there exists a group homomorphism $\Br(W_\tau)\rightarrow \Aut(\tUi)$, $\br_i\mapsto \tTT_i$ for $i\in\I$.
	\end{theorem}
	
	\begin{proof}
		While the proof is a rerun of the arguments given in \cite{Z23} and \cite{WZ23}, it does not formally reduce to them because of renormalization here, and so let us outline it. 
		
		In view of the proof of \cite[Theorem 3.1]{Z23}, it suffices to prove the existence of elements $\widetilde{\TT}_i^{-1}(x)$ and $\widetilde{\TT}_i(x)$ satisfying \eqref{eq:relbraid1} and \eqref{eq:relbraid2} on generators of $\tUi$. The case for $x=B_i,B_{\tau i}$ or $x\in\tU^{\imath 0}$ is straightforward (cf. \cite[\S4.5, 4.7]{WZ23}). For $x=B_j$, $j\notin\I_i$, the proof of \cite[Theorems 4.11, 5.17, 6.14]{Z23} can be adapted to our case. Note that the definition of $\widetilde{\mathscr{T}}_i$ is different from \cite{WZ23,Z23} because we are working with root vectors in $\U^+$, while they considered root vectors in $\U^-$. The recursive relations used in \cite{Z23} should also be altered accordingly (cf. Lemmas~\ref{lem: f_m induction relation}, \ref{lem: f_m1m2 induction relation} and \ref{lem: f_abc induction relation}), and the results of \cite[Subsections 4.3, 5.3, 6.2]{Z23} remain valid. The identity $\widetilde{\TT}_i^{-1}=\sigma^\imath\circ \widetilde{\TT}_i\circ \sigma^\imath$ and the relative braid relations for $\widetilde{\TT}_i$ can be checked as in \cite{WZ23,Z23}.
	\end{proof}
	
	Denote by $\tau_i$ the diagram involution of $\I_i:=\{i,\tau i\}$ defined by 
	\begin{align} \label{taui}
		r_i(\alpha_i)=-\alpha_{\tau_i(i)}, 
		\qquad
		r_i(\alpha_{\tau i})=-\alpha_{\tau_i(\tau i)}.
	\end{align}
	
	\begin{proposition} [{\cite[Proposition 4.11, Theorem 4.14]{WZ23}}]
		\label{prop:TiBi}
		For $i,j\in\I$, we have $\widetilde{\TT}_i(\K_j)=\K_{r_i(\alpha_j)}$ and
		\[\widetilde{\TT}_i(B_i)=v^{\frac{1}{2}(\alpha_i-\alpha_{\tau i},\alpha_i)}\K_{\tau_i(i)}^{-1}B_{\tau_i(\tau i)},\quad \widetilde{\TT}_i(B_{\tau i})=v^{\frac{1}{2}(\alpha_i-\alpha_{\tau i},\alpha_i)}\K_{\tau_i(\tau i)}^{-1}B_{\tau_i(i)}.
		\]
	\end{proposition}
	
	\begin{proof}
		We give a direct proof in our greatly simplified quasi-split setting. 
		From the definition of $\K_\alpha$ \eqref{eq:bbKi}, \cite[Proposition 4.11]{WZ23} implies that $\widetilde{\TT}_i(\K_j)=\K_{r_i(\alpha_j)}$. The formula of $\widetilde{\TT}_i(B_i)$ follows from a direct computation. In fact,  we have
		\begin{align*}
			\widetilde{\mathscr{T}}_{r_i}^{-1}(B_i)&=\widetilde{\mathscr{T}}_{r_i}^{-1}(F_i+E_{\tau i}K_i')\\
			&=v_i^{-1}K'^{-1}_{\tau_i(i)}E_{\tau_i(i)}+v_i^{-(\alpha_i,\alpha_{\tau i})}v_i^{-1}K_{\tau_i(\tau i)}^{-1}F_{\tau_i(\tau i)}K'^{-1}_{\tau_i(i)}\\
			&=v_i^{-1}(K_{\tau_i(\tau i)}K_{\tau_i(i)}')^{-1}(K_{\tau_i(\tau i)}E_{\tau_i(i)}+F_{\tau_i(\tau i)})\\
			&=v^{\frac{1}{2}(\alpha_{\tau i}-\alpha_i,\alpha_i)}\K_{\tau_i(\tau i)}^{-1}B_{\tau_i(\tau i)}^\sigma,
		\end{align*}
		where $B_i^\sigma=F_i+K_iE_{\tau i}$ by definition.
		Since $\widetilde{\Upsilon}_iB_{\tau_i(\tau i)}^\sigma=B_{\tau_i(\tau i)}\widetilde{\Upsilon}_i$ (see Proposition \ref{prop:Kmatrix}), we conclude that the element
		\[\widetilde{\TT}_i^{-1}(B_i)=v^{\frac{1}{2}(\alpha_{\tau i}-\alpha_i,\alpha_i)}\K_{\tau_i(\tau i)}^{-1}B_{\tau_i(\tau i)}\]
		satisfies \eqref{eq:relbraid1}. The formula of $\widetilde{\TT}_i(B_i)$ then follows by applying $\sigma^\imath$; cf. Theorem \ref{thm: relative T_i conjugate}.
	\end{proof}
	
	\begin{theorem}[{\cite[Theorem 3.7]{Z23}}]\label{thm: relative T_i generator}
		Let $i,j\in\I$ be such that $j\notin\{i,\tau i\}$.
		\begin{enumerate}
			\item If $i=\tau i$, then
			\begin{equation}
				\begin{aligned}\label{eq: T_i(B_j) for i=taui}
					\tTT_i(B_j)&=\sum_{r+s=-c_{ij}}(-1)^{r}v_i^{-r-\frac{1}{2}c_{ij}}(v_i-v_i^{-1})^{c_{ij}}B_{i,\ov{p}}^{(s)}B_jB_{i,\ov{c_{ij}}+\ov{p}}^{(r)}\\
					&\quad +\sum_{u\geq 1}\sum_{\stackrel{r+s=-c_{ij}-2u}{\ov{r}=\ov{p}}}(-1)^{r}v_i^{u-r+\frac{1}{2}(-c_{ij}-2u)}(v_i-v_i^{-1})^{c_{ij}+2u}   B_{i,\ov{p}}^{(s)}B_jB_{i,\ov{c_{ij}}+\ov{p}}^{(r)}\tk_i^u.
				\end{aligned}
			\end{equation}
			\item If $c_{i,\tau i}=0$, then
			\begin{equation}
				\begin{aligned}\label{eq: T_i(B_j) for c_itaui=0}
					\tTT_i(B_j)=\sum_{u=0}^{-\max\{c_{ij},c_{\tau i,j}\}}&\sum_{\substack{r_1+s_1=-c_{ij}-u\\ r_2+s_2=-c_{\tau i,j}-u}}(-1)^{r_1+r_2}v_i^{\frac{1}{2}(-c_{ij}-c_{\tau i,j}-2u)}(v_i-v_i^{-1})^{c_{ij}+c_{\tau i,j}+2u}\times\\
					&v_i^{-(r_1+r_2)+(r_1-r_2)u}B_i^{(s_1)}B_{\tau i}^{(s_2)}B_jB_{\tau i}^{(r_2)}B_{i}^{(r_1)}\tilde{k}_{\tau i}^u.
				\end{aligned}
			\end{equation}
			\item If $c_{i,\tau i}=-1$, then
			\begin{equation}
				\begin{aligned}\label{eq: T_i(B_j) for c_itaui=-1}
					\widetilde{\TT}_i(B_j)&=\sum_{w=0}^{-c_{\tau i,j}}\sum_{u=0}^{-c_{ij}}\sum_{\substack{r_1+s_1=-c_{\tau i,j}-w\\ r_2+s_2=-c_{ij}-c_{\tau i,j}-u-w\\ r_3+s_3=-c_{ij}-u}}(-1)^{r_1+r_2+r_3}v_i^{-c_{ij}-c_{\tau i,j}-u-w}(v_i-v_i^{-1})^{2(c_{ij}+c_{\tau i,j}+u+w)}\times\\
					&v_i^{2wr_1-(r_1+r_2+r_3)+(2u-w)r_2-ur_3-uw+\frac{u(u-1)+w(w-1)}{2}}B_i^{(s_1)}B_{\tau i}^{(s_2)}B_i^{(s_3)}B_jB_{i}^{(r_3)}B_{\tau i}^{(r_2)}\tilde{k}_{i}^u B_i^{(r_1)}\tilde{k}_{\tau i}^{w}.
				\end{aligned}
			\end{equation}
		\end{enumerate}
	\end{theorem}
	
	\begin{proof}
		The explicit formulas of $\widetilde{\TT}_i(B_j)$ are given in \cite[Theorem 3.7]{Z23}, the claim then follows from a renormalization according to our conventions on $\widetilde{\mathscr{T}}_{r_i}$.
	\end{proof}
	
	Similar to Lemma \ref{lem:QGbraid-bar}, we have the following.
	
	\begin{lemma}
		\label{lem:bar-invar-braid-Ui}
		The braid group actions $\tTT_{i}$ commute with the bar-involution, i.e., $\ov{\tTT_i(u)}=\tTT_i(\ov{u})$ for any $u\in\tUi$.
	\end{lemma}
	
	Corresponding to Lusztig's braid group symmetries $\widetilde{T}_{i,e}'$, $\widetilde{T}_{i,e}''$ on $\tU$, as in \cite{LW22b,WZ23,Z23}, we define
	\begin{align}
		\tTT_{i,1}''&=\tTT_{i},\quad \tTT_{i,-1}'=\tTT_{i}^{-1},
		\\
		\tTT_{i,-1}''&=\psi^\imath\circ \tTT_i\circ \psi^\imath,
		\quad \tTT_{i,1}'=\psi^\imath\circ \tTT_i^{-1}\circ \psi^\imath.
	\end{align}
	Moreover, we have 
	$$\tTT_{i,e}'=\sigma^\imath\circ \TT_{i,-e}''\circ \sigma^\imath,\quad e\in\{+1,-1\}, i\in\I.$$
	Then all the braid group actions $\tTT_{i,e}',\tTT_{i,e}''$ commute with the bar-involution since $\sigma^\imath,\psi^\imath$ commute with the bar-involution.


	\section{$\mathrm{i}$Quantum groups as $\mathrm{i}$Hopf algebras on the Borel}
	\label{sec:iQG=iHopfB}
	
	In this section, we show that a universal iquantum group can be realized as the iHopf algebra defined on the Borel quantum group.
	
	\subsection{iHopf algebra defined on $\tB$} 
	\label{subsec:iHopf:ff}
	
	Recall the Cartan matrix $C=(c_{ij})$ and $D=\diag(d_i\mid i\in\I)$. Let $'\ff$ be the free associative $\Q(v^{\frac12})$-algebra with generators $\theta_i$ ($i\in\I$); see \cite[Chap. 1]{Lus93}. We denote a rescaled version of $\theta_i$ by 
	\begin{align} \label{rescale:theta}
		\vartheta_i=(v_i -v_i^{-1}) \theta_i,
	\end{align}
	i.e., $\theta_i=(v_i -v_i^{-1}) ^{-1} \vartheta_i$. Let $\ff$ be the quotient algebra of $'\ff$ by the ideal generated by
	\begin{align}\label{Serre rel for ff}
		& \sum_{r=0}^{1-c_{ij}} (-1)^r \left[ \begin{array}{c} 1-c_{ij} \\r \end{array} \right]_{v_i}  \vartheta_i^r \vartheta_j  \vartheta_i^{1-c_{ij}-r},\qquad \forall i\neq j.
	\end{align}
	We endow $\ff$ with an $\N^\I$-grading by setting $\wt(\vartheta_i)=\alpha_i$. Let $\ff_\mu$ be the homogeneous subspace of degree $\mu$. Then $\ff=\bigoplus_{\mu\in\N^\I}\ff_\mu$.
	
	Let $'\hB$ be the $\Q(v^{\frac12})$-algebra generated by
	$\vartheta_i,\h_i$ $(i\in\I)$ subject to 
	\begin{align*}
		[\h_i,\h_j]=0, \qquad h_i\vartheta_j=v_i^{c_{ij}} \vartheta_j h_i.  
	\end{align*}
	and $\hB$ be the quotient algebra of $'\hB$ by the ideal generated by \eqref{Serre rel for ff}. 
	
	Let $r:\ff\to\ff\otimes\ff$ be the homomorphism defined by Lusztig \cite[1.2.6]{Lus93}. In Sweedler notation, we write $r(x)=\sum x_{(1)}\otimes x_{(2)}$ for any $x\in\ff$. Then the coproduct of $\hB$ satisfies
	\begin{align*}
		\Delta(x)=\sum x_{(1)}h_{\wt(x_{(2)})}\otimes x_{(2)}, \qquad \Delta^2(x)=\sum x_{(1)}h_{\wt(x_{(2)})}h_{\wt(x_{(3)})}\otimes x_{(2)}h_{\wt(x_{(3)})}\otimes x_{(3)}.
	\end{align*}
	This convention greatly improves the clarity of the computation and will be adopted throughout this paper.
	
	We identify (cf. \eqref{eq:Udj-gen}) 
	\begin{align} \label{ffUU}
		\ff \stackrel{\cong}{\longrightarrow}\U^+, \,  \vartheta_i \mapsto \vartheta_i^+ :=E_i,  \qquad
		\ff \stackrel{\cong}{\longrightarrow}\U^-, \,\vartheta_i \mapsto \vartheta_i^- := F_i. 
	\end{align}
	Then we have $\hB\cong \U^+\otimes \Q(v^{\frac12})[h_i\mid i\in\I]$. 
	
	Let $\tB$ (resp.  $'\tB$) be the algebra constructed from $\hB$ (resp. $'\hB$) with $h_i$ invertible for $i\in\I$. Then we can define the coproduct, counit and antipode by
	\begin{gather}
		\Delta(\vartheta_i)= \vartheta_i\otimes 1 +  h_i\otimes \vartheta_i, \qquad
		\Delta(\h_{i}) = \h_{i} \otimes \h_{i},\quad \forall i\in\I;
		\\
		\varepsilon(\vartheta_i)= 0, \qquad
		\varepsilon(\h_{i}) = 1=\varepsilon(\h_i^{-1}),\quad \forall i\in\I;\\
		S(\vartheta_i)=-h_i^{-1}\vartheta_i,\quad S(h_i)=h_i^{-1},\quad \forall i\in\I.
	\end{gather}
	We also have  $\tB\cong \U^+\otimes \Q(v^{\frac12})[h_i^{\pm 1}\mid i\in\I]$. 
	In this way, ${'\tB},\tB$ are Hopf algebras, and ${'\hB},\hB$ are bialgebras. Define a Hopf pairing $\varphi$ by
	\begin{align}
		\label{eq:hopf-pairing}
		\varphi(\vartheta_i,\vartheta_j)=\delta_{ij}(v_i-v_i^{-1}),\quad \varphi(\h_i,\h_j)=v_i^{c_{ij}},\quad\varphi(\vartheta_i,\h_j)=0,\quad \forall i,j\in\I,x,y\in\tB.
	\end{align}
	Then it gives (symmetric) Hopf pairings on the Hopf algebras $'\tB,\tB$, and bialgebra pairings on the bialgebras $'\hB,\hB$. Moreover, the pairing $\varphi$ is non-degenerate on $\tB$ and $\hB$; see \cite[Corollary 1.13]{Sch12}, cf. \cite[\S1.2.12]{Lus93}.  We can view $\hB$ as a subalgebra of $\tB$ in the following.
	
	Consider the iHopf algebras of diagonal type, $(\hB\otimes\hB )^\imath$ and $(\tB\otimes\tB )^\imath$, associated to $(\hB,\varphi)$ and $(\tB,\varphi)$, respectively. Here $(\tB\otimes\tB )^\imath$ is a Hopf algebra by the constructions in Section~\ref{sec:i-Hopf}, and we also consider $(\hB\otimes\hB )^\imath$ as a subalgebra of $(\tB\otimes\tB )^\imath$.

	\begin{lemma}\label{lem:U=iHopfBB}
		We have a Hopf algebra isomorphism $\widetilde{\Phi}_{\sharp}:\tU\to (\tB\otimes\tB )^\imath$ given by
		\begin{align*}
			E_i\mapsto\vartheta_i\otimes 1,&\quad F_i\mapsto 1\otimes \vartheta_i,\quad K_i \mapsto  h_i\otimes1,\quad K_i'\mapsto 1\otimes h_i,\quad \forall i\in\I.
		\end{align*}
		Moreover, $\widetilde{\Phi}_{\sharp}(\hU)=(\hB\otimes\hB )^\imath$.
	\end{lemma}
	
	\begin{proof}
		It is well known that the Drinfeld double of $\tB$ is isomorphic to $\tU$ via
		\begin{align*}
			E_i\mapsto\vartheta_i\otimes 1,&\quad F_i\mapsto 1\otimes \vartheta_i,\quad K_i \mapsto  h_i\otimes1,\quad K_i'\mapsto 1\otimes h_i,\quad \forall i\in\I.
		\end{align*}
		Composed with the isomorphism $\Phi^\texttt{D}: D(\tB) \rightarrow (\tB\otimes\tB)^\imath$ from Proposition~\ref{prop:double-diagonal}, we obtain the desired isomorphism and $\widetilde{\Phi}_{\sharp}$. The second claim is now clear.
	\end{proof}
	
	Let $\tau$ be an involution in $\Inv(C)$. Clearly $\tau$ preserves the Hopf pairing $\varphi$. Denote by
	\[
	\hB^\imath_\tau = \text{iHopf}\, \big(\hB, \varphi \circ (\tau \otimes 1)\big), \qquad
	\tB^\imath_\tau =\text{iHopf}\, \big(\tB, \varphi\circ (\tau \otimes 1)\big)
	\]
	the iHopf algebras defined on $\big(\hB, \varphi \circ (\tau \otimes 1)\big)$ and $\big(\tB, \varphi\circ (\tau \otimes 1)\big)$, respectively. Again we may consider $\hB^\imath_\tau$ as a subalgebra of $\tB^\imath_\tau$.
	
	\begin{lemma}
		\label{lem:fvartheta}
		Let $\tau$ be an involution of $C$. Then there exists a unique $\tau$-twisted compatible map $\chi:\tB\rightarrow \Q(v^{\frac12})$ such that $\chi(1)=1$, $\chi(\vartheta_i)=0$, $\chi(h_i)=\varphi(h_i,h_{\tau i})$, for $i\in\I$.
	\end{lemma}
	
	\begin{proof}
		We first note that for $a,b,c\in{'\tB}$,
		\[\chi((ab)c)=\chi(a(bc))=\sum\chi(a_{(1)})\chi(b_{(2)})\varphi(\tau(a_{(2)}),b_{(1)})\chi(c_{(3})\varphi(\tau(a_{(3)}),c_{(1)})\varphi(\tau(b_{(3)}),c_{(2)}).\]
		Therefore, we can define a $\tau$-twisted compatible map $\chi:{'\tB}\to \Q(v^{\frac12})$ using the initial values $\chi(1)=1$, $\chi(\vartheta_i)=0$, $\chi(h_i)=\varphi(h_i,h_{\tau i})$, for $i\in\I$. To see that $\chi$ descends to $\tB$, we need to prove that it vanishes on the ideal generated by
		\begin{align*}
			u_{ij}:=\sum_{r+s=1-c_{ij}} (-1)^r\vartheta_i^{(r)}\vartheta_j\vartheta_i^{(s)},\qquad \forall i\neq j\in\I.
		\end{align*}
		
		We shall use the following identity:
		\[\Delta(u_{ij})=u_{ij}\otimes 1+h_{\beta}\otimes u_{ij}\]
		where $\beta=(1-c_{ij})\alpha_i+\alpha_j$. Since $u_{ij}$ belongs to the radical of $\varphi$, we obtain
		\begin{align*}
			\chi(au_{ij}b)=\sum\chi(a_{(1)})\chi(u_{ij})\chi(b_{(3)})\varphi(\tau b_{(2)},1)\varphi(\tau(h_\beta),a_{(2)})\varphi(\tau(b_{(1)}),a_{(3)}),\quad\forall a,b\in{'\tB}.
		\end{align*}
		It therefore suffices to show that $\chi(u_{ij})=0$. To this end, note that, for $r+s=1-c_{ij}$,
		\begin{equation}\label{eq: vartheta on Serre-1}
			\begin{aligned}
				&\chi(\vartheta_i^{(r)}\vartheta_j\vartheta_i^{(s)})\\
				&=\sum\chi((\vartheta_i^{(r)})_{(1)})\chi((\vartheta_j)_{(2)})\chi((\vartheta_i^{(s)})_{(3)}) \varphi(\tau(\vartheta_i^{(r)})_{(2)},(\vartheta_j)_{(1)})\\
				&\quad\qquad\cdot\varphi(\tau(\vartheta_i^{(r)})_{(3)},(\vartheta_i^{(s)})_{(1)})\varphi(\tau(\vartheta_j)_{(3)},(\vartheta_i^{(s)})_{(2)})\\
				&=\sum\chi((\vartheta_i^{(r)})_{(1)})\chi(1)\chi((\vartheta_i^{(s)})_{(3)})\varphi(\tau(\vartheta_i^{(r)})_{(2)},\vartheta_j)\varphi(\tau(\vartheta_i^{(r)})_{(3)},(\vartheta_i^{(s)})_{(1)})\varphi(1,(\vartheta_i^{(s)})_{(2)})\\
				&+\sum\chi((\vartheta_i^{(r)})_{(1)})\chi(h_j)\chi((\vartheta_i^{(s)})_{(3)}) \varphi(\tau(\vartheta_i^{(r)})_{(2)},h_j)\varphi(\tau(\vartheta_i^{(r)})_{(3)},(\vartheta_i^{(s)})_{(1)})\varphi(\vartheta_{\tau j},(\vartheta_i^{(s)})_{(2)})\\
				&=\sum\chi((\vartheta_i^{(r)})_{(1)})\chi((\vartheta_i^{(s)})_{(3)})\varphi(\tau(\vartheta_i^{(r)})_{(2)},\vartheta_j)\varphi(\tau(\vartheta_i^{(r)})_{(3)},(\vartheta_i^{(s)})_{(1)})\varphi(1,(\vartheta_i^{(s)})_{(2)})\\
				&+\sum\chi((\vartheta_i^{(r)})_{(1)})\chi(h_j)\chi((\vartheta_i^{(s)})_{(3)}) \varphi(\tau(\vartheta_i^{(r)})_{(2)},h_j)\varphi(\tau(\vartheta_i^{(r)})_{(3)},(\vartheta_i^{(s)})_{(1)})\varphi(\vartheta_{\tau j},(\vartheta_i^{(s)})_{(2)}),
			\end{aligned}
		\end{equation}
		where we have used $\Delta^{(2)}(\vartheta_j)=\vartheta_j\otimes 1\otimes 1+h_j\otimes\vartheta_j\otimes 1+h_j\otimes h_j\otimes \vartheta_j$. Now for $\varphi(\tau(\vartheta_i^{(r)})_{(2)},\vartheta_j)$ (and also $\varphi(\vartheta_{\tau j},(\vartheta_i^{(s)})_{(2)})$) to produce a nonzero term, it is necessary that $\tau j=i$, so let us assume this. Then \eqref{eq: vartheta on Serre-1} simplifies to
		\begin{align*}
			\chi(\vartheta_i^{(r)}\vartheta_j\vartheta_i^{(s)})&=\sum\chi((\vartheta_i^{(r)})_{(1)})\chi((\vartheta_i^{(s)})_{(2)})\varphi(\tau(\vartheta_i^{(r)})_{(2)},\vartheta_j)\varphi(\tau(\vartheta_i^{(r)})_{(3)},(\vartheta_i^{(s)})_{(1)})\\
			&\quad+\sum\chi((\vartheta_i^{(r)})_{(1)})\chi(h_j)\chi((\vartheta_i^{(s)})_{(3)}) \varphi(h_{\tau\wt((\vartheta_i^{(r)})_{(2)})},h_j)\varphi(\tau(\vartheta_i^{(r)})_{(2)},(\vartheta_i^{(s)})_{(1)})\\
			&\quad\qquad\cdot\varphi(\vartheta_{\tau j},(\vartheta_i^{(s)})_{(2)}).
		\end{align*}
		Since $j=\tau i\neq i$ under our assumption, the term $\varphi(\tau(\vartheta_i^{(r)})_{(3)},(\vartheta_i^{(s)})_{(1)})$ is nonzero only if $\tau(\vartheta_i^{(r)})_{(3)}=1$ (similar for the term $\varphi(\tau(\vartheta_i^{(r)})_{(2)},(\vartheta_i^{(s)})_{(1)})$), and this implies
		\begin{equation}\label{eq: vartheta on Serre-2}
			\begin{aligned}
				\chi(\vartheta_i^{(r)}\vartheta_j\vartheta_i^{(s)})&=\sum\chi((\vartheta_i^{(r)})_{(1)})\chi(\vartheta_i^{(s)})\varphi(\tau(\vartheta_i^{(r)})_{(2)},\vartheta_j)\\
				&\quad+\sum\chi(\vartheta_i^{(r)})\chi(h_j)\chi((\vartheta_i^{(s)})_{(2)})\varphi(\vartheta_{\tau j},(\vartheta_i^{(s)})_{(1)}).
			\end{aligned}
		\end{equation}
		As $i\neq\tau i$, an easy induction gives $\chi(\vartheta_i^n)=0$ for any $n\geq 1$. Also we have \cite[\S 3.1.5]{Lus93}
		\[
		\Delta(\vartheta_i^{(n)})=\sum_{a+b=n}v_i^{ab}\vartheta_i^{(a)}h_i^{b}\otimes\vartheta_i^{(b)}.
		\]
		Combining all of these establishes the following formula for $\chi(\vartheta_i^{(r)}\vartheta_j\vartheta_i^{(s)})$:
		\[\chi(\vartheta_i^{(r)}\vartheta_j\vartheta_i^{(s)})=\begin{cases}
			v_i^{-c_{ij}}\chi(\vartheta_i^{(-c_{ij})})\chi(h_i)\varphi(\vartheta_i,\vartheta_{\tau j})&\text{if $r=1-c_{ij},s=0$},\\
			v_i^{-c_{ij}}\chi(\vartheta_i^{(-c_{ij})})\chi(h_j)\varphi(\vartheta_i,\vartheta_{\tau j})&\text{if $r=0,s=1-c_{ij}$},\\
			0&\text{otherwise}.
		\end{cases}\]
		Note that the right-hand side is nontrivial only if $c_{ij}=0$, in which case $\chi(\vartheta_i\vartheta_j)=\chi(\vartheta_j\vartheta_i)$ since $\chi(h_i)=\chi(h_j)$. This implies $\chi(u_{ij})=0$, which completes the proof of the lemma. 
	\end{proof}
	
	Below we make the homomorphism ${\xi}_{\tau,\chi}$ \eqref{xi} in Proposition \ref{prop:embedding-iHopf-diagonal} explicit and drop the index $\chi$ as it only depends on $\tau$.
	
	\begin{lemma} \label{lem:iHembed}
		There is an embedding of algebras
		\begin{equation} \label{eq:xi_theta embedding of iHopf}
			\widetilde{\xi}_\tau: \tB^\imath_\tau\longrightarrow(\tB\otimes\tB)^\imath
		\end{equation}
		which send
		\begin{align}
			\label{eq:Xitheta}
			\vartheta_{i}\mapsto 1\otimes\vartheta_i+(\vartheta_{\tau i}\otimes 1)\ast(1\otimes h_i), 
			\quad
			h_{i}\mapsto (h_{\tau i}\otimes 1)\ast(1\otimes h_i),
			\qquad\forall i\in \I.
		\end{align}
		Moreover, $\widetilde{\xi}_\tau(\hB^\imath_\tau)\subseteq (\hB\otimes\hB)^\imath$.
	\end{lemma}
	
	\begin{proof}
		By Proposition \ref{prop:embedding-iHopf-diagonal}, the $\tau$-twisted compatible map $\chi:\tB\rightarrow \Q(v^{\frac12})$ defined in Lemma~\ref{lem:fvartheta} gives rise to an algebra homomorphism $\widetilde{\xi}_\tau: \tB^\imath_\tau\rightarrow(\tB\otimes\tB)^\imath$. More precisely, for $x\in\hB$, by Proposition \ref{prop:embedding-iHopf-diagonal},  we have
		\[
		\widetilde{\xi}_\tau(x)=\sum\chi(x_{(2)}h_{\wt(x_{(3)})})\tau (x_{(3)})\otimes x_{(1)}h_{\wt(x_{(2)})}h_{\wt(x_{(3)})}.
		\]
		Note that the right-hand side has a leading term $1\otimes x$. It then follows by induction on degree that $\widetilde{\xi}_\tau$ is an embedding of algebras. One also checks that 
		\begin{align*}
			\widetilde{\xi}_\tau(\vartheta_i)&=1\otimes \vartheta_i+\chi(h_i)\vartheta_{\tau i}\otimes h_i=1\otimes\vartheta_i+\varphi(h_i,h_{\tau i})\vartheta_{\tau i}\otimes h_i\\
			&=1\otimes\vartheta_i+(\vartheta_{\tau i}\otimes 1)\ast(1\otimes h_i).
		\end{align*}
		The formula for $\widetilde{\xi}_\tau(h_i)$ is similarly verified. Finally, the second claim is immediate.
	\end{proof}
	
	Recall from \eqref{eq:Bi} the generators $B_i= F_i +  E_{\btau i} \tK_i'$ of the universal iquantum groups $\hUi$ and $\tUi$.
	
	\begin{theorem}\label{thm:iH=iQG}
		We have an algebra isomorphism $\widetilde{\Phi}^\imath:\tB^\imath_\tau\longrightarrow \tUi$ given by
		\begin{align*} 
			\vartheta_i &\mapsto B_i,\qquad h_i\mapsto \tk_{\tau i},\quad \forall i\in\I.
		\end{align*}
		Moreover, $\widetilde{\Phi}^\imath(\hB^\imath_\tau)=\hUi$.
	\end{theorem}
	
	\begin{proof}
		A simple computation shows that \begin{align*}
			\widetilde{\Phi}_\sharp^{-1}\circ\widetilde{\xi}_\tau(\vartheta_i)&=\widetilde{\Phi}_\sharp^{-1}(1\otimes\vartheta_i+(\vartheta_{\tau i}\otimes 1)\ast(1\otimes h_i))=F_i+E_{\tau i}K_i',\\
			\widetilde{\Phi}_\sharp^{-1}\circ\widetilde{\xi}_\tau(h_i)&=K_i'K_{\tau i}
		\end{align*}
		by using \eqref{eq:Xitheta}. Thus there exists an algebra homomorphism 
		$\widetilde{\Phi}^\imath:\tB^\imath_\tau\rightarrow \tUi$, $\vartheta_i\mapsto B_i,h_i\mapsto \tk_{\tau i}$ ($i\in\I$), such that the following  diagram commutes:
		\begin{equation} \label{CD:fU}
			\begin{tikzcd}
				\tB^\imath_\tau\ar[r,"\widetilde{\xi}_\tau"]\ar[d,dashed,swap,"\widetilde{\Phi}^\imath"]&(\tB\otimes\tB)^\imath\ar[d,"\widetilde{\Phi}_\sharp^{-1}"]\\
				\tUi\ar[r,hook]&\tU
			\end{tikzcd}
		\end{equation}
		In particular, $\widetilde{\Phi}^\imath$ sends $\vartheta_i$ to $B_i$ and $h_i$ to $\tk_{\tau i}$, hence is surjective. Since $\widetilde{\Phi}_\sharp$ is an isomorphism and $\widetilde{\xi}_\tau$ is an embedding, $\widetilde{\Phi}^\imath$ is also injective and hence an isomorphism. The last statement is clear.
	\end{proof}

	\subsection{A recursive formula}
	
	There exist linear maps known as skew-derivations (cf. \cite{Lus93}) 
	\[
	\partial_i^R:\ff\longrightarrow \ff, \qquad
	\partial_i^L:\ff\longrightarrow \ff
	\]
	such that $\partial_i^R(1)=\partial_i^L(1)=0$, $\partial_i^R(\vartheta_j)=\delta_{ij}=\partial_i^L(\vartheta_j)$, and
	\begin{align*}
		&\quad    
		\partial_i^R(fg)=\partial_i^R(f)g+v^{(\alpha_i,\mu)}f\partial_i^R(g),
		\\
		&\quad \partial_i^L(fg)=v^{(\alpha_i,\nu)}\partial_i^L(f)g+f\partial_i^L(g),
	\end{align*}
	for any $j\in\I, f\in\ff_\mu, g\in\ff_\nu$.
	Then, for $x\in\ff_\mu$, we have 
	\begin{align} \label{derivations}
		\begin{split}
			\Delta(x)&= x\otimes 1+\sum_{i\in\I}\partial_i^L(x)h_i\otimes \vartheta_i+(\text{rest})
			\\
			&= h_\mu\otimes x+\sum_{i\in\I}\vartheta_ih_{\mu-\alpha_i}\otimes\partial_i^R(x)+(\text{rest}).
		\end{split}
	\end{align}
	
	Recall the two algebras $(\widehat\BB, \cdot)$ and $(\hB^\imath_\tau, *)$ have the same underlying vector space (which contains $\ff$ as a subspace). 
	
	\begin{lemma}
		\label{lem:recursive}
		In $\hB^\imath_\tau$ (and $\tB^\imath_\tau$), for $x\in\ff$ and $i\in \I$, we have
		\begin{align*}
			\vartheta_i*x&=\vartheta_i\cdot x+(v_i-v_i^{-1})\partial_{\tau i}^L(x)\cdot h_{\tau i},
			\\
			x*\vartheta_i&=x\cdot \vartheta_i+(v_i-v_i^{-1})\partial_{\tau i}^R(x)\cdot h_{i}.
		\end{align*}
	\end{lemma}
	
	\begin{proof}
		Using $\Delta(\vartheta_i)=\vartheta_i\otimes 1+h_i\otimes\vartheta_i$ and \eqref{derivations}, we compute by using the multiplication formula \eqref{star product} that 
		\begin{align*}
			\vartheta_i*x &=\vartheta_i\cdot x+\varphi(\vartheta_i,\vartheta_{i })\partial_{\tau i}^L(x)\cdot h_{\tau i}
			\\
			&=\vartheta_i\cdot x+(v_i-v_i^{-1})\partial_{\tau i}^L(x)\cdot h_{\tau i}.
		\end{align*}
		The second formula is proved in a similar way.
	\end{proof}
	
	\begin{remark}  \label{rem:KY}
		Let $\mathcal{A}_\tau$ be the subalgebra of $\U$ generated by $F_i$, $K_iK_{\tau i}^{-1}$, for $i\in\I$. The subalgebras $\widehat{\ca}_\tau\subset \hU$ and $\widetilde{\ca}_\tau\subset\tU$ can be similarly defined via generators 
		$F_i$, $K_iK_{\tau i}'$, for $i\in\I$. 
		Kolb and Yakimov defined a star product $\star$ on $\mathcal{A}_\tau$ via the quasi R-matrix, which is uniquely determined by recursive formulas \cite[Lemma 4.7]{KY20}. Up to a straightforward modification of their construction by use of $\widehat{\mathcal{A}}_\tau$ and  $\widetilde{\ca}_\tau$ in place of $\mathcal{A}_\tau$, their recursive formulas will match with Lemma \ref{lem:recursive}. Braid group symmetries and dual canonical bases were not discussed {\em loc. cit.} It is possible to reestablish in our setting  a main result of that paper (see \cite[(6.1)]{KY20}) connecting the quasi K-matrix to quasi R-matrix, but we refrain from getting into the detail here as this requires much notational adjustment and renormalization.
	\end{remark}

	
	\section{$\mathrm{i}$Braid group symmetries on $\tUi$ via $\mathrm{i}$Hopf algebra} 
	\label{sec:braid}
	
	In this section, we shall give an iHopf algebra interpretation of the ibraid group action on iquantum groups, providing a new connection to braid group action of quantum groups. 
	
	\subsection{Connecting 2 braid group actions via iHopf}
	
	Let $j\neq i,\tau i$ in $\I$ in this section. Define the root vectors in $\hB$:
	\begin{align}
		f_{i,j;m}=\sum_{r+s=m}(-1)^rv_i^{r(c_{ij}+m-1)+\frac{1}{2}m}(v_i-v_i^{-1})^{-m}\vartheta_i^{(r)}\vartheta_j\vartheta_i^{(s)},\label{eq:root vector i=taui-1}
		\\
		f'_{i,j;m}=\sum_{r+s=m}(-1)^rv_i^{r(c_{ij}+m-1)+\frac{1}{2}m}(v_i-v_i^{-1})^{-m}\vartheta_i^{(s)}\vartheta_j\vartheta_i^{(r)},\label{eq:root vector i=taui-2}
	\end{align}
	which are slightly normalized versions of Lusztig's definition \cite{Lus93}. By the same proof of \cite[Proposition 37.2.5]{Lus93}, the rescaled braid group symmetry $\widetilde{T}_i$ satisfies
	\begin{align}\label{eq: tU T_i action on f_m}
		\widetilde{T}_i(f_{i,j;m}^+)=f_{i,j,-c_{ij}-m}'^+,\qquad \forall m,n\in\Z.
	\end{align}
	
	Let $\ad:\ff\to\ff$ be the adjoint action via the identification $\ff\cong\U^+$; it is given by
	\begin{align}\label{eq: ad action def}
		\ad(\vartheta_i)(x)=\vartheta_ix-h_ixh_i^{-1}\vartheta_i.
	\end{align}
	Recalling the anti-involution $\sigma$ from Lemma~\ref{lem:anti-involut-QG}, we have
	\begin{align} \label{rootvector2}
		f_{i,j;m}=v_i^{\frac{1}{2}m}(v_i-v_i^{-1})^{-m}\sigma\big(\ad(\vartheta_i^{(m)})(\vartheta_j)\big),
		\quad 
		f'_{i,j;m}=v_i^{\frac{1}{2}m}(v_i-v_i^{-1})^{-m}\ad(\vartheta_i^{(m)})(\vartheta_j).
	\end{align}
	
	The ibraid group symmetry $\tTT_i$ of $\tUi$ can be transported to $\tB^\imath_\tau$ via the isomorphism $\widetilde{\Phi}^\imath:\tB^\imath_\tau \stackrel{\cong}{\rightarrow} \tUi$ in Theorem \ref{thm:iH=iQG}. Recalling $\K_\alpha\in \tUi$ from \eqref{eq:bbKi}, we define
	\begin{align}
		\K_\alpha:=v^{\frac{1}{2}(\alpha,\tau\alpha)}h_{\tau\alpha} \in \tB^\imath_\tau
	\end{align}
	so that $\widehat{\Phi}^\imath(\K_\alpha)=\K_\alpha$. Then $\tTT_i(\K_\alpha)=\K_{r_i\alpha}$ in $\tB^\imath_\tau$; see  \eqref{def:ri} for definition of $r_i$. Denote by 
	\begin{align} \label{embed:f}
		\iota:\tB \longrightarrow \tB^\imath_\tau
	\end{align}
	the canonical embedding. The following is the main result of this section. 
	
	\begin{theorem}
		\label{thm:braidiHopf}
		For any $i,j\in \I$ such that $i\neq j,\tau j$, we have $$\tTT_i(\vartheta_j)=\iota\big(\widetilde{T}_{r_i}(\vartheta_j)\big).$$
	\end{theorem}
	
	The proof of the theorem is long,  and divided into three cases in \S\ref{subsec: root vector taui=i}--\ref{subsec: root vector c_i,taui=-1} below, depending on the rank one type.

	\subsection{The case for $c_{i,\tau i}=2$}\label{subsec: root vector taui=i}
	
	For $i\in \I$ with $\btau i= i$, generalizing the constructions in \cite{BW18b, BeW18}, we define the {\em idivided powers} of $\vartheta_i$ in $\tB_\tau^\imath$ to be (see also \cite{CLW21a})
	\begin{eqnarray}
		&&\vartheta_{i,\odd}^{(m)}=\frac{1}{[m]_{v_i}^!}\left\{ \begin{array}{ccccc} \vartheta_i*\prod_{s=1}^k (\vartheta_i*\vartheta_i+v_i(v_i-v_i^{-1})^2h_{i}[2s-1]_{v_i}^2 ) & \text{if }m=2k+1,\\
			\prod_{s=1}^k (\vartheta_i*\vartheta_i+v_i(v_i-v_i^{-1})^2h_i[2s-1]_{v_i}^2) &\text{if }m=2k; \end{array}\right.
		\label{eq:iDPodd} \\
		&&\vartheta_{i,\ev}^{(m)}= \frac{1}{[m]_{v_i}^!}\left\{ \begin{array}{ccccc} \vartheta_i*\prod_{s=1}^k (\vartheta_i*\vartheta_i+v_i(v_i-v_i^{-1})^2h_i[2s]_{v_i}^2 ) & \text{if }m=2k+1,\\
			\prod_{s=1}^{k} (\vartheta_i*\vartheta_i+v_i(v_i-v_i^{-1})^2h_i[2s-2]_{v_i}^2) &\text{if }m=2k. \end{array}\right.
		\label{eq:iDPev}
	\end{eqnarray}
	
	These idivided powers satisfy the following recursive relations:
	\begin{align}
		\label{lem:dividied power}
		\vartheta_{i}*\vartheta_{i,\ov{p}}^{(r)}=\left\{\begin{array}{llll}
			&[r+1]_{v_i}\vartheta_{i,\ov{p}}^{(r+1)}                                 & \text{if}\ \ov{p}\neq \ov{r};\\
			&[r+1]_{v_i}\vartheta_{i,\ov{p}}^{(r+1)}-v_i(v_i-v_i^{-1})^2[r]_{v_i}h_i\vartheta_{i,\ov{p}}^{(r-1)} & \text{if}\ \ov{p}= \ov{r}.
		\end{array}\right.
	\end{align}
	
	Recall the definition of $f_{i,j;m}$ and $f'_{i,j;m}$ in Section~\ref{sec:braid}. By comparing the our definition of root vectors with \cite[\S 37.2]{Lus93}, one can deduce the following relations.
	
	\begin{lemma}[{\cite[Lemma 7.1.2]{Lus93}}]\label{lem: f_m commutation with E_iF_i}
		For $m\geq 0$, in $\tU$ we have:  
		\begin{align*}
			f_{i,j;m}^+E_i-v_i^{c_{ij}+2m}E_if_{i,j;m}^+&=v_i^{-\frac{1}{2}}(v_i-v_i^{-1})[m+1]_{v_i}f_{i,j;m+1}^+,\\
			[F_i,f_{i,j;m}^+]&=v_i^{\frac{1}{2}}(v_i-v_i^{-1})[-c_{ij}-m+1]_{v_i}K_if_{i,j;m-1}^+.
		\end{align*}
	\end{lemma}
	
	Similar recursive formulas can be formulated in $\tB^\imath_\tau$, cf. \cite[\S 6.1]{CLW21b}.
	
	\begin{lemma}\label{lem: f_m induction relation}
		Let $m\in\Z$,  then in $\tB^\imath_\tau$,
		\begin{align*}
			&f_{i,j;m}\ast\vartheta_i-v_i^{c_{ij}+2m}\vartheta_i\ast f_{i,j;m}\\
			&=v_i^{-\frac{1}{2}}(v_i-v_i^{-1})[m+1]_{v_i}f_{i,j;m+1}-v_i^{c_{ij}+2m+\frac{1}{2}}(v_i-v_i^{-1})[-c_{ij}-m+1]_{v_i}h_i\ast f_{i,j;m-1},\\
			&\vartheta_i\ast f'_{i,j;m}-v_i^{c_{ij}+2m}f'_{i,j;m}\ast\vartheta_i\\
			&=v_i^{-\frac{1}{2}}(v_i-v_i^{-1})[m+1]_{v_i}f'_{i,j;m+1}-v_i^{c_{ij}+2m+\frac{1}{2}}(v_i-v_i^{-1})[-c_{ij}-m+1]_{v_i}h_i\ast f'_{i,j;m-1}.
		\end{align*}
	\end{lemma}
	\begin{proof}
		For simplicity let us write $f_m:=f_{i,j;m}$, $f'_m:=f'_{i,j;m}$. By \cite[Lemma 38.1.7]{Lus93}, we get
		\begin{align*}
			\Delta(f_m)&=h_{i}^m h_j\otimes f_m+\sum_{t=0}^{m}\prod_{h=0}^{m-t-1}\frac{1-v_i^{2(m+c_{ij}-h-1)}}{v_i-v_i^{-1}}v_i^{(t+\frac{1}{2})(m-t)}f_th_i^{m-t}\otimes\vartheta_i^{(m-t)},\\
			\Delta(f'_m)&=f_m'\otimes 1+\sum_{t=0}^{m}\prod_{h=0}^{m-t-1}\frac{1-v_i^{2(m+c_{ij}-h-1)}}{v_i-v_i^{-1}}v_i^{(t+\frac{1}{2})(m-t)}\vartheta_i^{(m-t)}h_i^th_j\otimes f'_t.
		\end{align*}
		The assertion then follows by a direct computation in $\tB^\imath_\tau$.
	\end{proof}
	
	For $m\in\Z$, we denote by $\ov{m}$ the class of $m$ in $\Z_2$. By comparing Lemma~\ref{lem: f_m induction relation} with \cite[Theorem 6.2]{CLW21b}, one can find the following explicit expressions of $f_{i,j;m}$ and $f'_{i,j;m}$.
	
	\begin{proposition}[{\cite[Theorem 6.2]{CLW21b}}]\label{prop: f_m expression by i-div}
		Let $m\in\Z$, $\ov{p}\in\Z_2$ and $\delta_{\ov{r},\ov{p}}$ be the Kronecker symbol.
		\begin{enumerate}
			\item[(a)] If $\ov{m}\neq\ov{c_{ij}}$, then 
			\begin{align}
				f_{i,j;m}=\sum_{u\geq 0}\sum_{r+s+2u=m}(-1)^{r+u}v_i^{\frac{1}{2}(m-2u)}&(v_i-v_i^{-1})^{-m+2u}\times v_i^{r(c_{ij}+m-1)+u(c_{ij}+m+2\delta_{\ov{r},\ov{p}})+u(u-1)}\notag\\ \label{eq:f_m expression-1}
				&\times\qbinom{\frac{c_{ij}+m-1}{2}}{u}_{v_i^2}h_i^u\ast\vartheta_{i,\ov{p}}^{(r)}\ast \vartheta_j\ast\vartheta_{i,\ov{p}+\ov{c_{ij}}}^{(s)},
			\end{align}
			and
			\begin{align}
				f'_{i,j;m}=\sum_{u\geq 0}\sum_{r+s+2u=m}(-1)^{r+u}v_i^{\frac{1}{2}(m-2u)}&(v_i-v_i^{-1})^{-m+2u}\times v_i^{r(c_{ij}+m-1)+u(c_{ij}+m+2\delta_{\ov{r},\ov{p}})+u(u-1)}\notag\\\label{eq:f'_m expression-1}
				&\times\qbinom{\frac{c_{ij}+m-1}{2}}{u}_{v_i^2}h_i^u\ast\vartheta_{i,\ov{p}}^{(s)}\ast \vartheta_j\ast\vartheta_{i,\ov{p}+\ov{c_{ij}}}^{(r)}.
			\end{align}
			\item[(b)] If $\ov{m}=\ov{c_{ij}}$, then
			\begin{equation}\label{eq:f_m expression-2}
				\begin{aligned}
					f_{i,j;m}=\sum_{u\geq 0}\sum_{r+s+2u=m}(-1)^{r+u}v_i^{\frac{1}{2}(m-2u)}&(v_i-v_i^{-1})^{-m+2u}\times v_i^{r(c_{ij}+m-1)+u(c_{ij}+m+1)+u(u-1)}\\
					&\times\qbinom{\frac{c_{ij}+m-2\delta_{\ov{r},\ov{p}}}{2}}{u}_{v_i^2}h_i^u\ast\vartheta_{i,\ov{p}}^{(r)}\ast \vartheta_j\ast\vartheta_{i,\ov{p}+\ov{c_{ij}}}^{(s)},
				\end{aligned}
			\end{equation}
			and
			\begin{equation}\label{eq:f'_m expression-2}
				\begin{aligned}
					f'_{i,j;m}=\sum_{u\geq 0}\sum_{r+s+2u=m}(-1)^{r+u}v_i^{\frac{1}{2}(m-2u)}&(v_i-v_i^{-1})^{-m+2u}\times v_i^{r(c_{ij}+m-1)+u(c_{ij}+m+1)+u(u-1)}\\
					&\times\qbinom{\frac{c_{ij}+m-2\delta_{\ov{r},\ov{p}}}{2}}{u}_{v_i^2}h_i^u\ast\vartheta_{i,\ov{p}}^{(s)}\ast \vartheta_j\ast\vartheta_{i,\ov{p}+\ov{c_{ij}}}^{(r)}.
				\end{aligned}
			\end{equation}
		\end{enumerate}
		In particular, we have $\tTT_i(\vartheta_j)=f'_{i,j;-c_{ij}}$.
	\end{proposition}
	\begin{proof}
		In \cite[\S6.1]{CLW21b} the authors defined elements $\tilde{y}_{i,j;1,m,\bar{p},\bar{t},-1}$, $\tilde{y}'_{i,j;1,m,\bar{p},\bar{t},-1}$ satisfy the following recursive relations: 
		\begin{align*}
			&\tilde{y}_{i,j;1,m,\bar{p},\bar{t},-1}B_i-v_i^{2m+c_{ij}}B_i\tilde{y}_{i,j;1,m,\bar{p},\bar{t},-1}\\
			&=[m+1]_i\tilde{y}_{i,j;1,m+1,\bar{p},\bar{t},-1}+[-c_{ij}-m+1]_iv_i^{2m+c_{ij}}\tilde{k}_i\tilde{y}_{i,j;1,m-1,\bar{p},\bar{t},-1},\\
			&B_i\tilde{y}'_{i,j;1,m,\bar{p},\bar{t},-1}-v_i^{2m+c_{ij}}\tilde{y}'_{i,j;1,m,\bar{p},\bar{t},-1}B_i\\
			&=[m+1]_i\tilde{y}'_{i,j;1,m+1,\bar{p},\bar{t},-1}+[-c_{ij}-m+1]_iv_i^{2m+c_{ij}}\tilde{k}_i\tilde{y}'_{i,j;1,m-1,\bar{p},\bar{t},-1}.
		\end{align*}
		The general formulas of $\tilde{y}_{i,j;1,m,\bar{p},\bar{t},-1}$ and $\tilde{y}'_{i,j;1,m,\bar{p},\bar{t},-1}$ are given in \cite[(6.1)--(6.4)]{CLW21b}. By Lemma~\ref{lem: f_m induction relation}, $f_{i,j;m}$ and $f'_{i,j;m}$ satisfies the same recursive relation (up to some normalization on coefficients), so we can use the proof of \cite[Theorem 6.2]{CLW21b} to deduce \eqref{eq:f_m expression-1}--\eqref{eq:f'_m expression-2} (cf. \cite[\S 6.4, 6.5]{CLW21b}).
		
		Finally, setting $m=-c_{ij}$ we see that
		\begin{align*}
			f'_{i,j;-c_{ij}}&=\sum_{r+s=-c_{ij}}(-1)^{r}v_i^{-r-\frac{1}{2}c_{ij}}(v_i-v_i^{-1})^{c_{ij}}\vartheta_{i,\ov{p}}^{(s)}\ast\vartheta_j\ast\vartheta_{i,\ov{c_{ij}}+\ov{p}}^{(r)}\\
			&+\sum_{u\geq 1}\sum_{\stackrel{r+s=-c_{ij}-2u}{\ov{r}=\ov{p}}}(-1)^{r}v_i^{\frac{1}{2}(-c_{ij}-2u)}(v_i-v_i^{-1})^{c_{ij}+2u}\times v_i^{-r+u+u(u-1)} h_i^u\ast\vartheta_{i,\ov{p}}^{(s)}\ast\vartheta_j\ast\vartheta_{i,\ov{c_{ij}}+\ov{p}}^{(r)}.
		\end{align*}
		Under the isomorphism $\widetilde{\Phi}^\imath:\tB^\imath_\tau\to \tUi$, the right-hand side is mapped to $\tTT_i(B_j)$ (see \eqref{eq: T_i(B_j) for i=taui}), so the last statement follows.
	\end{proof}
	
	Using the recursive formulas of Lemma~\ref{lem: f_m induction relation} together with the observation $\tTT_i(f_{i,j;0})=\tTT_i(\vartheta_j)=f'_{i,j;-c_{ij}}$, we are ready to prove the following general result for the action of $\tTT_i$ on root vectors, cf. \cite[Theorem 7.1]{Lus93}.
	
	\begin{proposition}\label{prop:ibraid:split}
		Let $\tau i=i\neq j$ and $m\in\Z$. Then
		\begin{equation}\label{eq: rel braid action on root vector i=taui}
			\tTT_i(f_{i,j;m})=f'_{i,j;-c_{ij}-m}.
		\end{equation}
	\end{proposition}
	
	\begin{proof}
		For simplicity we set $f_m=f_{i,j;m}$ and $f'_m=f'_{i,j;m}$. The case $m=0$ is already known. Now we proceed by induction on $m$. Assume that $\tTT_i(f_t)=f'_{-c_{ij}-t}$ for $0\leq t\leq m$. Now applying $\tTT_i$ to
		\begin{equation}\label{eq:f(m) and theta_i commutator}
			f_m\vartheta_i-v^{c_{ij}+2m}\vartheta_if_m=v_i^{-\frac{1}{2}}(v_i-v_i^{-1})[m+1]_{v_i}f_{m+1}
		\end{equation}
		gives
		\begin{equation}\label{eq: T_if(m) induction}
			v^{-\frac{1}{2}}(v-v^{-1})[m+1]_{v_i}\tTT_i(f_{m+1})=\tTT_i(f_m\vartheta_i)-v^{c_{ij}+2m}\tTT_i(\vartheta_if_m).
		\end{equation}
		Using Lemma~\ref{lem: f_m induction relation} we then deduce that
		\begin{align*}
			\tTT_i(f_m\vartheta_i)&=\tTT_i(f_m\ast\vartheta_i)=\K_i^{-1}\ast f'_{-c_{ij}-m}\ast\vartheta_i\\
			&=\K_i^{-1}\ast(f'_{-c_{ij}-m}\vartheta_i+v_i^{-\frac{1}{2}}(v_i-v_i^{-1})[m+1]_{_i}\K_i\ast f'_{-c_{ij}-m-1})\\
			&=\K_i^{-1}\ast (f'_{-c_{ij}-m}\vartheta_i)+v_i^{-\frac{1}{2}}(v_i-v_i^{-1})[m+1]_{v_i}f'_{-c_{ij}-m-1},
		\end{align*}
		and
		\begin{align*}
			\tTT_i(\vartheta_if_m)&=\tTT_i(\vartheta_i\ast f_m-v_i^{-\frac{1}{2}}(v_i-v_i^{-1})[-c_{ij}-m+1]_{v_i}\K_i\ast f_{m-1})\\
			&=\tTT_i(\vartheta_i)\ast \tTT_i(f_m)-v_i^{-\frac{1}{2}}(v_i-v_i^{-1})[-c_{ij}-m+1]_{v_i}\K_i^{-1}\ast f'_{-c_{ij}-m+1})\\
			&=\K_i^{-1}\ast \vartheta_i\ast f'_{-c_{ij}-m}-v_i^{-\frac{1}{2}}(v_i-v_i^{-1})[-c_{ij}-m+1]_{v_i}\K_i^{-1}\ast f'_{-c_{ij}-m+1})\\
			&=\K_i^{-1}\ast (\vartheta_if'_{-c_{ij}-m})-v_i^{-\frac{1}{2}}(v_i-v_i^{-1})[-c_{ij}-m+1]_{v_i}\K_i^{-1}\ast f'_{-c_{ij}-m+1}).
		\end{align*}
		Now \eqref{eq: T_if(m) induction} can be computed as
		\begin{equation}\label{eq:T_i of f(m) compute}
			\begin{aligned}
				&v_i^{-\frac{1}{2}}(v_i-v_i^{-1})[m+1]_{v_i}\tTT_i(f_{m+1})\\
				&=\K_i^{-1}\ast(f'_{-c_{ij}-m}\vartheta_i-v_i^{c_{ij}+2m}\vartheta_if'_{-c_{ij}-m})+v_i^{-\frac{1}{2}}(v_i-v_i^{-1})[m+1]_{v_i}f'_{-c_{ij}-m-1}\\
				&\quad+v_i^{c_{ij}+2m-\frac{1}{2}}(v_i-v_i^{-1})[-c_{ij}-m+1]_{v_i}\K_i^{-1}\ast f'_{-c_{ij}-m+1}.
			\end{aligned}
		\end{equation}
		Note that from \eqref{eq:f(m) and theta_i commutator} we also get
		\[f'_{-c_{ij}-m}\vartheta_i-v_i^{c_{ij}+2m}\vartheta_if'_{-c_{ij}-m}=-(v_i-v_i^{-1})v_i^{c_{ij}+2m-\frac{1}{2}}[-c_{ij}-m+1]_{v_i}f'_{-c_{ij}-m+1}.\]
		Plugging this into \eqref{eq:T_i of f(m) compute}, we see that
		\[v_i^{-\frac{1}{2}}(v_i-v_i^{-1})[m+1]_{v_i}\tTT_i(f_{m+1})=v_i^{-\frac{1}{2}}(v_i-v_i^{-1})[m+1]_{v_i}f'_{-c_{ij}-m-1}.\]
		Since $v_i^{-\frac{1}{2}}(v_i-v_i^{-1})[m+1]_{v_i}\neq 0$, we deduce that $\tTT_i(f_{m+1})=f'_{-c_{ij}-m-1}$. This completes the induction step.
	\end{proof}

	\subsection{The case for $c_{i,\tau i}=0$}\label{subsec: root vector c_i,taui=0}
	
	For $i\neq\tau i$ and $c_{i,\tau i}=0$, we define the root vectors by 
	\begin{align}
		\begin{split}f_{i,\tau i,j;m,n}&=\sum_{\substack{r_1+s_1=m\\ r_2+s_2=n}}(-1)^{r_1+r_2}v_i^{r_1(c_{ij}+m-1)+r_2(c_{\tau i,j}n-1)+\frac{1}{2}(m+n)}(v_i-v_i^{-1})^{-m-n}\\
			&\qquad\qquad \times\vartheta_i^{(r_1)}\vartheta_{\tau i}^{(r_2)}\vartheta_j\vartheta_{\tau i}^{(s_2)}\vartheta_i^{(s_1)},\label{eq:root vector c_i,taui=0-1}
		\end{split}
		\\
		\begin{split}
			f'_{i,\tau i,j;m,n}&=\sum_{\substack{r_1+s_1=m\\ r_2+s_2=n}}(-1)^{r_1+r_2}v_i^{r_1(c_{ij}+m-1)+r_2(c_{\tau i,j}n-1)+\frac{1}{2}(m+n)}(v_i-v_i^{-1})^{-m-n}
			\\
			&\qquad\qquad\times\vartheta_i^{(s_1)}\vartheta_{\tau i}^{(s_2)}\vartheta_j\vartheta_{\tau i}^{(r_2)}\vartheta_i^{(r_1)}.\label{eq:root vector c_i,taui=0-2}
		\end{split}
	\end{align}
	
	The following lemma provides a counterpart of \eqref{eq: tU T_i action on f_m} for $f_{i,\tau i,j;m,n}$ and $f'_{i,\tau i,j;m,n}$.
	
	\begin{lemma}\label{eq: tU T_ri action on f_mn}
		For $m,n\in\Z$, we have in $\tU$,
		\begin{align}
			\widetilde{T}_{r_i}(f_{i,\tau i,j;m,n}^+)&=f'^+_{i,\tau i,j;-c_{ij}-m,-c_{\tau i,j}-n}.
		\end{align}
	\end{lemma}
	\begin{proof}
		Recall that $\widetilde{T}_i(f_{i,j;m}^+)=f'^+_{i,j;-c_{ij}-m}$. Now
		\begin{align*}
			&\widetilde{T}_{r_i}(f_{i,\tau i,j;m,n}^+)\\
			&=\widetilde{T}_{\tau i}\widetilde{T}_i\Big(\sum_{\substack{r_1+s_1=m\\ r_2+s_2=n}}(-1)^{r_1+r_2}v_i^{-r_1(-c_{ij}-m+1)-r_2(-c_{\tau i,j}-n+1)+\frac{1}{2}(m+n)}\\
			&\hspace{2cm}\times(v_i-v_i^{-1})^{-m-n}E_i^{(r_1)}E_{\tau i}^{(r_2)}E_jE_{\tau i}^{(s_2)}E_i^{(s_1)}\Big)\\
			&=\widetilde{T}_{\tau i}\Big(\sum_{r_2+s_2=n}(-1)^{r_2}v_i^{-r_2(-c_{\tau i,j}-n+1)+\frac{1}{2}n}(v_i-v_i^{-1})^{c_{\tau i,j}}E_{\tau i}^{(r_2)}T_i(f_{i,j;m}^+)E_{\tau i}^{(s_2)}\Big)\\
			&=\widetilde{T}_{\tau i}\Big(\sum_{r_2+s_2=n}(-1)^{r_2}v_i^{-r_2(-c_{\tau i,j}-n+1)+\frac{1}{2}n}(v_i-v_i^{-1})^{c_{\tau i,j}}E_{\tau i}^{(r_2)}f'^+_{i,j;-c_{ij}-m}E_{\tau i}^{(s_2)}\Big)\\
			&=\sum_{r_1+s_1=-c_{ij}-m}(-1)^{r_1}v_i^{-r_1(m+1)+\frac{1}{2}(-c_{ij}-m)}E_i^{(s_1)}\widetilde{T}_{\tau i}(f_{\tau i;n}^+)E_i^{(r_1)}\\
			&=\sum_{r_1+s_1=-c_{ij}-m}(-1)^{r_1}v_i^{-r_1(m+1)+\frac{1}{2}(-c_{ij}-m)}E_i^{(s_1)}f_{\tau i;-c_{\tau i,j}-n}^+E_i^{(r_1)}\\
			&=f'^+_{i,j;-c_{ij}-m,-c_{\tau i,j}-n}.\qedhere
		\end{align*}
	\end{proof}

	Since $\vartheta_i$ commutes with $\vartheta_{\tau i}$ in our case, from Lemma~\ref{lem: f_m commutation with E_iF_i} we immediately deduce the following.
	
	\begin{lemma}\label{lem:f_mcommutateE_iF_i2}
		For $m,n\in\Z$, we have
		\begin{align*}
			f_{i,\tau i,j;m,n}^+E_i-v_i^{c_{ij}+2m}E_if_{i,\tau i,j;m,n}^+&=v_i^{-\frac{1}{2}}(v_i-v_i^{-1})[m+1]_{v_i}f_{i,\tau i,j;m+1,n},\\
			f_{i,\tau i,j;m,n}^+E_{\tau i}-v_i^{c_{\tau i,j}+2n}E_{\tau i}f_{i,\tau i,j;m,n}^+&=v_i^{-\frac{1}{2}}(v_i-v_i^{-1})[n+1]_{v_i}f_{i,\tau i,j;m,n+1},\\
			[F_i,f_{i,\tau i,j;m,n}^+]&=v_i^{\frac{1}{2}}(v_i-v_i^{-1})K_if_{i,\tau i,j;m-1,n}^+,\\
			[F_{\tau i},f_{i,\tau i,j;m,n}^+]&=v_i^{\frac{1}{2}}(v_i-v_i^{-1})K_{\tau i}f_{i,\tau i,j;m,n-1}^+.
		\end{align*}
	\end{lemma}
	
	Similarly the following recursive formulas can be formulated in $\tB^\imath_\tau$, cf. \cite[\S 5.1]{Z23}.
	
	\begin{lemma}\label{lem: f_m1m2 induction relation}
		For $m,n\in\Z$, set $f_{m,n}:=f_{i,\tau i,j;m,n}$ and $f'_{m,n}:=f'_{i,\tau i,j;m,n}$. Then the following relations holds in $\tB^\imath_\tau$:
		\begin{align*}
			&f_{m,n}\ast\vartheta_i-v_i^{c_{ij}+2m}\vartheta_i\ast f_{m,n}\\
			&=v_i^{-\frac{1}{2}}(v_i-v_i^{-1})[m+1]_{v_i}f_{m+1,n}-v_i^{c_{ij}+2m+\frac{1}{2}}(v_i-v_i^{-1})[-c_{\tau i,j}-n+1]_{v_i}h_{\tau i}\ast f_{m,n-1},\\
			&f_{m,n}\ast\vartheta_{\tau i}-v_i^{c_{\tau i,j}+2n}\vartheta_{\tau i}\ast f_{m,n}\\
			&=v_i^{-\frac{1}{2}}(v_i-v_i^{-1})[n+1]_{v_i}f_{m,n+1}-v_i^{c_{\tau i,j}+2n+\frac{1}{2}}(v_i-v_i^{-1})[-c_{ij}-m+1]_{v_i}h_i\ast f_{m-1,n},
		\end{align*}
		\begin{align*}
			&\vartheta_i\ast f'_{m,n}-v_i^{c_{ij}+2m}f'_{m,n}\ast\vartheta_i\\
			&=v_i^{-\frac{1}{2}}(v_i-v_i^{-1})[m+1]_{v_i}f'_{n,m+1,n}-v_i^{c_{ij}+2m+\frac{1}{2}}(v_i-v_i^{-1})[-c_{\tau i,j}-n+1]_{v_i}f'_{m,n-1}\ast h_i,\\
			&\vartheta_{\tau i}\ast f'_{m,n}-v_i^{c_{\tau i,j}+2n}f'_{m,n}\ast\vartheta_{\tau i}\\
			&=v_i^{-\frac{1}{2}}(v_i-v_i^{-1})[n+1]_{v_i}f'_{m,n+1}-v_i^{c_{\tau i,j}+2n+\frac{1}{2}}(v_i-v_i^{-1})[-c_{ij}-m+1]_{v_i}f'_{m-1,n}\ast h_{\tau i}.
		\end{align*}
	\end{lemma}
	
	\begin{proof}
		Using the computations of \cite[Lemma 38.1.7]{Lus93}, one can deduce the following identities:
		\begin{align*}
			\Delta(f_{m,n})&=h_i^{m}h_{\tau i}^{n}h_j\otimes f_{m,n}\\
			&+\sum_{t_1=0}^{m}\sum_{t_2=0}^{n}\prod_{p_1=0}^{m-t_1-1}\prod_{p_2=0}^{n-t_2-1}\frac{(1-v_i^{2(m+c_{ij}-p_1-1)})(1-v_i^{2(n+c_{ij}-p_2-1)})}{(v_i-v_i^{-1})^2}\\
			&\cdot v_i^{(t_1+\frac{1}{2})(m-t_1)+(t_2+\frac{1}{2})(n-t_2)}f_{t_1,t_2}h_i^{m-t_1}h_{\tau i}^{n-t_2}\otimes\vartheta_i^{(m-t_1)}\vartheta_{\tau i}^{(n-t_2)},\\
			\Delta(f'_{m,n})&=f'_{m,n}\otimes 1\\
			&+\sum_{t_1=0}^{m}\sum_{t_2=0}^{n}\prod_{p_1=0}^{m-t_1-1}\prod_{p_2=0}^{n-t_2-1}\frac{(1-v_i^{2(m+c_{ij}-p_1-1)})(1-v_i^{2(n+c_{ij}-p_2-1)})}{(v_i-v_i^{-1})^2}\\
			&\cdot v_i^{(t_1+\frac{1}{2})(m-t_1)+(t_2+\frac{1}{2})(n-t_2)}\vartheta_i^{(m-t_1)}\vartheta_{\tau i}^{(n-t_2)}h_i^{t_1}h_{\tau i}^{t_2}h_j\otimes f'_{t_1,t_2}.
		\end{align*}
		From these the lemma follows by a direct computation.
	\end{proof}
	
	By comparing Lemma~\ref{lem: f_m1m2 induction relation} with the recursive relations in \cite[Section 5.1]{Z23}, we obtain the following explicit expressions for $f_{i,\tau i,j;m,n}$ and $f'_{i,\tau i,j;m,n}$, cf. \cite[Proposition 5.11]{Z23}.
	
	\begin{proposition}\label{prop:f_m1m2 expression}
		For $m,n\in\Z$, in $\tB^\imath_\tau$ we have \begin{equation}\label{eq: f_m1m2 expansion}
			\begin{aligned}
				f_{i,\tau i,j;m,n}&=\sum_{u=0}^{\min\{m,n\}}\sum_{\substack{r_1+s_1=m-u\\ r_2+s_2=n-u}}(-1)^{r_1+r_2}v_i^{\frac{1}{2}(m+n-2u)}(v_i-v_i^{-1})^{-(m+n)+2u}\times\\
				&\qquad v_i^{r_1(c_{ij}+m-1)+r_2(c_{\tau i,j}+n-1)+u(r_1-r_2+c_{ij}+m)}\qbinom{-c_{\tau i,j}-n+u}{u}_{v_i}\times\\
				&\qquad h_{\tau i}^u\ast\vartheta_i^{(r_1)}\ast\vartheta_{\tau i}^{(r_2)}\ast\vartheta_j\ast\vartheta_{\tau i}^{(s_2)}\ast\vartheta_i^{(s_1)},
			\end{aligned}
		\end{equation}
		and
		\begin{equation}\label{eq: f'_m1m2 expansion}
			\begin{aligned}
				f'_{i,\tau i,j;m,n}&=\sum_{u=0}^{\min\{m,n\}}\sum_{\substack{r_1+s_1=m-u\\ r_2+s_2=n-u}}(-1)^{r_1+r_2}v_i^{\frac{1}{2}(m+n-2u)}(v_i-v_i^{-1})^{-(m+n)+2u}\times\\
				&\qquad v_i^{r_1(c_{ij}+m-1)+r_2(c_{\tau i,j}+n-1)+u(r_1-r_2+c_{ij}+m)}\qbinom{-c_{\tau i,j}-n+u}{u}_{v_i}\times\\
				&\qquad\vartheta_i^{(s_1)}\ast\vartheta_{\tau i}^{(s_2)}\ast \vartheta_j\ast\vartheta_{\tau i}^{(r_2)}\ast\vartheta_{i}^{(r_1)}\ast h_i^u.
			\end{aligned}
		\end{equation}
		In particular, $\widetilde{\TT}_i(\vartheta_j)=f'_{i,\tau i,j;-c_{ij},-c_{\tau i,j}}$.
	\end{proposition}
	
	\begin{proof}
		Let us write $f_{m,n}:=f_{i,\tau i,\tau i,j;m,n}$, $f'_{m,n}:=f'_{i,\tau i,\tau i,j;m,n}$. Then 
		\[\vartheta_{\tau i}\ast f'_{0,n}-v_i^{c_{\tau i,j}+2n}f'_{0,n}\ast\vartheta_{\tau i}=v_i^{-\frac{1}{2}}(v_i-v_i^{-1})[n+1]_{v_i}f'_{0,n+1}.\]
		From this one can deduce 
		\[f'_{0,n}=v_i^{\frac{1}{2}n}(v_i-v_i^{-1})^{-n}\sum_{r+s=n}(-1)^rv_i^{r(c_{\tau i,j}+n-1)}\vartheta_{\tau i}^{(r)}\ast\vartheta_j^{(n)}\ast\vartheta_{\tau i}^{(s)}.\]
		Now $f'_{m,n}$ satisfies the following recursive relation:
		\begin{equation}\label{eq: f'_m1m2 induction}
			\begin{aligned}
				&\vartheta_i\ast f'_{m,n}-v_i^{c_{ij}+2m}f'_{m,n}\ast\vartheta_i\\
				&=v_i^{-\frac{1}{2}}(v_i-v_i^{-1})[m+1]_{v_i}f'_{m+1,n}-v_i^{c_{ij}+2m+\frac{1}{2}}(v_i-v_i^{-1})[-c_{\tau i,j}-n+1]_{v_i}f'_{m,n-1}\ast h_i,
			\end{aligned}
		\end{equation}
		and \eqref{eq: f'_m1m2 expansion} is equivalent to
		\begin{equation}\label{eq: f'_m1m2 induction on x_0m-2}
			\begin{aligned}
				f'_{m,n}&=\sum_{u=0}^{m}\sum_{r+s=m-u}(-1)^{r}v_i^{r(c_{ij}+m-1+u)+u(c_{ij}+m)}v_i^{\frac{1}{2}(m-u)}(v_i-v_i^{-1})^{-m+u}\times\\
				&\qbinom{-c_{\tau i,j}-n+u}{u}_i \vartheta_i^{(s)}\ast f'_{0,n-u}\ast\vartheta_{i}^{(r)}\ast h_i^u.
			\end{aligned}
		\end{equation}
		It then suffices to show that \eqref{eq: f'_m1m2 induction on x_0m-2} satisfies the recursive relation \eqref{eq: f'_m1m2 induction}, which can be checked similar to \cite[Proposition 5.11]{Z23}. The formula for $f_{m,n}$ can be obtained similarly.
		
		Finally, setting $m=-c_{ij}$ and $n=-c_{\tau i,j}$ in \eqref{eq: f'_m1m2 expansion} gives
		\begin{align*}
			f'_{i,\tau i,j;-c_{ij},-c_{\tau i,j}}&=\sum_{u=0}^{-\max\{c_{ij},c_{\tau i,j}\}}\sum_{\substack{r_1+s_1=-c_{ij}-u\\ r_2+s_2=-c_{\tau i,j}-u}}(-1)^{r_1+r_2}v_i^{\frac{1}{2}(-c_{ij}-c_{\tau i,j}-2u)}(v_i-v_i^{-1})^{c_{i,j}+c_{\tau i,j}+2u}\times\\
			&\qquad v_i^{-r_1-r_2+u(r_1-r_2)}\vartheta_i^{(s_1)}\ast\vartheta_{\tau i}^{(s_2)}\ast \vartheta_j\ast\vartheta_{\tau i}^{(r_2)}\ast\vartheta_{i}^{(r_1)}\ast h_i^u.
		\end{align*}
		Under the isomorphism $\widetilde{\Phi}^\imath:\tB^\imath_\tau\to \tUi$, the right-hand side is mapped to $\tTT_i(B_j)$ (see \eqref{eq: T_i(B_j) for c_itaui=0}), so the last statement follows.
	\end{proof}
	
	Using the recursive formulas of Lemma~\ref{lem: f_m1m2 induction relation}, we are now ready to prove the following result for the action of $\tTT_i$ on root vectors, cf. \cite[Theorem 7.3]{Z23}.
	
	\begin{proposition}\label{prop:ibraid:diag}
		For any $j\neq i,\tau i$ and $m,n\in\Z$, we have 
		\[\tTT_i(f_{i,\tau i,j;m,n})=f'_{i,\tau i,j;-c_{ij}-m,-c_{\tau i,j}-n}.\]
	\end{proposition}
	\begin{proof}
		For simplicity we write
		\[f_{m,n}:=f_{i,\tau i,j;m,n},\quad f'_{m,n}:=f'_{i,\tau i,j;m,n}.\] 
		The claim is known for $m=n=0$, so we proceed by induction on $(m,n)$. Assume that $\tTT_i(f_{t_1,t_2})=f'_{-c_{ij}-t_1,-c_{\tau i,j}-t_2}$ for any $0\leq t_1\leq m,0\leq t_2\leq n$. Recall that we have
		\begin{equation}\label{eq:x(m,n) and theta_i commutator}
			f_{m,n}\vartheta_i-v_i^{c_{ij}+2m}\vartheta_if_{m,n}=v_i^{-\frac{1}{2}}(v_i-v_i^{-1})[m+1]_{v_i}f_{m+1,n}.
		\end{equation}
		Apply $\tTT_i$ to \eqref{eq:x(m,n) and theta_i commutator}, we get 
		\begin{equation}\label{eq:T_ix(m,n) induction}
			v_i^{-\frac{1}{2}}(v_i-v_i^{-1})[m+1]_{v_i}\tTT_i(f_{m+1,n})=\tTT_i(f_{m+1,n}\vartheta_i)-v_i^{c_{ij}+2m}\tTT_i(\vartheta_i f_{m,n}).
		\end{equation}
		From the recursive relations in Lemma~\ref{lem: f_m1m2 induction relation}, one can deduce
		\begin{align*}
			\tTT_i(f_{m,n}\vartheta_i)&=\tTT_i(f_{m,n})\ast\tTT_i(\vartheta_i)=v_if'_{-c_{ij}-m,-c_{\tau i,j}-n}\ast\K_i^{-1}\ast\vartheta_{\tau i}\\
			&=v_i^{(c_{ij}-c_{\tau i,j})+2(m-n)+1}\K_i^{-1}\ast(f'_{-c_{ij}-m,-c_{\tau i,j}-n}\vartheta_{\tau i})\\
			&\quad+v_i^{-\frac{1}{2}}(v_i-v_i^{-1})[m+1]_{v_i}f'_{-c_{ij}-m-1,-c_{\tau i,j}-n},
		\end{align*}
		and
		\begin{align*}
			\tTT_i(\vartheta_if_{m,n})&=\tTT_i(\vartheta_i\ast f_{m,n})-v_i^{\frac{1}{2}}(v_i-v_i^{-1})[-c_{\tau i,j}-n+1]_{v_i}\K_i\ast f_{m,n-1}\\
			&=v_i\K_i^{-1}\ast(\vartheta_{\tau i}f'_{-c_{ij}-m,-c_{\tau i,j}-n})\\
			&\quad-v_i^{\frac{1}{2}}(v_i-v_i^{-1})[-c_{\tau i,j}-n+1]_{v_i}\K_i^{-1}\ast f'_{-c_{ij}-m,-c_{\tau i,j}-n+1}.
		\end{align*}
		Now \eqref{eq:T_ix(m,n) induction} can be computed as 
		\begin{align*}
			&v_i^{-\frac{1}{2}}(v_i-v_i^{-1})[m+1]_{v_i}\tTT_i(f_{m+1,n})\\
			&=v_i^{(c_{ij}-c_{\tau i,j})+2(m-n)+1}\K_i^{-1}\ast(f'_{-c_{ij}-m,-c_{\tau i,j}-n}\vartheta_{\tau i}-v_i^{c_{\tau i,j}+2n}\vartheta_{\tau i}f'_{-c_{ij}-m,-c_{\tau i,j}-n})\\
			&\quad+v_i^{-\frac{1}{2}}(v_i-v_i^{-1})[m+1]_{v_i}f'_{-c_{ij}-m-1,-c_{\tau i,j}-n}\\
			&\quad+v_i^{c_{ij}+2m+\frac{1}{2}}(v_i-v_i^{-1})[-c_{\tau i,j}-n+1]_{v_i}\K_i^{-1}\ast f'_{-c_{ij}-m,-c_{\tau i,j}-n+1}.
		\end{align*}
		Note that we also have
		\begin{align*}
			&f'_{-c_{ij}-m,-c_{\tau i,j}-n}\vartheta_{\tau i}-v_i^{c_{\tau i,j}+2n}\vartheta_{\tau i}f'_{-c_{ij}-m,-c_{\tau i,j}-n}\\
			&=-v_i^{c_{\tau i,j}+2n-\frac{1}{2}}(v_i-v_i^{-1})[-c_{\tau i,j}-n+1]_{v_i}f'_{-c_{ij}-m,-c_{\tau i,j}-n+1}.
		\end{align*}
		From this we conclude
		\[v_i^{-\frac{1}{2}}(v_i-v_i^{-1})[m+1]_{v_i}\tTT_i(f_{m+1,n})=v_i^{-\frac{1}{2}}(v_i-v_i^{-1})[m+1]_{v_i}f'_{-c_{ij}-m-1,-c_{\tau i,j}-n}.\]
		Since $v_i^{-\frac{1}{2}}(v_i-v_i^{-1})[m+1]_{v_i}\neq 0$, this implies $\tTT_i(f_{m+1,n})=f'_{-c_{ij}-m-1,-c_{\tau i,j}-n}$. Using the same idea we can show that $\tTT_i(f_{m,n+1})=f'_{-c_{ij}-m,-c_{\tau i,j}-n-1}$. This completes the induction step.
	\end{proof}

	\subsection{The case for $c_{i,\tau i}=-1$}
	\label{subsec: root vector c_i,taui=-1}
	
	For $i\neq\tau i$ and $c_{i,\tau i}=-1$, we define the root vectors: 
	\begin{equation}\label{eq:root vector c_i,taui=-1-1}
		\begin{aligned}
			f_{i,\tau i,j;a,b,c}&=\sum_{\substack{r_1+s_1=a\\ r_2+s_2=b\\ r_3+s_3=c}}(-1)^{r_1+r_2+r_3}v_i^{\frac{1}{2}(a+b+c)}(v_i-v_i^{-1})^{-(a+b+c)}\\
			&\times v_i^{r_1(c_{ij}+a-1)+r_2(c_{\tau i,j}+b-1)+r_3(c_{ij}+c-1)+r_1(2c-b)-r_2c}\vartheta_i^{(r_1)}\vartheta_{\tau i}^{(r_2)}\vartheta_i^{(r_3)}\vartheta_j\vartheta_i^{(s_3)}\vartheta_{\tau i}^{(s_2)}\vartheta_i^{(s_1)};
		\end{aligned}
	\end{equation}
	\begin{equation}\label{eq:root vector c_i,taui=-1-2}
		\begin{aligned}
			f'_{i,\tau i,j;a,b,c}&=\sum_{\substack{r_1+s_1=a\\ r_2+s_2=b\\ r_3+s_3=c}}(-1)^{r_1+r_2+r_3}v_i^{\frac{1}{2}(a+b+c)}(v_i-v_i^{-1})^{-(a+b+c)}\\
			&\times v_i^{r_1(c_{ij}+a-1)+r_2(c_{\tau i,j}+b-1)+r_3(c_{ij}+c-1)+r_1(2c-b)-r_2c}\vartheta_i^{(s_1)}\vartheta_{\tau i}^{(s_2)}\vartheta_i^{(s_3)}\vartheta_j\vartheta_i^{(r_3)}\vartheta_{\tau i}^{(r_2)}\vartheta_i^{(r_1)}.
		\end{aligned}
	\end{equation}
	
	Comparing this with the definition of $x_{i,\tau i,j;a,b,c}$ and $x'_{i,\tau i,j;a,b,c}$ in \cite[Definition 6.1]{Z23}, we obtain the following lemma.
	
	\begin{lemma}[{\cite[Lemma 6.5]{Z23}}]\label{lem: f_abc relation with theta}
		We have, for $a,b,c\geq 0$,
		\begin{align*}
			&[f_{i,\tau i,j;a,b,c}^+,E_i]_{v_i^{2(a+c)-b+c_{ij}}}=v_i^{-\frac{1}{2}}(v_i-v_i^{-1})[a+1]_{v_i}f_{i,\tau i,j;a+1,b,c}^+,\\
			&[f_{i,\tau i,j;a,b,c}^+,E_{\tau i}]_{v_i^{2b-(a+c)+c_{\tau i,j}}}=v_i^{-\frac{1}{2}}(v_i-v_i^{-1})([b-a+1]_{v_i}f_{i,\tau i,j;a,b+1,c}^+\\
			&\hspace{5.1cm}+[c+1]_{v_i}f_{i,\tau i,j;a-1,b+1,c+1}^+),\\
			&[F_i,f_{i,\tau i,j;a,b,c}^+]=v_i^{\frac{1}{2}}(v_i-v_i^{-1})([-c_{ij}-a+b-2c+1]_{v_i}K_if_{i,\tau i,j;a-1,b,c}^+\\
			&\hspace{2.9cm}+[-c_{ij}-c+1]_{v_i}K_if_{i,\tau i,j;a,b,c-1}^+),\\
			&[F_{\tau i},f_{i,\tau i,j;a,b,c}^+]=v_i^{\frac{1}{2}}(v_i-v_i^{-1})[-c_{\tau i,j}-b+c+1]_{v_i}K_{\tau i}f_{i,\tau i,j;a,b-1,c}^+.
		\end{align*}
	\end{lemma}
	
	The following lemma generalizes the result of \cite[Lemma 6.13]{Z23}, and can be viewed as a counterpart of \eqref{eq: tU T_i action on f_m}.
	
	\begin{lemma}\label{eq: tU T_ri action on f_abc}
		For $a,b,c\in\Z$, the following holds in $\tU$:
		\begin{equation}\label{eq: T_bri of f_abc}
			\widetilde{T}_{r_i}(f_{i,\tau i,j;a,b,c}^+)=f'^+_{i,\tau i,j;-c_{\tau i,j}-b+c,-c_{\tau i,j}-c_{ij}-a-c,-c_{ij}-c}.
		\end{equation}
	\end{lemma}
	\begin{proof}
		For simplicity we write $f_{a,b,c}:=f_{i,\tau i,j;a,b,c}$ and $f'_{a,b,c}:=f'_{i,\tau i,j;a,b,c}$. We first prove that
		\begin{equation}\label{eq: T_bri of f_a00}
			\widetilde{T}_{r_i}(f_{a,0,0}^+)=f'^+_{-c_{\tau i,j},-c_{\tau i,j}-c_{ij}-a,-c_{ij}}.
		\end{equation}
		Note that from \cite[Lemma 6.13]{Z23} we know
		\[\widetilde{T}_{r_i}(f_{0,0,0}^+)=f'^+_{-c_{\tau i,j},-c_{\tau i,j}-c_{ij},-c_{ij}}.\]
		Now assume that \eqref{eq: T_bri of f_a00} holds for some $a\geq 0$. Then applying $\widetilde{T}_{r_i}$ to
		\[f_{a,0,0}^+E_i-v_i^{2a+c_{ij}}E_if_{a,0,0}^+=v_i^{-\frac{1}{2}}(v_i-v_i^{-1})[a+1]_{v_i}f_{a+1,0,0}^+,\]
		we find that 
		\begin{align*}
			&v_i^{-\frac{1}{2}}(v_i-v_i^{-1})[a+1]_{v_i}\widetilde{T}_{r_i}(f_{a+1,0,0}^+)\\
			&=v_i^{-1}f'^+_{-c_{\tau i,j},-c_{\tau i,j}-c_{ij}-a,-c_{ij}}F_{\tau i}K_{\tau i}-v_i^{2a+c_{ij}-1}F_{\tau i}K_{\tau i}f'^+_{-c_{\tau i,j},-n(c_{\tau i,j}+c_{ij})-a,-c_{ij}}\\
			&=v_i^{-1}[f'_{-c_{\tau i,j},-c_{\tau i,j}-c_{ij}-a,-c_{ij}},F_{\tau i}]K_{\tau i}\\
			&=v_i^{-\frac{1}{2}}(v_i-v_i^{-1})[a+1]_{v_i}f'^+_{-c_{\tau i,j},-c_{\tau i,j}-c_{ij}-a-1,-c_{ij}}.
		\end{align*}
		Since $v_i^{-\frac{1}{2}}(v_i-v_i^{-1})[a+1]_{v_i}\neq 0$, this implies $\widetilde{T}_{r_i}(f_{a+1,0,0}^+)=f'^+_{-c_{\tau i,j},-c_{\tau i,j}-c_{ij}-a-1,-c_{ij}}$. 
		
		We now prove the general case by induction on $(b,c)$. Assume that \eqref{eq: T_bri of f_abc} is true for any $f_{a,s,t}$ with $s\leq b,t\leq c$ and $(s,t)\neq (b,c)$. Now applying $\widetilde{T}_{r_i}$ to the identity
		\begin{align*}
			[f_{a+1,b-1,c-1}^+,E_{\tau i}]_{v_i^{2b-2-(a+c)+c_{\tau i,j}}}=v_i^{-\frac{1}{2}}(v_i-v_i^{-1})([b-a-1]_{v_i}f_{a+1,b,c-1}^++[c]_{v_i}f_{a,b,c}^+),
		\end{align*}
		we get
		\begin{align*}
			&v_i^{-\frac{1}{2}}(v_i-v_i^{-1})[c]_{v_i}\widetilde{T}_{r_i}(f_{a,b,c}^+)
			\\
			&=v_i^{-1}[f'^+_{-c_{\tau i,j}-b+c,-c_{\tau i,j}-c_{ij}-a-c,-c_{ij}-c+1},F_i]K_i\\
			&\quad-v_i^{-\frac{1}{2}}(v_i-v_i^{-1})[b-a-1]_{v_i}f'^+_{-c_{\tau i,j}-b+c-1,-c_{\tau i,j}-c_{ij}-a-c,-c_{ij}-c+1}\\
			&=v_i^{-\frac{1}{2}}(v_i-v_i^{-1})[c]_{v_i}f'^+_{-c_{\tau i,j}-b+c,-c_{\tau i,j}-c_{ij}-a-c,-c_{ij}-c}.
		\end{align*}
		Since $v_i^{-\frac{1}{2}}(v_i-v_i^{-1})[c]_{v_i}\neq 0$, we conclude 
		\[\widetilde{T}_{r_i}(f_{a,b,c}^+)=f'^+_{-c_{\tau i,j}-b+c,-c_{\tau i,j}-c_{ij}-a-c,-c_{ij}-c}.\]
		This completes the induction step.
	\end{proof}
	
	\begin{lemma}\label{lem: f_abc induction relation}
		For $a,b,c\in\Z$, write $f_{a,b,c}:=f_{i,\tau i,j;a,b,c}$, $f'_{a,b,c}:=f'_{i,\tau i,j;a,b,c}$. Then
		\begin{equation}\label{eq: f_abc induction with theta_i}
			\begin{aligned}
				&f_{a,b,c}\ast\vartheta_i-v_i^{c_{ij}+2(a+c)-b}\vartheta_i\ast f_{a,b,c}=v_i^{-\frac{1}{2}}(v_i-v_i^{-1})[a+1]_{v_i}f_{a+1,b,c}\\
				&\qquad-v_i^{c_{ij}+2(a+c)-b+\frac{1}{2}}(v_i-v_i^{-1})[-c_{\tau i,j}-b+c+1]_{v_i}h_{\tau i}\ast f_{a,b-1,c},
			\end{aligned}
		\end{equation}
		\begin{equation}\label{eq: f_abc induction with theta_taui}
			\begin{aligned}
				&f_{a,b,c}\ast \vartheta_{\tau i}-v_i^{c_{\tau i,j}+2b-(a+c)}\vartheta_{\tau i}\ast f_{a,b,c}=v_i^{-\frac{1}{2}}(v_i-v_i^{-1})([b-a+1]_{v_i}f_{a,b+1,c}
				\\
				&\qquad\qquad+[c+1]_{v_i}f_{a-1,b+1,c+1})-v_i^{c_{\tau i,j}+2b-(a+c)+\frac{1}{2}}(v_i-v_i^{-1})h_i\\
				&\qquad\qquad\ast([-c_{ij}-a-2c+b+1]_{v_i}f'_{a-1,b,c}+[-c_{ij}-c+1]_{v_i}f_{a,b,c-1}),
			\end{aligned}
		\end{equation}
		\begin{equation}\label{eq: f'_abc induction with theta_i}
			\begin{aligned}
				&\vartheta_i\ast f'_{a,b,c}-v_i^{c_{ij}+2(a+c)-b}f'_{a,b,c}\ast\vartheta_i=v_i^{-\frac{1}{2}}(v_i-v_i^{-1})[a+1]_{v_i}f'_{a+1,b,c}
				\\
				&\qquad\qquad-v_i^{c_{ij}+2(a+c)-b+\frac{1}{2}}(v_i-v_i^{-1})[-c_{\tau i,j}-b+c+1]_{v_i}f'_{a,b-1,c}\ast h_i,
			\end{aligned}
		\end{equation}
		\begin{equation}\label{eq: f'_abc induction with theta_taui}
			\begin{aligned}
				&\vartheta_{\tau i}\ast f'_{a,b,c}-v_i^{c_{\tau i,j}+2b-(a+c)}f'_{a,b,c}\ast\vartheta_i=v_i^{-\frac{1}{2}}(v_i-v_i^{-1})([b-a+1]_{v_i}f'_{a,b+1,c}
				\\
				&\qquad\qquad+[c+1]_{v_i}f'_{a-1,b+1,c+1})-v_i^{c_{\tau i,j}+2b-(a+c)+\frac{1}{2}}(v_i-v_i^{-1})
				\\
				&\qquad\qquad([-c_{ij}-a-2c+b+1]_{v_i}f'_{a-1,b,c}+[-c_{ij}-c+1]_{v_i}f'_{a,b,c-1})\ast h_{\tau i}.
			\end{aligned}
		\end{equation}
	\end{lemma}
	\begin{proof}
		Using the definitions $f_{a,b,c}$ and $f'_{a,b,c}$, one can compute that 
		\begin{align*}
			\Delta(f_{a,b,c})&=f_{a,b,c}\otimes 1+v_i^{c_{ij}+2(a+c-1)-b+\frac{1}{2}}[-c_{ij}-a-2c+b+1]_{v_i}f'_{a-1,b,c}h_i\otimes\vartheta_i\\
			&+v_i^{c_{\tau i,j}+2(b-1)-(a+c)+\frac{1}{2}}[-c_{\tau i,j}-b+c+1]_{v_i}f'_{a,b-1,c}h_{\tau i}\otimes\vartheta_{\tau i}\\
			&+v_i^{c_{ij}+2(a+c-1)-b}[-c_{ij}-c+1]_{v_i}f'_{a,b,c-1}h_i\otimes\vartheta_i+\cdots+h_i^{a+c}h_{\tau i}^bh_j^n\otimes f'_{a,b,c},\\
			\Delta(f'_{a,b,c})&=f'_{a,b,c}\otimes 1+v_i^{c_{ij}+2(a+c-1)-b+\frac{1}{2}}[-c_{ij}-a-2c+b+1]_{v_i}\vartheta_ih_i^{a+c-1}h_{\tau i}^bh_j^n\otimes f'_{a-1,b,c}\\
			&+v_i^{c_{\tau i,j}+2(b-1)-(a+c)+\frac{1}{2}}[-c_{\tau i,j}-b+c+1]_{v_i}\vartheta_{\tau i}h_i^{a+c}h_{\tau i}^{b-1}h_j^n\otimes f'_{a,b-1,c}\\
			&+v_i^{c_{ij}+2(a+c-1)-b}[-c_{ij}-c+1]_{v_i}\vartheta_ih_i^{a+c-1}h_{\tau i}^bh_j^n\otimes f'_{a,b,c-1}+\cdots+h_i^{a+c}h_{\tau i}^bh_j^n\otimes f'_{a,b,c}.
		\end{align*}
		From these the assertion follows from Lemma~\ref{lem: f_abc relation with theta}.
	\end{proof}
	
	By comparing Lemma \ref{lem: f_abc induction relation} with the recursive relations in \cite[\S 6.1]{Z23}, we get the following explicit expression for $f_{i,\tau i,j;a,b,c}$ and $f'_{i,\tau i,j;a,b,c}$.
	
	\begin{proposition}\label{prop: f_abc expression}
		For $a,b,c\in\Z$, in $\tB^\imath_\tau$
		\begin{align*}
			f_{i,\tau i,j;a,b,c}&=\sum_{w=0}^{\min\{a,b\}}\sum_{u=0}^{\min\{b-w,c\}}\sum_{\substack{r_1+s_1=a-w\\ r_2+s_2=b-u-w\\ r_3+s_3=c-u}}(-1)^{r_1+r_2+r_3}v_i^{\frac{1}{2}(a+b+c-2u-2w)}(v_i-v_i^{-1})^{-(a+b+c)+2u+2w}\times\\
			&v_i^{r_1(c_{ij}-b+2c+a-1+2w)+r_2(c_{\tau i,j}+b-w-c+2u-1)+r_3(c_{ij}+c-u-1)+w(c_{ij}-b+2c+a)+u(c_{\tau i,j}+b-w-c)}\\
			&\qbinom{-c_{\tau i,j}-b+c+w}{w}_{v_i}\qbinom{-c_{ij}-c+u}{u}_{v_i}\times\\
			&h_i^w\ast\vartheta_i^{(r_1)}\ast h_{\tau i}^u\ast\vartheta_{\tau i}^{(r_2)}\ast\vartheta_i^{(r_3)}\ast\vartheta_j\ast\vartheta_i^{(s_3)}\ast\vartheta_{\tau i}^{(s_2)}\ast\vartheta_i^{(s_1)},
		\end{align*}
		and
		\begin{align*}
			f'_{i,\tau i,j;a,b,c}&=\sum_{w=0}^{\min\{a,b\}}\sum_{u=0}^{\min\{b-w,c\}}\sum_{\substack{r_1+s_1=a-w\\ r_2+s_2=b-u-w\\ r_3+s_3=c-u}}(-1)^{r_1+r_2+r_3}v_i^{\frac{1}{2}(a+b+c-2u-2w)}(v_i-v_i^{-1})^{-(a+b+c)+2u+2w}\times\\
			&v_i^{r_1(c_{ij}-b+2c+a-1+2w)+r_2(c_{\tau i,j}+b-w-c+2u-1)+r_3(c_{ij}+c-u-1)+w(c_{ij}-b+2c+a)+u(c_{\tau i,j}+b-w-c)}\\
			&\qbinom{-c_{\tau i,j}-b+c+w}{w}_{v_i}\qbinom{-c_{ij}-c+u}{u}_{v_i}\times\\
			&\vartheta_i^{(s_1)}\ast\vartheta_{\tau i}^{(s_2)}\ast\vartheta_i^{(s_3)}\ast\vartheta_j\ast\vartheta_i^{(r_3)}\ast\vartheta_{\tau i}^{(r_2)}\ast h_{\tau i}^u\ast\vartheta_i^{(r_1)}\ast h_i^{w}.
		\end{align*}
		In particular, we have $\widetilde{\TT}_i(\vartheta_j)=f'_{i,\tau i,j;-c_{\tau i,j},-c_{ij}-c_{\tau i,j},-c_{ij}}$.
	\end{proposition}
	
	\begin{proof}
		For simplicity we write $f_{a,b,c}:=f_{i,\tau i,j;a,b,c}$ and $f'_{a,b,c}:=f'_{i,\tau i,j;a,b,c}$. Setting $a=b=0$ in \eqref{eq: f'_abc induction with theta_i} gives
		\[\vartheta_i\ast f_{0,0,c}-v_i^{2c+c_{ij}}=v_i^{-\frac{1}{2}}(v_i-v_i^{-1})[c+1]_{v_i}f_{0,0,c+1},\]
		and from this one can deduce 
		\[f'_{0,0,c}=v_i^{\frac{1}{2}c}(v_i-v_i^{-1})^{-c}\sum_{r+s=c}(-1)^rv_i^{r(c_{ij}+c-1)}\vartheta_i^{(r)}\ast\vartheta_j\ast\vartheta_i^{(s)}.\]
		Now set $a=0$ in \eqref{eq: f'_abc induction with theta_taui}, we get the following recursive relation for $f'_{0,b,c}$:
		\begin{align*}
			&\vartheta_{\tau i}\ast f'_{0,b,c}-v_i^{c_{\tau i,j}+2b-c}f'_{0,b,c}\ast\vartheta_{\tau i}\\
			&=v_i^{-\frac{1}{2}}(v_i-v_i^{-1})[b+1]_{v_i}f'_{0,b+1,c}-v_i^{c_{\tau i,j}+2b-c+\frac{1}{2}}(v_i-v_i^{-1})[-c_{ij}-c+1]_{v_i}f'_{0,b,c-1}\ast h_{\tau i},
		\end{align*}
		which implies
		\begin{align*}
			f'_{0,b,c}&=\sum_{u=0}^{\min\{b,c\}}\sum_{r+s=b-u}(-1)^{r}v_i^{r(c_{\tau i,j}+b-c-1+2u)+u(c_{\tau i,j}+b-c)}\times\\
			&\qquad v_i^{\frac{1}{2}(b-u)}(v_i-v_i^{-1})^{-b+u}\qbinom{-c_{ij}-c+u}{u}_{v_i} \vartheta_{\tau i}^{(s)}\ast f'_{0,0,c-u}\ast\vartheta_{\tau i}^{(r)}\ast h_{\tau i}^u.
		\end{align*}
		Finally, using \eqref{eq: f'_abc induction with theta_i} one can prove that 
		\begin{align*}
			f'_{a,b,c}&=\sum_{w=0}^{-c_{\tau i,j}}\sum_{u=0}^{-c_{ij}}\sum_{\substack{r_1+s_1=-c_{\tau i,j}-w\\ r_2+s_2=-c_{\tau i,j}-c_{ij}-u-w\\ r_3+s_3=-c_{ij}-u}}(-1)^{r_1+r_2+r_3}v_i^{\frac{1}{2}(-2c_{ij}-2c_{\tau i,j}-2u-2w)}(v_i-v_i^{-1})^{-(a+b+c)+2u+2w}\times\\
			&v_i^{r_1(c_{ij}-b+2c+a-1+2w)+r_2(c_{\tau i,j}+b-w-c+2u-1)+r_3(c_{ij}+c-u-1)+w(c_{ij}-b+2c+a)+u(c_{\tau i,j}+b-w-c)}\times\\
			&\qbinom{-c_{\tau i,j}-b+c+w}{w}_{v_i}\qbinom{-c_{ij}-c+u}{u}_{v_i}\times\\
			&\vartheta_i^{(s_1)}\ast\vartheta_{\tau i}^{(s_2)}\ast\vartheta_i^{(s_3)}\ast\vartheta_j\ast\vartheta_i^{(r_3)}\ast\vartheta_{\tau i}^{(r_2)}\ast h_{\tau i}^u\ast\vartheta_i^{(r_1)}\ast h_i^{w}.
		\end{align*}
		Combining all of these gives the expression for $f'_{a,b,c}$. The expression of $f_{a,b,c}$ can be proved similarly.
		
		Finally, setting $a=-c_{\tau i,j}$, $b=-c_{\tau i,j}-c_{ij}$ and $c=-c_{ij}$ in \eqref{eq: f'_abc induction with theta_taui} gives
		\begin{align*}
			&f'_{i,\tau i,j;-c_{\tau i,j},-c_{ij}-c_{\tau i,j},-c_{ij}}\\
			&=\sum_{w=0}^{-c_{\tau i,j}}\sum_{u=0}^{-c_{ij}}\sum_{\substack{r_1+s_1=-c_{\tau i,j}-w\\ r_2+s_2=-c_{ij}-c_{\tau i,j}-u-w\\ r_3+s_3=-c_{ij}-u}}(-1)^{r_1+r_2+r_3}v_i^{-c_{ij}-c_{\tau i,j}-u-w}(v_i-v_i^{-1})^{2(c_{ij}+c_{\tau i,j}+u+w)}\\
			&\times v_i^{-(r_1+r_2+r_3)+2wr_1+(2u-w)r_2-ur_3-uw}\vartheta_i^{(s_1)}\ast\vartheta_{\tau i}^{(s_2)}\ast\vartheta_i^{(s_3)}\ast\vartheta_j\ast\vartheta_i^{(r_3)}\ast\vartheta_{\tau i}^{(r_2)}\ast h_{\tau i}^u\ast\vartheta_i^{(r_1)}\ast h_i^{w}.
		\end{align*}
		Under the isomorphism $\widetilde{\Phi}^\imath:\tB^\imath_\tau\to \tUi$, the right-hand side is mapped to $\tTT_i(B_j)$ (see \eqref{eq: T_i(B_j) for c_itaui=-1}), so the last statement follows.
	\end{proof}
	
	\begin{proposition}\label{prop:ibraid:quasisplit}
		We have, for $a,b,c\geq 0$,
		\begin{align}
			\tTT_i(f_{i,\tau i,j;a,b,c})=f'_{i,\tau i,j;-c_{\tau i,j}-b+c,-c_{\tau i,j}-c_{ij}-a-c,-c_{ij}-c}.
		\end{align}
	\end{proposition}
	\begin{proof}
		We write $f_{a,b,c}:=f_{i,\tau i,j;1,a,b,c}$, $f'_{a,b,c}:=f'_{i,\tau i,j;1,a,b,c}$. The case $a=b=c=0$ has already been proved. Now we use induction on $a$ to show that 
		\begin{align}\label{eq: T_i of f_a00}
			\tTT_i(f_{a,0,0})=f'_{-c_{\tau i,j},-c_{\tau i,j}-c_{ij}-a,-c_{ij}}.
		\end{align}
		For this, assume that \eqref{eq: T_i of f_a00} is true for some $a\in\N$. Applying $\tTT_i$ to the identity
		\begin{align}
			f_{a,0,0}\ast\vartheta_i-v_i^{2a+c_{ij}}\vartheta_i\ast f_{a,0,0}=v_i^{-\frac{1}{2}}(v_i-v_i^{-1})[a+1]_{v_i}f_{a+1,0,0},
		\end{align}
		we get
		\begin{align*}
			&v_i^{-\frac{1}{2}}(v_i-v_i^{-1})[a+1]_{v_i}\tTT_i(f_{a+1,0,0})\\
			&=\tTT_i(f_{a,0,0}\ast\vartheta_i-v_i^{2a+c_{ij}}\vartheta_i\ast f_{a,0,0})\\
			&=v^{c_{ij}-c_{\tau i,j}+3a+\frac{3}{2}}\K_{\tau i}^{-1}\ast(f'_{-c_{\tau i,j},-c_{ij}-c_{\tau i,j}-a,-c_{ij}}\ast\vartheta_i-v_i^{c_{\tau i,j}-a}\vartheta_i\ast f'_{-c_{\tau i,j},-c_{ij}-c_{\tau i,j}-a,-c_{ij}})\\
			&=-v^{c_{ij}-c_{\tau i,j}+3a+1}(v_i-v_i^{-1})[-c_{\tau i,j}+1]_{v_i}\K_{\tau i}^{-1}\ast f'_{-c_{\tau i,j}+1,-c_{ij}-c_{\tau i,j}-a,-c_{ij}}\\
			&+v_i^{-\frac{1}{2}}(v_i-v_i^{-1})[a+1]_{v_i}f'_{-c_{\tau i,j},-c_{ij}-c_{\tau i,j}-a-1,-c_{ij}}.
		\end{align*}
		Note that $f'_{-c_{\tau i,j}+1,-c_{ij}-c_{\tau i,j}-a,-c_{ij}}$ vanishes since 
		\[\widetilde{T}_{r_i}^{-1}(f'_{-c_{\tau i,j}+1,-c_{ij}-c_{\tau i,j}-a,-c_{ij},-c_{ij}})=f_{a,-1,0}=0\]
		and this implies $\tTT_i(f_{a+1,0,0})=f'_{-c_{\tau i,j},-c_{ij}-c_{\tau i,j}-a-1,-c_{ij}}$.
		
		We now prove the general case by induction on $(b,c)$. Assume that the claim is true for $f_{a,s,t}$ with $s\leq b,t\leq c$ and $(s,t)\neq(b,c)$. Then applying $\widetilde{T}_{r_i}$ to the identity
		\begin{align*}
			&f_{a+1,b-1,c-1}\ast\vartheta_{\tau i}-v_i^{c_{\tau i,j}+2b-(a+c)-2}f_{a+1,b-1,c-1}\ast\vartheta_{\tau i}\\
			&=v_i^{-\frac{1}{2}}(v_i-v_i^{-1})\big([b-a-1]_{v_i}f_{a+1,b,c-1}+[c]_{v_i}f_{a,b,c}\big) -v_i^{c_{\tau i,j}+2b-(a+c)-1}(v_i-v_i^{-1})\times\\
			&\qquad\K_{\tau i}\ast \big([-c_{ij}-a-2c+b+1]_{v_i}f_{a,b-1,c-1}+[-c_{ij}-c+2]_{v_i}f_{a+1,b-1,c-2}\big)
		\end{align*}
		and using \eqref{eq: f_abc induction with theta_taui}, one finds that
		\begin{align*}
			&v_i^{-\frac{1}{2}}(v_i-v_i^{-1})\big([b-a-1]_{v_i}f'_{-c_{\tau i,j}-b+c-1,-(c_{\tau i,j}+c_{ij})-a-c,-c_{ij}-c+1}+[c]_{v_i}\widetilde{T}_{r_i}(f_{a,b,c})\big)\\
			&=v_i^{-\frac{1}{2}}(v_i-v_i^{-1})\big([b-a+1]_{v_i}f'_{-c_{\tau i,j}-b+c-1,-(c_{\tau i,j}+c_{ij})-a-c,-c_{ij}-c+1}\\
			&\quad+[c]_{v_i}f'_{-c_{\tau i,j}-b+c,-(c_{\tau i,j}+c_{ij})-a-c,-c_{ij}-c}\big).
		\end{align*}
		This completes the induction step.
	\end{proof}


\end{document}